# Box products in nilpotent normal form theory: The factoring method

James Murdock *

*Department of Mathematics, Iowa State University, Ames, IA 50011, United States*

**Abstract**

Let $N$ be a nilpotent matrix and consider vector fields $\dot{\mathbf{x}} = N\mathbf{x} + \mathbf{v}(\mathbf{x})$ in normal form. Then $\mathbf{v}$ is equivariant under the flow $e^{N^*t}$ for the inner product normal form or $e^{Mt}$ for the $\mathfrak{sl}_2$ normal form. These vector equivariants can be found by finding the scalar invariants for the Jordan blocks in $N^*$ or $M$; taking the *box product* of these to obtain the invariants for $N^*$ or $M$ itself; and then *boosting* the invariants to equivariants by another box product. These methods, developed by Murdock and Sanders in 2007, are here given a self-contained exposition with new foundations and new algorithms yielding improved (simpler) Stanley decompositions for the invariants and equivariants. Ideas used include transvectants (from classical invariant theory), Stanley decompositions (from commutative algebra), and integer cones (from integer programming). This approach can be extended to covariants of $\mathfrak{sl}_2^k$ for $k > 1$, known as SLOCC in quantum computing.

## 1. Introduction

Let $N$ be an $n \times n$ real nilpotent matrix, and consider smooth systems of differential equations in $\mathbb{R}^n$ of the form

* Correspondence to: 1510 Roosevelt Ave., Ames, Iowa 50010, United States. Fax: +1 515 294 5454.
*E-mail address:* jmurdock@iastate.edu.

$$\dot{\mathbf{x}} = N\mathbf{x} + \mathbf{v}(\mathbf{x}), \tag{1.1}$$

where $\mathbf{v}$ is **strictly nonlinear** ($\mathbf{v}(\mathbf{0}) = \mathbf{0} \in \mathbb{R}^n$ and $\mathbf{v}'(\mathbf{0}) = 0 \in \mathfrak{gl}_n$). We usually regard (1.1) as a **formal system**, that is, $\mathbf{v}(\mathbf{x})$ is expanded in a formal power series (which will contain only terms of degree $\geq 2$). Choose a **normal form style** for such systems, either the $\mathfrak{sl}_2$ style, the inner product style, or the simplified style (see Section 2). For a given $N$ and a given normal form style, the **weak description problem** calls for a description of the set of all $\mathbf{v}$ such that (1.1) is in normal form; the **strong description problem** calls for a procedure that generates all such systems exactly once, without repetition. For instance, if $N = N_2$ is the $2 \times 2$ nilpotent matrix in upper Jordan form, the solution of both description problems for the inner product normal form style can be given by the expression

$$\mathbf{v} \in \mathbb{R}[[x_1]]x_1 \begin{bmatrix} x_1 \\ x_2 \end{bmatrix} \oplus \mathbb{R}[[x_1]]x_1^2 \begin{bmatrix} 0 \\ 1 \end{bmatrix}, \tag{1.2}$$

called a **Stanley decomposition**. The meaning of (1.2) is that:

1. A system (1.1) with $N = N_2$ is in inner product normal form if and only if it can be written as

$$\mathbf{v}(\mathbf{x}) = f(x_1)x_1 \begin{bmatrix} x_1 \\ x_2 \end{bmatrix} + g(x_1)x_1^2 \begin{bmatrix} 0 \\ 1 \end{bmatrix}, \tag{1.3}$$

   where $f$ and $g$ are formal power series in $x_1$. Therefore (1.2) solves the weak description problem.
2. No two distinct systems of the form (1.3) will be equal (as formal power series). Thus all systems in normal form are generated exactly once as $f$ and $g$ range over $\mathbb{R}[[x_1]]$, and (1.2) also solves the strong description problem.

**1.1. Remark.** For the more familiar semisimple case, where $N$ is replaced by a matrix $S$ that is diagonalizable over $\mathbb{C}$, the weak description problem is solved by requiring that $\mathbf{v}$ contain only resonant vector monomials (in appropriate complex variables). This does not solve the strong description problem, since there are relations among the resonant monomials allowing the same vector field to be written in more than one way. Finding Stanley decompositions in the semisimple case has only been addressed in a few cases. The semisimple Hamiltonian $(1, 1, \ldots, 1)$ resonance is studied in [5], and the $(p, q, r)$ Hamiltonian resonance is treated in [29, Thm. 10.7.1]. The $(p, q)$ resonance for two coupled nonlinear oscillators, not necessarily Hamiltonian, is studied in [22, §4.5].

Cushman and Sanders [10] solve the strong description problem for nilpotent $N$ in the $\mathfrak{sl}_2$ style for $n \leq 6$ (except for the $6 \times 6$ nilpotent matrix having only one Jordan block), and also for any even $n$ if all Jordan blocks of $N$ are $2 \times 2$. Their method, which is now unnecessarily complicated, is based on the fact that vector fields $\mathbf{v}$ in $\mathfrak{sl}_2$ normal form are equivariant with respect to a certain one-parameter group; the method applies in principle to any nilpotent $N$, except that the calculations become too large.

Since the publication of [10], progress has been made in two directions: **1.** In [21], a method called **boosting** was given that obtains a description of the equivariant vector fields $\mathbf{v}$ from the (easier) description of the invariant scalar fields for the same group. The same paper shows that the solutions for the inner product and simplified styles follow from the solution for the $\mathfrak{sl}_2$ style. **2.** In [26] and [29, Ch. 12] a method called **box products** was given that obtains a description of the *invariants* for a given $N$ from separate descriptions of the invariants for each Jordan block in $N$. Very roughly, if

$$N = \begin{bmatrix} \acute{N} & 0 \\ 0 & \grave{N} \end{bmatrix}, \tag{1.4}$$

then the algebra $\mathcal{I}$ of invariants associated with $N$ can be viewed as a “box product” of the algebras of invariants associated with $\acute{N}$ and $\grave{N}$:

$$\mathcal{I} = \acute{\mathcal{I}} \boxtimes \grave{\mathcal{I}}. \tag{1.5}$$

(This becomes precise when viewed as a product of Stanley decompositions for the respective algebras; $\boxtimes$ distributes over direct sums in the decompositions.) There are close connections between the boosting and box product methods, since tensor product considerations underly both. Using the box product method we were able to solve the description problem for some $7 \times 7$ matrices $N$ such as $N_{223}$ and $N_{34}$, where the subscripts indicate Jordan block sizes. Thus the solution of the description problem breaks into three steps:

1. Solve the description problem for invariants of nilpotent matrices with one Jordan block. These are known for blocks up to $5 \times 5$, and are stated in Section 4 below.
2. Combine these by box products to describe invariants of nilpotent matrices with several blocks. A new version of this procedure is developed in Sections 6–16.
3. Boost the result from invariants to equivariants. A new version of this is given in Section 17.

The present paper was motivated by certain drawbacks, discussed in Remark 10.3, to the *expansion method* for computing box products given in [26] and [29]. The *factoring method* for box products, developed here, overcomes these drawbacks. For example, the Stanley decomposition in Example 16.4 below is far superior to the one given in [26] for the same $N$. Since this approach clarifies the foundations of both box products and boosting, and uses several new ideas and notations almost from the beginning, it seems time for a new exposition of the whole subject. Therefore this paper is written so as to be independent of [21,26], and [29].

The ultimate goal of this type of research (from the viewpoint of dynamical systems) is to obtain an algorithmic way to determine the bifurcation behavior of rest points with large nilpotent linear part. Such a program falls into three parts: describing the normal forms, unfolding them, and “blowing up” the unfoldings to detect the bifurcations. This paper is confined to the first part of such a program. The second (unfolding) is addressed in [23], and the third part (the blow-ups) should probably be addressed by specialists in the techniques of Takens, Bogdanov, A.D. Bruno [6,7], and F. Dumortier [13].

This paper is organized in three parts. Part 1 (Sections 2–5) gives the required background from $\mathfrak{sl}_2$ representation theory, invariant theory, and commutative algebra in the form and notation that we will use, collected from scattered sources. The box product is developed in Part 2 (Sections 6–16). Boosting is treated in Part 3 (Section 17). There are certain notation changes to be aware between the three parts. Part 1 requires only one coordinate system, $\mathbf{x} \in \mathbb{R}^n$, one nilpotent matrix $N$, and one algebra of invariants, $\mathcal{I}$. Part 2 requires three of each: $\mathbf{x} \in \mathbb{R}^n$, $\mathbf{y} \in \mathbb{R}^m$, and $(\mathbf{x}, \mathbf{y}) \in \mathbb{R}^{n+m}$; $\acute{N}$, $\grave{N}$, and $N$, related as in (1.4); and $\acute{\mathcal{I}}$, $\grave{\mathcal{I}}$, and $\mathcal{I}$, as in (1.5). The full set of notation changes for Part 2 is summarized in a table in Section 6. Part 3 involves a similar set of notation changes, listed in a table in Section 17. The object of boosting is to pass from $\mathcal{I}$ to $\mathcal{E}$, and the fundamental formula replacing (1.5) is

$$\mathcal{E} = \mathcal{I} \boxtimes \mathbb{R}^n. \tag{1.6}$$

The operation $\boxtimes$ in (1.6) is a variant of the box product in which the second factor, $\mathbb{R}^n$, is only a vector space and not an algebra.

The unusual use of accents in the notations of Part 2 calls for a comment. The algebras $\acute{\mathcal{I}}$ and $\grave{\mathcal{I}}$ in (1.5) are called the **left and right inputs** to the box product, and $\mathcal{I}$ is called the **output**. The accents ´ and ` are used on several symbols besides $N$ and $\mathcal{I}$. They can be regarded as pointing downward to the left and to the right, indicating that a symbol is associated with the left or right input, respectively; the absence of an accent indicates that the symbol is associated with the output.

A more detailed breakdown of Part 2 is as follows. **External transvectants** are defined in Section 7 and used to define the **pre-box product** in Section 8. The **replacement theorem** in Section 9 is used to convert the pre-box product into the box product. Since this is a fundamental idea, it is illustrated with easy examples in Section 10 before being developed rigorously in Section 11, which reduces the computation of box products to a problem in integer programming. Sections 12–15 develop algorithms that solve this integer programming problem; our algorithms seem to be simpler than the general algorithms of integer programming, probably because the transvectant cone is simpler than a general integer cone. Specifically, Section 12 illustrates the earliest factoring method that we discovered; this section could be skipped, but is helpful for motivation. Section 13 introduces a notation for sets of transvectants needed in the algorithms to follow. Sections 14 and 15 present the algorithms for box products. Section 16 contains examples of calculations using these algorithms.

Here is a guide to the location of the main definitions and notations.

1. Section 2 for normal form styles and the Lie algebra $\mathfrak{sl}_2$. Notations $\mathcal{D}_A, \mathsf{L}_A, X, Y, Z, \mathcal{X}, \mathcal{Y}, \mathcal{Z}, \mathsf{X}, \mathsf{Y}, \mathsf{Z}, \mathcal{I}, \mathcal{E}, N, M, H$.
2. Section 3 for Hilbert bases, preferred sets, standard preferred sets, Stanley decompositions, block decompositions, and matrices denoting sets of monomials. These "**b-matrices**" can be recognized because they are surrounded by brackets and contain two rows of integers (and sometimes one or more rows of column headings above a horizontal line). Notations $\boldsymbol{\alpha}, p, \widehat{\boldsymbol{\alpha}}, \mathbf{k}, \mathbf{t}, A, B, \mu$.
3. Section 5 for internal transvectants and Section 7 for external. Notation $(f, g)^{(s)}$.
4. Section 6 for notations involved with left inputs, right inputs, and outputs to the box product (some using accents, some using distinct letters: $\boldsymbol{\beta}, q, \boldsymbol{\ell}, \boldsymbol{\vartheta}, t, \mathbf{K}, \nu$). More detail about these is given in Section 11.
5. Section 8 for $\boldsymbol{\alpha\beta}$-transvectants, admissible transvectants, transvectant monomials, pre-box product. Notations $[\mathbf{k}; \boldsymbol{\ell}; s]$,$(\mathbf{k}; \boldsymbol{\ell}; s)$, $\widehat{\boxtimes}$ .
6. Section 11 for transvectant cone, transvectant map, primes, spanning sets, weak spanning sets, factoring order, preferred factorization, box product. Notations $\mathcal{C}, \mathcal{C}^*, \mathbb{T}, \boldsymbol{\pi}, \boldsymbol{\vartheta}, t, r, \Omega, \boxtimes_{\boldsymbol{\vartheta}}$.
7. Section 13 for matrices denoting sets of transvectants. These "**p-matrices**" look like b-matrices but are surrounded by parentheses instead of brackets. Notations $E, F, G, H$.

8. Section 15 for filtration, reconstruction, assembly, and special subsets.
9. Section 17 for notations involved with boosting from invariants to equivariants.

The sequel to this paper, [19], will include a technique (called the method of tokens) that simplifies the calculation of box products when either of the inputs has several Hilbert basis elements with the same weight; theorems about the location of prime transvectants in the transvectant cone; and theorems about relations among external transvectants.

# Part 1. Background

## 2. Basic definitions

Let $A$ be any real $n \times n$ matrix. A smooth scalar function $f : \mathbb{R}^n \to \mathbb{R}$ is **invariant** under the one-parameter group $e^{At}$ if $f(e^{At}\mathbf{x}) = f(\mathbf{x})$ for all $t \in \mathbb{R}$. This is equivalent to the vanishing of the directional derivative

$$\mathcal{D}_A f(\mathbf{x}) := f'(\mathbf{x})A\mathbf{x} = 0, \tag{2.1}$$

where $f'(\mathbf{x})$ is the row vector of partial derivatives of $f$. Similarly, a vector field $\mathbf{v}$ is **equivariant** under $e^{At}$ if the following Lie derivative vanishes:

$$\mathsf{L}_A\mathbf{v} := \mathbf{v}'(\mathbf{x})A\mathbf{x} - A\mathbf{v}(\mathbf{x}) = 0, \tag{2.2}$$

where $\mathbf{v}'$ is the matrix of partial derivatives of $\mathbf{v}$. Both of these definitions are also applicable when $f$ and $\mathbf{v}$ are formal power series. For each degree $k$, $\mathsf{L}_A$ maps the space $\mathcal{V}_k$ of homogeneous vector fields of degree $k$ into itself, so although $\mathsf{L}_A$ is a map of infinite-dimensional vector spaces, it can also be considered as a collection of maps of finite-dimensional spaces. Since

$$\mathsf{L}_A(f\mathbf{v}) = (\mathcal{D}_A f)\mathbf{v} + f\mathsf{L}_A\mathbf{v}, \tag{2.3}$$

the set $\mathcal{E}_A = \ker \mathsf{L}_A$ of formal equivariants is a module over the algebra $\mathcal{I}_A = \ker \mathcal{D}_A$ of formal invariants.

Again, let $A$ be any real $n \times n$ matrix. A **normal form style** for systems of the form

$$\dot{\mathbf{x}} = A\mathbf{x} + \mathbf{v}(\mathbf{x}) \tag{2.4}$$

is a choice of a complement $\mathcal{N}$ to the image of $\mathsf{L}_A$; for any such $\mathcal{N}$, every system (2.4) can be transformed by a change of coordinates into one such that $\mathbf{v} \in \mathcal{N}$ [22, Ch. 4]. If $A$ is semisimple, there is only one useful normal form style, the **semisimple normal form style**, defined by $\mathcal{N} = \mathcal{E}_A = \ker \mathsf{L}_A$. Since the term $A\mathbf{x}$ in (2.4) is already equivariant, taking $\mathbf{v}$ to be equivariant also makes the whole system invariant, which has important dynamical consequences. In the nilpotent case $A = N$, $\mathcal{E}_N$ is not complementary to im $\mathsf{L}_N$ and does not define a normal form style for (1.1). Instead, three normal form styles are important.

1. The **inner product normal form style** is defined by $\mathcal{N} = \mathcal{E}_{N^*} = \ker \mathsf{L}_{N^*}$, where $N^*$ is the transpose of $N$.

2. The $\mathfrak{sl}_2$ **normal form style** is defined by $\mathcal{N} = \mathcal{E}_M = \ker \mathsf{L}_M$, where $(N, M, H)$ is an $\mathfrak{sl}_2$ triad. This will be discussed in more detail below.
3. The most useful normal form style for applications is probably the **simplified normal form style**, which is defined based on the inner product style but is modified so that $\mathbf{v}$ is zero except in rows corresponding to bottom rows of nilpotent blocks of $N$. For instance, (1.3) is replaced by

$$\mathbf{v}(\mathbf{x}) = f(x_1)x_1 \begin{bmatrix} 0 \\ x_2 \end{bmatrix} + g(x_1)x_1^2 \begin{bmatrix} 0 \\ 1 \end{bmatrix} = \begin{bmatrix} 0 \\ f(x_1)x_1 + g(x_1)x_1^2 \end{bmatrix} \tag{2.5}$$

in simplified normal form. This greatly simplifies the appearance of the normal form, but loses the equivariance property of the inner product and $\mathfrak{sl}_2$ styles. However, the equivariance property seems to have no dynamical consequences, since the linear term $N\mathbf{x}$ in (1.1) does not share the equivariance of the nonlinear terms in normal form. As noted in the introduction, the usefulness of nilpotent normal forms results from unfolding and blowing-up arguments unrelated to equivariance.

Since our box product and boosting arguments depend heavily on $\mathfrak{sl}_2$, we use the $\mathfrak{sl}_2$ style here, but set things up in a way that yields correct results for the inner product style as well. Once a description for the inner product style has been found, it is trivial to obtain a description for the simplified style as well.

**2.1. Remark.** It is important to note that this “zeroing operation” **does not** map a system in inner product normal form to the simplified normal form *of the same system*, but is a one-to-one map of the set of *all* systems in inner product normal form onto the set of all systems in simplified normal form, so it **does** solve the strong description problem for the simplified normal form. See [22, §4.6] and [21] for more information.

**2.2. Remark.** The inner product normal form was discovered by Belitskii [4] and rediscovered by Elphick et al [15]. The simplified normal form occurred first in [15] in a few nilpotent examples. It was generalized to matrices that are not necessarily nilpotent in [22, §4.6], and used in [21] and [23]; for non-nilpotent cases the simplified normal form retains its (useful) equivariance with respect to the semisimple part of the linear term. (The earlier paper [20], and Section 6.4 of [22], are obsolete and are replaced by [23]. See the note on [20] in the references for more details.)

The Lie algebra $\mathfrak{sl}_2 = \mathfrak{sl}_2(\mathbb{R})$ is the span (over $\mathbb{R}$) of the matrices

$$\mathsf{x} = \begin{bmatrix} 0 & 1 \\ 0 & 0 \end{bmatrix}, \qquad \mathsf{y} = \begin{bmatrix} 0 & 0 \\ 1 & 0 \end{bmatrix}, \qquad \mathsf{z} = \begin{bmatrix} 1 & 0 \\ 0 & -1 \end{bmatrix},$$

equipped with the commutator bracket $[A, B] = AB - BA$, under which

$$[\mathsf{x}, \mathsf{y}] = \mathsf{z}, \quad [\mathsf{z}, \mathsf{x}] = 2\mathsf{x}, \quad [\mathsf{z}, \mathsf{y}] = -2\mathsf{y}.$$

An $\mathfrak{sl}_2$ **triad** is an ordered triple $(X, Y, Z)$ of $n \times n$ matrices (or of linear operators on a vector space $V$, which may be infinite-dimensional) satisfying

$$[X, Y] = Z, \quad [Z, X] = 2X, \quad [Z, Y] = -2Y.$$

A triad is equivalent to a **representation** of $\mathfrak{sl}_2$ on $\mathbb{R}^n$ (or $V$). The following theorem collects well-known facts (see Remark 2.4 below for references).

**2.3. Theorem.** *Let $(X, Y, Z)$ be an $\mathfrak{sl}_2$ triad of $n \times n$ matrices. Then $X$ and $Y$ are nilpotent, and $Z$ is semisimple, with integer eigenvalues called* **weights** *and eigenvectors called* **weight vectors**. *The* **top subspace** $\ker X \subset \mathbb{R}^n$ *has a (non-unique) basis $\mathbf{v}_1, \ldots, \mathbf{v}_k$ of weight vectors, known as* **top weight vectors**. *The weight of $\mathbf{v}_j$ will be denoted by $\widehat{\mathbf{v}}_j$, so $Z\mathbf{v}_j = \widehat{\mathbf{v}}_j \mathbf{v}_j$.*

*For each $\mathbf{v}_j$, the vectors $Y^i\mathbf{v}_j$ for $i = 0, \ldots, \widehat{\mathbf{v}}_j$ are nonzero weight vectors with decreasing weights $\widehat{\mathbf{v}}_j - 2i$. These vectors form a* **chain** *of length $\ell_j = \widehat{\mathbf{v}}_j + 1$, terminating with $Y^{\widehat{\mathbf{v}}_j}\mathbf{v}_j$, of weight $-\widehat{\mathbf{v}}_j$, which lies in the* **bottom space** *as a* **bottom weight vector**. *Each such chain spans an irreducible subrepresentation of $\mathfrak{sl}_2$. The $k$ chains with top weight vectors $\mathbf{v}_1, \ldots, \mathbf{v}_k$ form a basis for $\mathbb{R}^n$, which we call a* **chain-weight basis**. *The action of $X$ and $Y$ on these basis vectors is given by*

$$X(Y^i\mathbf{v}_j) = \begin{cases} \mathbf{0} & (i = 0) \\ i(\widehat{\mathbf{v}}_j + 1 - i)Y^{i-1}\mathbf{v}_j & (i = 1, \ldots, \widehat{\mathbf{v}}_j) \end{cases}$$

$$Y(Y^i\mathbf{v}_j) = \begin{cases} Y^{i+1}\mathbf{v}_j & (i = 0, \ldots, \widehat{\mathbf{v}}_j - 1) \\ \mathbf{0} & (i = \widehat{\mathbf{v}}_j) \end{cases}$$

*Thus $X$ maps up the chains, with its kernel at the top, while $Y$ maps down the chains, with its kernel at the bottom. Furthermore $Y$ is in lower Jordan form (with respect to this basis), and $X$ is in* **modified upper Jordan form**, *with the* **pressure factors** *$i(\widehat{\mathbf{v}}_j + 1 - i)$ in place of the off-diagonal ones that occur in (ordinary) upper Jordan form. The pressure on $Y^i\mathbf{v}_j$ is equal to the sum of the weights of the elements above $Y^i\mathbf{v}_j$ in the chain. These pressures are also the eigenvalues of $XY$, which is semisimple.*

**2.4. Remark.** These results are usually proved for $\mathfrak{sl}_2(\mathbb{C})$ rather than $\mathfrak{sl}_2(\mathbb{R})$, because the algebraic closure of $\mathbb{C}$ gives the initial existence of the eigenvalues of $Z$. See, for instance, [16]. However, every real $\mathfrak{sl}_2$ triad is also a complex $\mathfrak{sl}_2$ triad, and once the proof gets to the point that the weights are integers (and therefore real), it is clear that $Z$ is diagonalizable *over the reals* when $(X, Y, Z)$ is a real triad, so the whole proof goes through for $\mathfrak{sl}_2(\mathbb{R})$ as well. (There are great differences between the representation theory of $\mathfrak{sl}_2(\mathbb{R})$ and $\mathfrak{sl}_2(\mathbb{C})$, but these differences concern only infinite-dimensional representations.) The terminology of *pressure* in item 4 was introduced in [22] and used in [21] and [26].

Given an $\mathfrak{sl}_2$ triad $(X, Y, Z)$ on $\mathbb{R}^n$, a triad $(\mathcal{X}, \mathcal{Y}, \mathcal{Z})$ of differential operators on $\mathbb{R}[[\mathbf{x}]]$ is defined by:

$$\mathcal{X} := \mathcal{D}_Y, \quad \mathcal{Y} := \mathcal{D}_X, \quad \mathcal{Z} := \mathcal{D}_Z. \tag{2.6}$$

The interchange of $X$ and $Y$ enables $(\mathcal{X}, \mathcal{Y}, \mathcal{Z})$ to be an $\mathfrak{sl}_2$ triad (rather than an "anti-triad" in which $[\mathcal{Z}, \mathcal{X}] = -2\mathcal{X}$ and $[\mathcal{Z}, \mathcal{Y}] = 2\mathcal{Y}$). The operators $(\mathcal{X}, \mathcal{Y}, \mathcal{Z})$ preserve degree, and therefore restrict to triads on the finite-dimensional subspaces $\mathbb{R}[\mathbf{x}]_d$ of homogeneous polynomials in $\mathbf{x}$ of degree $d$. Thus it is not necessary to invoke functional analysis in the study of these triads, even though $\mathbb{R}[[\mathbf{x}]]$ is an infinite-dimensional vector space. A weight *vector* of $(\mathcal{X}, \mathcal{Y}, \mathcal{Z})$ is a *scalar* formal power series (usually a polynomial) $f$ that is an eigenvector of $\mathcal{Z}$; its weight will be denoted $\widehat{f}$, thus

$$\mathcal{Z}f = \widehat{f}f. \tag{2.7}$$

**2.5. Remark.** A more natural way to obtain a triad on $\mathbb{R}[[\mathbf{x}]]$ is to consider

$$\frac{d}{dt} f(e^{-Xt}\mathbf{x})\Big|_{t=0} = -f'(\mathbf{x})X\mathbf{x} = -\mathcal{D}_X f(\mathbf{x}),$$

and similarly for $Y$ and $Z$. Here the exponential maps from the Lie algebra $\mathfrak{sl}_2$ to the Lie group $\mathrm{SL}_2$, which acts on $R[[\mathbf{x}]]$ by $g \cdot f = f \circ g^{-1}$. (Using $g^{-1}$ makes this a group representation rather than an anti-representation.) The derivative returns to the Lie algebra. The resulting triad is $(-\mathcal{D}_X, -\mathcal{D}_Y, -\mathcal{D}_Z)$, with no interchange of $X$ and $Y$. It would then be necessary to use the chain bottoms (in $\ker(-\mathcal{D}_Y) = \ker \mathcal{D}_Y$), rather than tops, as the basis for $\mathcal{I}$.

**2.6. Remark.** Note that $\mathbb{R}[\mathbf{x}]_d \subset \mathbb{R}[\mathbf{x}] \subset \mathbb{R}[[\mathbf{x}]]$. As a vector space, $\mathbb{R}[\mathbf{x}]$ is the direct sum of the subspaces $\mathbb{R}[\mathbf{x}]_d$, and a natural vector space basis for $\mathbb{R}[\mathbf{x}]$ is the set of monomials in $\mathbf{x}$; $\mathbb{R}[[\mathbf{x}]]$ is the direct *product* of the same subspaces, but its elements are written in "sum notation" (as formal power series) rather than "product notation" (as $\infty$-tuples of homogeneous monomials). There is no countable vector space basis for $\mathbb{R}[[\mathbf{x}]]$, but the countable set of monomials in $\mathbf{x}$ is a **formal vector space basis** for $\mathbb{R}[[\mathbf{x}]]$ in the sense that these monomials are linearly independent and every element of $\mathbb{R}[[\mathbf{x}]]$ is a **countable linear combination** of monomials. Since $\mathcal{X}$, $\mathcal{Y}$, and $\mathcal{Z}$ preserve degree, an element $f$ of $\mathbb{R}[[\mathbf{x}]]$, belongs to $\ker \mathcal{X}$ if and only if each homogeneous component of $f$ belongs to $\ker \mathcal{X}$. It follows that **the theories of polynomial invariants ($\mathcal{I} \cap \mathbb{R}[\mathbf{x}]$) and formal invariants ($\mathcal{I}$) are exactly the same, except for the use of countable linear combinations in the formal case**; we therefore use the word **span** ambiguously to include both cases. By the same token, we can use Gröbner basis arguments (that properly apply only to $\mathbb{R}[\mathbf{x}]$) and carry the results over to $\mathbb{R}[[\mathbf{x}]]$ (where *standard basis* arguments would usually be required instead); this is only possible because the basis elements are homogeneous. For more about formal vector spaces and formal bases, see [17]. For *standard bases*, see [8, §4.4].

According to (2.1), the algebra $\mathcal{I}$ of formal invariants of $e^{Yt}$ is given by

$$\mathcal{I} = \ker \mathcal{X} = \ker \mathcal{D}_Y; \tag{2.8}$$

these can also be called **invariants of** $\mathcal{X}$ and **seminvariants of** $(\mathcal{X}, \mathcal{Y}, \mathcal{Z})$. (Seminvariant means half-invariant. An **invariant of** $(\mathcal{X}, \mathcal{Y}, \mathcal{Z})$ is required to satisfy $\mathcal{X}f = 0$, $\mathcal{Y}f = 0$, and $\mathcal{Z}f = 0$, but the latter follows from the first two. So if $\mathcal{X}f = 0$, $f$ is half-invariant under the triad.)

**2.7. Remark.** The terminology of seminvariants goes back to Cayley and Sylvester, and belongs to the "classical" invariant theory of binary forms. The primary objects of interest in classical invariant theory were the *covariants*, which we will not use or define. The leading term of a covariant, called the *source* of the covariant, is a seminvariant, and this association gives an isomorphism between the algebras of seminvariants and of covariants. The result is originally due to M. Roberts. An easy classical explanation of the covariant-seminvariant correspondence is found in [12]; for the attribution to M. Roberts, see [14, §§109, 111]. Another proof is given in [18, Lecture XII]. A good modern introduction to classical invariant theory is [27], but the discussion of seminvariants requires corrections that have been posted online (see the note in the reference). Some modern authors, such as [11], confuse seminvariants with another classical notion, that of *relative invariants*.

Similarly, given $(X, Y, Z)$ we can form the triad $(\mathsf{X}, \mathsf{Y}, \mathsf{Z})$ of Lie operators

$$\mathsf{X} := \mathsf{L}_Y, \quad \mathsf{Y} := \mathsf{L}_X, \quad \mathsf{Z} := \mathsf{L}_Z \tag{2.9}$$

acting on vector fields $\mathbf{v} : \mathbb{R}^n \to \mathbb{R}^n$. Then $\mathcal{E} = \ker \mathsf{X}$ is the space of **equivariants** of $\mathsf{X}$ or $\mathsf{L}_Y$ or $e^{Yt}$, also called **semi-equivariants** of $(\mathsf{X}, \mathsf{Y}, \mathsf{Z})$. As before, $\mathcal{E}$ is a module over $\mathcal{I} = \ker \mathcal{X}$. The system of differential equations

$$\dot{\mathbf{x}} = X\mathbf{x} + \mathbf{v}(\mathbf{x}) \tag{2.10}$$

is in $\mathfrak{sl}_2$ **normal form** with respect to $(X, Y, Z)$ if $\mathbf{v} \in \mathcal{E}$; note that $\mathfrak{sl}_2$ normal form is not defined unless the whole triad $(X, Y, Z)$ is specified, because $X$ does not uniquely determine $\mathsf{X} = \mathsf{L}_Y$. As with inner product normal forms, every system (2.10) can be put into $\mathfrak{sl}_2$ normal form with respect to $(X, Y, Z)$ by a near-identity transformation of variables [22, §4.8]. The description problem for $\mathfrak{sl}_2$ normal forms is to find a formula that gives each element of $\mathcal{E}$ in a unique manner, for a given choice of $(X, Y, Z)$.

Now we return to the original system of differential equations (1.1) with fixed nilpotent $N$. By the Jacobson–Morosov lemma [22, Alg. 2.7.2], there exist many choices of matrices $M$, $H$ such that $(N, M, H)$ is an $\mathfrak{sl}_2$ triad. For any such choice, it follows that $(M^*, N^*, H^*)$ is also an $\mathfrak{sl}_2$ triad, where $*$ denotes transpose. The definitions of inner product and $\mathfrak{sl}_2$ normal forms imply the following lemma.

**2.8. Lemma.** *The system* $\dot{\mathbf{x}} = N\mathbf{x} + \mathbf{v}(\mathbf{x})$ *is in inner product normal form if and only if the system* $\dot{\mathbf{x}} = M^*\mathbf{x} + \mathbf{v}(\mathbf{x})$ *is in* $\mathfrak{sl}_2$ *normal form with respect to* $(M^*, N^*, H^*)$.

**Proof.** Let $(X, Y, Z) = (M^*, N^*, H^*)$. Then $\dot{\mathbf{x}} = M^*\mathbf{x} + \mathbf{v}(\mathbf{x})$ is in $\mathfrak{sl}_2$ normal form with respect to $(X, Y, Z)$ if and only if $\mathbf{v} \in \ker \mathsf{X} = \ker \mathsf{L}_Y = \ker \mathsf{L}_{N^*}$. This is the same as the condition for (1.1) to be in inner product normal form with respect to $N$. □

**2.9. Remark.** It must be understood that the two systems $\dot{\mathbf{x}} = N\mathbf{x} + \mathbf{v}(\mathbf{x})$ and $\dot{\mathbf{x}} = M^*\mathbf{x} + \mathbf{v}(\mathbf{x})$ with $\mathbf{v} \in \mathcal{E}$ are in normal form (with different styles), *but will not be normal forms of the same system or of each other.* This is because the projections used in solving the computation problem (to put a particular system into normal form) differ for the two normal form styles. The projections for the inner product normal form are orthogonal with respect to an inner product, and the computation of these projections requires solution of systems of linear equations, which can be quite large. The projection computation for $\mathfrak{sl}_2$ normal forms requires only iteration of known operations that do not increase in size in higher degrees. See [22].

**For the remainder of this paper, we adopt the following notational conventions for all computational examples. In theoretical contexts these specifications need not hold.**

- $N$ is a nilpotent matrix in upper Jordan form. Subscripts indicate block sizes, so that

$$N_4 = \begin{bmatrix} 0 & 1 & 0 & 0 \\ 0 & 0 & 1 & 0 \\ 0 & 0 & 0 & 1 \\ 0 & 0 & 0 & 0 \end{bmatrix}, \qquad N_{23} = \begin{bmatrix} 0 & 1 & & & \\ 0 & 0 & & & \\ & & 0 & 1 & 0 \\ & & 0 & 0 & 1 \\ & & 0 & 0 & 0 \end{bmatrix}.$$

- $(N, M, H)$ is the specific $\mathfrak{sl}_2$ triad built from $N$ as follows. For each $k \times k$ block $N_k$ occurring in $N$, the matrices $M$ and $H$ contain the $k \times k$ blocks

$$M_k = \begin{bmatrix} 0 & & & & & \\ (1)(k-1) & 0 & & & & \\ & (2)(k-2) & 0 & & & \\ & & (3)(k-3) & & & \\ & & & \ddots & & \\ & & & & (k-1)(1) & 0 \end{bmatrix} \tag{2.11}$$

and

$$H_k = \begin{bmatrix} (k-1) & & & & \\ & (k-3) & & & \\ & & (k-5) & & \\ & & & \ddots & \\ & & & & -(k-1) \end{bmatrix}. \tag{2.12}$$

For instance,

$$M_4 = \begin{bmatrix} 0 & 0 & 0 & 0 \\ 3 & 0 & 0 & 0 \\ 0 & 4 & 0 & 0 \\ 0 & 0 & 3 & 0 \end{bmatrix}, \quad H_4 = \begin{bmatrix} 3 & 0 & 0 & 0 \\ 0 & 1 & 0 & 0 \\ 0 & 0 & -1 & 0 \\ 0 & 0 & 0 & -3 \end{bmatrix}.$$

- $(X, Y, Z) = (M^*, N^*, H)$. Since $H$ is real and diagonal, $H = H^*$. Thus

$$X_4 = \begin{bmatrix} 0 & 3 & 0 & 0 \\ 0 & 0 & 4 & 0 \\ 0 & 0 & 0 & 3 \\ 0 & 0 & 0 & 0 \end{bmatrix}, \quad Y_4 = \begin{bmatrix} 0 & 0 & 0 & 0 \\ 1 & 0 & 0 & 0 \\ 0 & 1 & 0 & 0 \\ 0 & 0 & 1 & 0 \end{bmatrix}, \quad Z_4 = \begin{bmatrix} 3 & 0 & 0 & 0 \\ 0 & 1 & 0 & 0 \\ 0 & 0 & -1 & 0 \\ 0 & 0 & 0 & -3 \end{bmatrix}.$$

- $(\mathcal{X}, \mathcal{Y}, \mathcal{Z})$ is given by

$$\mathcal{X} = \mathcal{D}_Y = \mathcal{D}_{N^*} = (N^*\mathbf{x}) \cdot \nabla$$
$$\mathcal{Y} = \mathcal{D}_X = \mathcal{D}_{M^*} = (M^*\mathbf{x}) \cdot \nabla$$
$$\mathcal{Z} = \mathcal{D}_Z = \mathcal{D}_H = (H^*\mathbf{x}) \cdot \nabla.$$

Subscripts on operators are again the block sizes of $N$. For instance

$$\mathcal{X}_4 = x_1\frac{\partial}{\partial x_2} + x_2\frac{\partial}{\partial x_3} + \frac{\partial}{\partial x_4},$$
$$\mathcal{Y}_4 = 3x_2\frac{\partial}{\partial x_1} + 4x_3\frac{\partial}{\partial x_2} + 3x_4\frac{\partial}{\partial x_3},$$
$$\mathcal{Z}_4 = 3x_1\frac{\partial}{\partial x_1} + x_2\frac{\partial}{\partial x_2} - x_3\frac{\partial}{\partial x_3} - 3x_4\frac{\partial}{\partial x_4}.$$

- $(\mathsf{X}, \mathsf{Y}, \mathsf{Z})$ is given by

$$\begin{aligned} \mathsf{X} &= \mathsf{L}_Y = \mathsf{L}_{N^*} = \mathcal{X} - N^* \\ \mathsf{Y} &= \mathsf{L}_X = \mathsf{L}_{M^*} = \mathcal{Y} - M^* \\ \mathsf{Z} &= \mathsf{L}_Z = \mathsf{L}_H = \mathcal{Z} - H, \end{aligned} \tag{2.13}$$

  where $\mathcal{X}$, $\mathcal{Y}$, and $\mathcal{Z}$ operate componentwise on vector fields $\mathbf{v}(\mathbf{x})$ and $N^*, M^*, H$ act by matrix multiplication on $\mathbf{v}(\mathbf{x})$.
- $\mathcal{I} = \ker \mathcal{X}$ and $\mathcal{E} = \ker \mathsf{X}$. Subscripts on $\mathcal{I}$ or $\mathcal{E}$ are the block sizes of $N$.

With these conventions, $\dot{\mathbf{x}} = N\mathbf{x} + \mathbf{v}(\mathbf{x})$ is in inner product normal form, and $\dot{\mathbf{x}} = M^*\mathbf{x} + \mathbf{v}(\mathbf{x})$ is in $\mathfrak{sl}_2$ normal form, if and only if $\mathbf{v} \in \mathcal{E}$. **Therefore a solution of the description problem for the $\mathfrak{sl}_2$ normal form with linear part $M^*$ is also a solution of the description problem for the inner product normal form with linear part $N$.**

## 3. Hilbert bases, preferred sets, block and Stanley decompositions

A **weight seminvariant** of $(\mathcal{X}, \mathcal{Y}, \mathcal{Z})$ is a seminvariant $f \in \ker \mathcal{X}$ that is also a weight vector (an eigenvector of $\mathcal{Z}$), so that $\mathcal{Z}f = \widehat{f} f$. For convenience, we will refer to such an $f$ as a **weight invariant** of $\mathcal{X}$, although the full triad must be specified in order for this to make sense. By Theorem 2.3, each finite-dimensional space $\mathcal{I} \cap \mathbb{R}[\mathbf{x}]_d$ has a basis consisting of homogeneous weight invariants of degree $d$. A **Hilbert basis** for the algebra $\mathcal{I} \cap \mathbb{R}[\mathbf{x}]$ of *polynomial* invariants of $\mathcal{X}$ is a finite ordered set $\boldsymbol{\alpha} = (\alpha_1, \dots, \alpha_p)$ of homogeneous weight invariants such that every element $f \in \mathcal{I} \cap \mathbb{R}[\mathbf{x}]$ is expressible as a polynomial function of $\boldsymbol{\alpha}$, or equivalently, as a finite linear combination of the monomials $\boldsymbol{\alpha}^{\mathbf{k}} = \alpha_1^{k_1} \cdots \alpha_p^{k_p}$ in the Hilbert basis elements. In view of Remark 2.6, we can extend this definition to say that $\boldsymbol{\alpha}$ is a Hilbert basis for the algebra $\mathcal{I}$ of *formal* invariants if every $f \in \mathcal{I}$ is a *countable* linear combination of the $\boldsymbol{\alpha}^{\mathbf{k}}$, and then observe that $\boldsymbol{\alpha}$ is a Hilbert basis for $\mathcal{I}$ if and only if it is a Hilbert basis for $\mathcal{I} \cap \mathbb{R}[\mathbf{x}]$. Associated with a Hilbert basis $\boldsymbol{\alpha}$ is its **row vector of weights** $\widehat{\boldsymbol{\alpha}} = \begin{bmatrix}\widehat{\alpha}_1 & \dots & \widehat{\alpha}_p\end{bmatrix}$. Each monomial $\boldsymbol{\alpha}^{\mathbf{k}} = \alpha_1^{k_1} \cdots \alpha_p^{k_p}$ is a homogeneous weight invariant of degree $k_1 \deg \alpha_1 + \cdots + k_p \deg \alpha_p$ with weight

$$\widehat{\boldsymbol{\alpha}^{\mathbf{k}}} = \widehat{\boldsymbol{\alpha}}\mathbf{k} = \widehat{\alpha}_1 k_1 + \cdots + \widehat{\alpha}_p k_p. \tag{3.1}$$

**3.1. Remark.** Let $\mathbb{R}[\mathbf{x}]_{dw}$ be the subspace of $\mathbb{R}[\mathbf{x}]_d$ that are also weight vectors with weight $w$. Then

$$\mathbb{R}[\mathbf{x}] = \bigoplus_d \bigoplus_w \mathbb{R}[\mathbf{x}]_{dw},$$

so $\mathbb{R}[\mathbf{x}]$ is doubly graded by degree and weight. The same is true for $\mathbb{R}[[\mathbf{x}]]$ if countable sums are allowed. This restricts to a double grading of $\mathcal{I}$.

**3.2. Remark.** In our examples, $Z$ is diagonal (and not merely diagonalizable). This makes it easy to recognize the weight functions of $\mathcal{Z}$. Let $Z = \mathrm{diag}(w_1, \ldots, w_n)$. Then every monomial $\mathbf{x}^{\mathbf{m}} = x_1^{m_1} \cdots x_n^{m_n}$ in $\mathbf{x}$ is a weight vector of weight $m_1 w_1 + \cdots + m_n w_n$, and a polynomial in $\mathbf{x}$ is a weight vector if and only if all of its monomials have the same weight. Such a polynomial is called **isobaric**. (This notion is classical, even though the classical definition of weight is different from the $\mathfrak{sl}_2$ definition; the same polynomials are isobaric under either definition. To add the confusion, there is a third, unrelated classical meaning of weight, also called index, that arises in connection with relative invariants and covariants.)

Although $\alpha_1^{k_1} \cdots \alpha_p^{k_p}$, when regarded as a *symbol string*, has the appearance of a monomial, as a *mathematical object* it is a polynomial in $\mathbf{x}$; each $\alpha_i$ is already a polynomial in $\mathbf{x}$, so the product of these is a polynomial. To avoid this and other difficulties with "monomials in $\boldsymbol{\alpha}$," it is customary to introduce a set of indeterminates $t_1, \ldots, t_p$ equal in number to $\alpha_1, \ldots, \alpha_p$, and to let $\Phi$ be the unique algebra homomorphism $\Phi : \mathbb{R}[\mathbf{t}] \to \mathbb{R}[\mathbf{x}]$ satisfying $\Phi(t_i) = \alpha_i$ for $i = 1, \ldots, p$. The kernel of $\Phi$ is an ideal $J \subset \mathbb{R}[\mathbf{t}] = \mathbb{R}[t_1, \ldots, t_p]$ called the **ideal of relations** among $\boldsymbol{\alpha}$. For instance, we will see later that there is a Hilbert basis with $p = 4$ for $\mathcal{I}_4$ in which $\alpha_3^2 - \alpha_2^3 + \alpha_1^2\alpha_4$ is identically zero as a polynomial in $\mathbf{x}$. While it is correct to write $\alpha_3^2 - \alpha_2^3 + \alpha_1^2\alpha_4 = 0$, the left-hand side of this equation, as a mathematical object rather than as a symbol string, *is* $0 \in \mathbb{R}[\mathbf{x}]$, so the equation reduces to $0 = 0$, which does not identify the intended relation among $\boldsymbol{\alpha}$ at all. It is better to say $t_3^2 - t_2^3 + t_1^2 t_4 \in J$, since $t_3^2 - t_2^3 + t_1^2 t_4$ is not zero and does contain the information that identifies the relation. But our practice will be to speak freely of monomials and polynomials in $\boldsymbol{\alpha}$, treating them sometimes as symbol strings and sometimes as polynomials in $\mathbf{x}$ according to context, except when to do so would create mathematical confusion. In such cases we will invoke the machinery of $\mathbf{t}$ and $\Phi$.

When $J$ is nontrivial, the representation of an invariant as a polynomial in $\boldsymbol{\alpha}$ is not unique. Therefore a Hilbert basis $\boldsymbol{\alpha}$ is usually not enough to solve the description problem for $\mathcal{I}$, which calls for a unique expression for each invariant. But there is a simple solution. The monomials $\boldsymbol{\alpha}^{\mathbf{k}}$ span $\mathcal{I}$, so any linearly independent subset $A$ of these monomials constitutes a formal vector space basis for $\mathcal{I}$. We call $A$ a **set of preferred monomials** in $\boldsymbol{\alpha}$, or simply a **preferred set**. *Any invariant is uniquely expressible as a (countable) linear combination of preferred monomials.* When $\mathbf{t}$ notation is used, $A$ is regarded as a subset of the monomials in $\mathbf{t}$ such that $\Phi(A)$ is a formal vector space basis for $\mathcal{I}$. (Thus $A$ is a set of monomials in $\boldsymbol{\alpha}$, or the corresponding set of monomials in $\mathbf{t}$, according to the context.)

**3.3. Remark.** Notice the many meanings of "unique". The Hilbert basis for $\mathcal{I}$ *is not* unique. The preferred set $A$ for $\mathcal{I}$ *is not* unique; it depends on a choice of the Hilbert basis $\boldsymbol{\alpha}$, and the choice of a linearly independent set $A$ of monomials in $\boldsymbol{\alpha}$. But after $A$ is chosen, the expression for any given polynomial invariant *is* unique.

Since there exist infinite sets that are not finitely describable, we now limit the preferred sets $A$ that we consider to those that can be described finitely using either *block decompositions* or *Stanley decompositions*. Let $\mathbb{N} = \mathbb{Z}_{\geq 0}$ be the set of nonnegative integers, and let $\mathbf{k} \in \mathbb{N}^p$ represent the monomial $\boldsymbol{\alpha}^{\mathbf{k}}$ or $\mathbf{t}^{\mathbf{k}}$; when used in this way, $\mathbb{N}^p$ is called a **Newton space**. If $a \in \mathbb{N}$ and $b \in \mathbb{N} \cup \{\infty\}$, let the **nonnegative integer interval** $[a, b]$ be defined by

$$[a,b] = \begin{bmatrix} b \\ a \end{bmatrix} = \begin{cases} \{j \in \mathbb{N} : a \leq j \leq b\} & \text{if} \quad b \neq \infty \\ \{j \in \mathbb{N} : a \leq j < \infty\} & \text{if} \quad b = \infty. \end{cases} \tag{3.2}$$

The column form of this notation is read from bottom to top, like the limits on an integral sign. Writing $[a, \infty]$ as if it were a closed interval (when $[a, \infty)$ would be more natural) allows us to handle the cases $b < \infty$ and $b = \infty$ with a single notation. Note that $[a, b]$ is empty if $b < a$. If $p$ such intervals, written as columns, are entered into a **b-matrix** (a matrix surrounded by square brackets), the matrix denotes the Cartesian product of the intervals:

$$B = \begin{bmatrix} \mathbf{b} \\ \mathbf{a} \end{bmatrix} = \begin{bmatrix} b_1 & \cdots & b_p \\ a_1 & \cdots & a_p \end{bmatrix} = [a_1, b_1] \times \cdots \times [a_p, b_p] \subset \mathbb{N}^p. \tag{3.3}$$

The letter $B$ is used ambiguously to denote the b-matrix itself, the Cartesian product as a subset of $\mathbb{N}^p$, and the corresponding sets of monomials $\boldsymbol{\alpha}^{\mathbf{k}}$ and $\mathbf{t}^{\mathbf{k}}$; it will be clear from the context which is intended. Such sets (of integer vectors or of monomials) will be called **blocks**.

**3.4. Remark.** In later examples we sometimes add "decorations" to the b-matrix to include additional information for convenience. See, for instance, (10.2) and (10.6).

**3.5. Remark.** Beginning in Section 13, it will be helpful to generalize this notation by allowing $a$ and $b$ in (3.2) to be negative, while still requiring $j \in \mathbb{N}$ on the right-hand side, so that $[a, b]$ is still a nonnegative interval. For instance, $[-2, 3] = \{0, 1, 2, 3\}$ and $[-5, -1] = \emptyset$. Then, in (3.3), the bottom row can contain negative entries but represents the same block as if these entries were zero.

A **block decomposition** of a preferred set $A$ is a representation of $A$ as a union of disjoint blocks:

$$A = B^1 \sqcup \cdots \sqcup B^\mu \tag{3.4}$$

Here $\sqcup$ denotes *disjoint union* (in the sense of a union of sets that are *already* disjoint), and $\mu$ is the number of blocks. It follows immediately from (3.4) that

$$\mathcal{I} = \text{Span } A = \text{Span } B^1 \oplus \cdots \oplus \text{Span } B^\mu, \tag{3.5}$$

where each block $B^i$ is to be understood as a set of $\boldsymbol{\alpha}^{\mathbf{k}}$, regarded as polynomials in $\mathbf{x}$, and the spans are *countable* spans whenever $B^i$ is an infinite set (see Remark 2.6).

Block decompositions were introduced in [24] and [25] as an alternative to Stanley decompositions. A **Stanley space** in $\mathbb{R}[\mathbf{t}]$ is a vector subspace of the form $\mathbb{R}[T]\varphi$, where $T$ is a subset of the indeterminates $\{t_1, \ldots, t_p\}$ and $\varphi$ is a monomial in $\mathbf{t}$; $\mathbb{R}[T]\varphi$ is the span of the set of multiples of $\varphi$ by polynomials in the indeterminates listed in $T$. To state this in terms of blocks, let $B$ be the block (3.3) with bottom row $\mathbf{a}$ such that $\mathbf{t}^{\mathbf{a}} = \varphi$ and top row $\mathbf{b}$ such that $b_i = \infty$ if $t_i \in T$, $b_i = a_i$ otherwise. Then we refer to $B$ as a **Stanley block**, and observe that Span $B = \mathbb{R}[T]\varphi$. In Newton space, the inner corner (closest to $\mathbf{0}$) of $B$ represents $\varphi$, and in every direction $\mathbf{e}_i$, the block either extends to infinity or has no thickness at all. A block that is *not* Stanley has at least one direction in which its thickness is finite but not zero. Any block can be written as a disjoint union of Stanley blocks, and there is a unique *minimal* way of doing this (using the fewest Stanley blocks); see [24, Alg. 2.2], and Example (4.11) below, illustrating a block rewritten as a disjoint union of two Stanley blocks; the original block is non-Stanley in direction $\mathbf{e}_3$, and the two Stanley blocks are "slices" of the given block perpendicular to $\mathbf{e}_3$. The principal difference between block decompositions and Stanley decompositions is that a block decomposition describes a set of monomials, while a Stanley decomposition describes the space spanned by those monomials. Thus a block decomposition is a *disjoint union* of blocks, whereas a Stanley decomposition is a *direct sum* of Stanley spaces, having the form

$$\mathbb{R}[T_1]\varphi_1 \oplus \cdots \oplus \mathbb{R}[T_\mu]\varphi_\mu. \tag{3.6}$$

All of this extends to $\mathbb{R}[[\mathbf{t}]]$ if ordinary spans are replaced by countable spans. For us, the advantages of block decompositions are that they are more concise, and more suitable for the algorithms in Section 14 below. The advantage of a Stanley decomposition is that it is readily converted into an expression for $\mathcal{I}$ containing finitely many arbitrary formal power series, as in (4.12) and (4.13) below. Therefore we do our box product computations using block decompositions, and convert these to Stanley decompositions afterwards.

A **monomial ideal** $I \subset \mathbb{R}[\mathbf{t}]$ is an ideal generated by monomials. The **standard monomials** of a monomial ideal $I$ are the monomials in $\mathbb{R}[\mathbf{t}]$ that are **not** in $I$. Conversely, we will call a set $A$ of monomials in $\mathbb{R}[\mathbf{t}]$ a **standard set** if $A$ is the set of standard monomials for some monomial ideal $I$. The following lemma characterizes standard sets without mentioning the associated ideal.

**3.6. Lemma.** *A set $A$ of monomials in $\mathbb{R}[\mathbf{t}]$ is standard if and only if every monomial in $\mathbb{R}[\mathbf{t}]$ that divides a monomial in $A$ belongs to $A$.*

**Proof.** Suppose that $A$ is standard, with associated monomial ideal $I$. Let $F \in A$ and let $g$ be a monomial that divides $F$, so that $F = fg$ for some monomial $f$. We claim that $g \in A$. Since $A$ is standard, if $g \notin A$ then $g \in I$. Then since $I$ is an ideal, $F = fg \in I$, contradicting $F \in A$.

Conversely, suppose every divisor of a monomial in $A$ belongs to $A$. Let $M$ be the set of monomials not in $A$, suppose $g \in M$, and let $f$ be any monomial. Then $fg \in M$, since if it belonged to $A$, its divisor $g$ would belong to $A$ rather than $M$. Therefore $M$ absorbs monomials under multiplication, and it follows easily that $I = \text{Span}\, M$ absorbs polynomials. Since $I$ is also closed under addition, it is an ideal, and since it is generated by monomials, it is a monomial ideal. □

**3.7. Remark.** There is an extensive literature on Stanley decompositions in the context of commutative algebra, where they were introduced as a tool for studying the concept of "depth", which does not concern us. The only part of this literature that we use is [32, pp. 278–279]; the method given there produces the Stanley decompositions for $\mathcal{I}_1$–$\mathcal{I}_5$ given in the next section, but is not used again in this paper. Here is a summary of the ideas in those pages. Let $J \subset \mathbb{R}[\mathbf{t}]$ be the ideal of relations among $\boldsymbol{\alpha}$. Then $\mathcal{I} \cap \mathbb{R}[\mathbf{x}]$ is isomorphic as an algebra to $\mathbb{R}[\mathbf{t}]/J$, where the isomorphism is induced by $\Phi : \mathbb{R}[\mathbf{t}] \to \mathbb{R}[\mathbf{x}]$. Choose a Gröbner term order for $\mathbb{R}[\mathbf{t}]$ and let $I \subset \mathbb{R}[\mathbf{t}]$ be the ideal of leading terms of $J$. Then $I$ is a monomial ideal, and $\mathcal{I} \cap \mathbb{R}[\mathbf{x}]$ is isomorphic as a vector space (but not as an algebra) to $\mathbb{R}[\mathbf{t}]/I$; see the discussion of equation (4.14) below. It follows that $\mathcal{I} \cap \mathbb{R}[\mathbf{x}]$ is spanned by the monomials $\boldsymbol{\alpha}^{\mathbf{k}} = \Phi(\mathbf{t}^{\mathbf{k}})$ as $\mathbf{t}$ ranges over the standard monomials of $I$, so these monomials form a preferred set $A$. In view of Remark 2.6, this applies to $\mathcal{I}$ also if ordinary span is replaced by countable span. To apply these ideas requires some algorithms. There is a well-known algorithm (see for instance [1, p. 81]) to compute a Gröbner basis for $J = \ker \Phi$ using an eliminating order on $\mathbb{R}[\mathbf{x}, \mathbf{t}]$. The leading monomials of the Gröbner basis for $J$ generate $I$. Lemma 2.4 in [32] gives an algorithm to produce a Stanley decomposition for the span of the standard monomials of $I$, which (using countable spans) is $\mathcal{I}$. Alternatively, algorithms in [24] and [25] could be used to produce a block decomposition.

Preferred sets obtained by Remark 3.7 are always standard. Since the factoring method for box products developed in this paper requires standard preferred sets as inputs, and produces a standard preferred set as its output, **we assume from here on that all preferred sets mentioned are standard unless otherwise noted**.

**3.8. Remark.** Lemma 3.6 shows that the preferred set given by the Stanley decomposition for $\mathcal{I}_{223}$ found by the expansion method in [26, p. 252] is not standard, since it describes a preferred set that contains $\alpha\beta\zeta$ (in the term $\mathcal{R}[[\alpha, \beta, \delta]]\alpha\zeta$) but not $\beta\zeta$ (since the only other term involving $\zeta$ is $\mathcal{R}[[\delta, \mu]]\zeta$). Thus the expansion method is prone to producing nonstandard preferred sets. The preferred set is still valid, but a nonstandard preferred set is never as simple as a standard one. We will give a better Stanley decomposition for $\mathcal{I}_{223}$ in Example 16.4.

## 4. Descriptions of $\mathcal{I}_1$ through $\mathcal{I}_5$

We are now prepared to give descriptions of $\mathcal{I}_1, \ldots, \mathcal{I}_5$. These are the building blocks from which box products will be computed. These results can be obtained using Remark 3.7; see [9] and [10], which use different conventions than ours for the $\mathfrak{sl}_2$ triad.

The case of $\mathcal{I}_1$ is trivial and is usually omitted; $\mathcal{X}_1$ is the zero operator and $\mathcal{I}_1$ is the entire space $\mathbb{R}[[\mathbf{x}]] = \mathbb{R}[[x_1]]$. The Hilbert basis is the single polynomial $\alpha_1 = x_1$, with degree and weight as shown in the table

| polynomial | degree | weight |
|---|---|---|
| $\alpha_1 = x_1$ | 1 | 0 |

. (4.1)

The block decomposition is

$$A_1 = \begin{bmatrix} \infty \\ 0 \end{bmatrix}, \tag{4.2}$$

and the Stanley decomposition is

$$\mathbb{R}[[\alpha_1]] = \mathbb{R}[[x_1]]. \tag{4.3}$$

(This case is also trivial at the level of equivariants; every system is in normal form.)

For $\mathfrak{I}_2$, the Hilbert basis is

| polynomial | degree | weight |
|---|---|---|
| $\alpha_1 = x_1$ | 1 | 1 |

(4.4)

the block decomposition is

$$A_2 = \begin{bmatrix} \infty \\ 0 \end{bmatrix}, \tag{4.5}$$

and the Stanley decomposition is

$$\mathfrak{I}_2 = \mathbb{R}[[\alpha_1]] = \mathbb{R}[[x_1]]. \tag{4.6}$$

This says that every invariant can be written uniquely in the form $f(x_1, x_2) = F(x_1)$, where $F$ is a formal power series. Although (4.6) looks like (4.3), the meaning is different because $\mathfrak{I}_2$ is no longer all of $\mathbb{R}[[\mathbf{x}]]$.

For $\mathfrak{I}_3$, the Hilbert basis is

| polynomial | degree | weight |
|---|---|---|
| $\alpha_1 = x_1$ | 1 | 2 |
| $\alpha_2 = x_2^2 - 2x_1x_3$ | 2 | 0 |

(4.7)

(The expression for $\alpha_2$ is a version of the familiar $b^2 - 4ac$ for the quadratic form. It has a different scaling because of our notational conventions.) The block decomposition is

$$A_3 = \begin{bmatrix} \infty & \infty \\ 0 & 0 \end{bmatrix} \tag{4.8}$$

and the Stanley decomposition is

$$\mathfrak{I}_3 = \mathbb{R}[[\alpha_1, \alpha_2]]. \tag{4.9}$$

Every invariant $f$ can be written uniquely in the form

$$f(x_1, x_2, x_3) = F(\alpha_1, \alpha_2),$$

where $F$ is a formal power series in two variables.

For $\mathfrak{I}_4$, the Hilbert basis is

| polynomial | degree | weight |
|---|---|---|
| $\alpha_1 = x_1$ | 1 | 3 |
| $\alpha_2 = x_2^2 - 2x_1x_3$ | 2 | 2 |
| $\alpha_3 = x_2^3 - 3x_1x_2x_3 + 3x_1^2x_4$ | 3 | 3 |
| $\alpha_4 = 9x_1^2x_4^2 - 3x_2^2x_3^2 - 18x_1x_2x_3x_4 + 8x_1x_3^3 + 6x_2^3x_4$ | 4 | 0 |

(4.10)

The block decomposition is

$$A_4 = \begin{bmatrix} \infty & \infty & 1 & \infty \\ 0 & 0 & 0 & 0 \end{bmatrix} = \begin{bmatrix} \infty & \infty & 0 & \infty \\ 0 & 0 & 0 & 0 \end{bmatrix} \sqcup \begin{bmatrix} \infty & \infty & 1 & \infty \\ 0 & 0 & 1 & 0 \end{bmatrix}. \tag{4.11}$$

The second version illustrates how a block decomposition can be rewritten using Stanley blocks, increasing the number of blocks used. The corresponding Stanley decomposition is

$$\mathcal{I}_4 = \mathbb{R}[[\alpha_1, \alpha_2, \alpha_4]] \oplus \mathbb{R}[[\alpha_1, \alpha_2, \alpha_4]]\alpha_3. \tag{4.12}$$

Every invariant $f$ can be written uniquely in the form

$$f(x_1, x_2, x_3, x_4) = F(\alpha_1, \alpha_2, \alpha_4) + G(\alpha_1, \alpha_2, \alpha_4)\alpha_3, \tag{4.13}$$

where $F$ and $G$ are formal power series in three variables. The number of formal functions required equals the number of Stanley spaces.

Since this is the first example that contains a relation, we use it to illustrate some points discussed in Section 3. The following relation holds for $\boldsymbol{\alpha}$:

$$\alpha_3^2 = \alpha_2^3 + \alpha_1^2\alpha_4. \tag{4.14}$$

By means of this relation, any occurrence of $\alpha_3$ to the second or higher power in an expression for $f \in \mathcal{I}_4$ can be eliminated; if this is not done, the expression for $f$ will not be unique. In fact the ideal of relations is $J = \langle t_3^2 - t_2^3 - t_1^2 t_4 \rangle$, and (for a particular choice of term order) the ideal of leading terms of $J$ is $I = \langle t_3^2 \rangle$. This makes clear why $\mathcal{I}_4 \cap \mathbb{R}[\mathbf{t}]$ is only isomorphic to $R[\mathbf{t}]/I$ as a vector space, not as an algebra (Remark 3.7); to reduce a polynomial modulo $I$ we need only erase terms in powers of $t_3$ beyond the first, whereas to reduce modulo $J$ (and preserve the product) we must replace $t_3^2$ by $t_2^3 + t_1^2 t_4$, and compute the replacements for higher powers accordingly.

The block decomposition and Stanley decomposition given above are not unique. By choosing a different term order in which $t_2^3$ is the leading term, one could, instead, eliminate powers of $\alpha_2$ higher than the second, obtaining a different block decomposition and Stanley decomposition. It would have three Stanley spaces and require three arbitrary formal functions.

For $\mathcal{I}_5$, the Hilbert basis is

| polynomial | degree | weight |
|---|---|---|
| $\alpha_1 = x_1$ | 1 | 4 |
| $\alpha_2 = x_2^2 - 2x_1x_3$ | 2 | 4 |
| $\alpha_3 = x_3^2 + 2x_1x_5 - 2x_2x_4$ | 2 | 0 |
| $\alpha_4 = x_2^3 - 3x_1x_2x_3 + 3x_1^2x_4$ | 3 | 6 |
| $\alpha_5 = 2x_3^3 - 12x_1x_3x_5 + 6x_2^2x_5$ | 3 | 0 |
| $\quad - 6x_2x_3x_4 + 9x_1x_4^2$ | | |

. (4.15)

There is a relation

$$\alpha_4^2 - \alpha_2^3 + 3\alpha_1\alpha_2\alpha_3 - \alpha_1^3\alpha_5 = 0.$$

A block decomposition is

$$A_5 = \begin{bmatrix} \infty & \infty & \infty & 1 & \infty \\ 0 & 0 & 0 & 0 & 0 \end{bmatrix} \tag{4.16}$$

and the Stanley decomposition is

$$\mathfrak{I}_5 = \mathbb{R}[[\alpha_1, \alpha_2, \alpha_3, \alpha_5]] \oplus \mathbb{R}[[\alpha_1, \alpha_2, \alpha_3, \alpha_5]]\alpha_4. \tag{4.17}$$

It is impractical to attempt a Stanley decomposition for $\mathfrak{I}_6$, since it is known [14] that there are 23 Hilbert basis elements and a large number of relations for the covariants of a quintic form (and therefore also for $\mathfrak{I}_6$).

## 5. Internal (semi)-transvectants

If $f$ and $g$ in $\mathfrak{I}$ are weight invariants with weights $\widehat{f}$ and $\widehat{g}$, the **internal transvectant** of $f$ and $g$ with **transvectant strength** $s$, denoted $(f, g)^{(s)}$, is defined for integers $s$ satisfying

$$0 \le s \le \widehat{f} \qquad \text{and} \qquad 0 \le s \le \widehat{g} \tag{5.1}$$

by the formula

$$(f, g)^{(s)} = \sum_{j=0}^{s} (-1)^j W_{fg}^{sj} (\mathfrak{Y}^j f)(\mathfrak{Y}^{s-j} g), \tag{5.2}$$

where

$$W_{fg}^{sj} = \binom{s}{j} \cdot \frac{(\widehat{f} - j)!}{(\widehat{f} - s)!} \cdot \frac{(\widehat{g} - s + j)!}{(\widehat{g} - s)!}. \tag{5.3}$$

**5.1. Remark.** In athletics, a person's strength limits how much weight the person can lift. With transvectants, it is the opposite: the two weights, $\widehat{f}$ and $\widehat{g}$, limit the strength by the inequalities (5.1). We say that the weights must each be large enough to **support** the strength.

**5.2. Remark.** The transvectant operations in classical invariant theory are defined for covariants. By the Roberts isomorphism (Remark 2.7), these can be transferred from covariants to seminvariants, and we have defined them that way. This transfer was done in [21, §6] and in greater detail in [29, §12.3]. More recently, it was done (independently) by Bedratyuk [2,3], who calls the transferred operation a *semi-transvectant*. We call the transvectants in this section *internal* because they take place *within* a single algebra $\mathfrak{I}$ of (semi)-invariants. Later we work with two algebras $\acute{\mathfrak{I}}$ and $\grave{\mathfrak{I}}$, and introduce transvectants *between* these. Those transvectants are called *external*.

**5.3. Theorem.** *If $f, g \in \mathfrak{I}$ are homogeneous weight invariants satisfying (5.1), then the internal transvectant $(f, g)^{(s)}$ is a homogeneous weight invariant in $\mathfrak{I}$ with degree and weight*

$$\deg(f, g)^{(s)} = \deg f + \deg g, \qquad \widehat{(f, g)^{(s)}} = \widehat{f} + \widehat{g} - 2s.$$

**Proof.** The degree of $(f,g)^{(s)}$ is self-evident, and the weight follows from Theorem 2.3, applied to the triad $(\mathcal{X},\mathcal{Y},\mathcal{Z})$ restricted to the finite-dimensional subspace $\mathbb{R}[\mathbf{x}]_{\deg f+\deg g}$). It remains to show that $\mathcal{X}(f,g)^{(s)}=0$. To do so, use

$$\begin{aligned}\mathcal{X}((\mathcal{Y}^i f)(\mathcal{Y}^{s-i} g)) &= \mathcal{X}(\mathcal{Y}^i f)(\mathcal{Y}^{s-i} g)+(\mathcal{Y}^i f)\mathcal{X}(\mathcal{Y}^{s-i} g)\\ &= i(\widehat{f}+1-i)(\mathcal{Y}^{i-1} f)(\mathcal{Y}^{s-i} g)\\ &\quad +(\mathcal{Y}^i f)(s-i)(\widehat{g}+1-s+i)(\mathcal{Y}^{s-i-1} g),\end{aligned} \tag{5.4}$$

which follows from the fact that $\mathcal{X}$ is a derivation, together with the pressure formula from Theorem 2.3. Thus, applying $\mathcal{X}$ to (5.2) doubles the number of terms from $s+1$ to $2(s+1)$; the first and last of these terms are zero, and the rest cancel in pairs, after the constants from (5.4) are multiplied by those from (5.2). □

# Part 2. Box products

## 6. A change of setting and notation

In Sections 2–5, we have studied invariants using the following set of notations:

a vector space $\mathbb{R}^n$ with variable $\mathbf{x}$;
a triad $(X,Y,Z)$ of $n\times n$ matrices acting on $\mathbb{R}^n$;
the algebra $\mathbb{R}[[\mathbf{x}]]$ of formal functions on $\mathbb{R}^n$;
an associated triad $(\mathcal{X},\mathcal{Y},\mathcal{Z})$ of operators on $\mathbb{R}[[\mathbf{x}]]$;
an algebra of invariants $\mathcal{I}=\ker\mathcal{X}\subset\mathbb{R}[[\mathbf{x}]]$;
a Hilbert basis $\boldsymbol{\alpha}=(\alpha_1,\dots,\alpha_p)$ for $\mathcal{I}$;
a Newton space $\mathbb{N}^p$ with variable $\mathbf{k}$, with $\boldsymbol{\alpha}^{\mathbf{k}}=\alpha_1^{k_1}\cdots\alpha_p^{k_p}$;
a preferred set $A$ of monomials in $\mathbb{R}[\boldsymbol{\alpha}]$;
a block decomposition $A=B^1\sqcup\cdots\sqcup B^\mu$ for $A$.

From here on, we will need three such settings at once, two for the *inputs* to the box product operation, and one for the *output*. The following table collects the notations that will be used, in the same order as the list above.

| left input | right input | output |
|---|---|---|
| $\mathbf{x}\in\mathbb{R}^n$ | $\mathbf{y}\in\mathbb{R}^m$ | $(\mathbf{x},\mathbf{y})\in\mathbb{R}^{n+m}$ |
| $(\acute{X},\acute{Y},\acute{Z})$ | $(\grave{X},\grave{Y},\grave{Z})$ | $(X,Y,Z)$ |
| $\mathbb{R}[[\mathbf{x}]]$ | $\mathbb{R}[[\mathbf{y}]]$ | $\mathbb{R}[[\mathbf{x},\mathbf{y}]]$ |
| $(\acute{\mathcal{X}},\acute{\mathcal{Y}},\acute{\mathcal{Z}})$ | $(\grave{\mathcal{X}},\grave{\mathcal{Y}},\grave{\mathcal{Z}})$ | $(\mathcal{X},\mathcal{Y},\mathcal{Z})$ |
| $\acute{\mathcal{I}}=\ker\acute{\mathcal{X}}$ | $\grave{\mathcal{I}}=\ker\grave{\mathcal{X}}$ | $\mathcal{I}=\ker\mathcal{X}$ |
| $\boldsymbol{\alpha}=(\alpha_1,\dots,\alpha_p)$ | $\boldsymbol{\beta}=(\beta_1,\dots,\beta_q)$ | $\boldsymbol{\vartheta}=(\vartheta_1,\dots,\vartheta_t)$ |
| $\mathbf{k}\in\mathbb{N}^p$ | $\boldsymbol{\ell}\in\mathbb{N}^q$ | $\mathbf{K}\in\mathbb{N}^t$ |
| $\boldsymbol{\alpha}^{\mathbf{k}}=\alpha_1^{k_1}\cdots\alpha_p^{k_p}$ | $\boldsymbol{\beta}^{\boldsymbol{\ell}}=\beta_1^{\ell_1}\cdots\beta_q^{\ell_q}$ | $\boldsymbol{\vartheta}^{\mathbf{K}}=\vartheta_1^{K_1}\cdots\vartheta_t^{K_t}$ |
| $\acute{A}\subset\mathbb{R}[\boldsymbol{\alpha}]$ | $\grave{A}\subset\mathbb{R}[\boldsymbol{\beta}]$ | $A=\acute{A}\boxtimes\grave{A}\subset\mathbb{R}[\boldsymbol{\vartheta}]$ |
| $\acute{A}=\acute{B}^1\sqcup\cdots\sqcup\acute{B}^\mu$ | $\grave{A}=\grave{B}^1\sqcup\cdots\sqcup\grave{B}^\nu$ | |

It is critical that the variables $x_i$ in $\mathbf{x}$ are independent of $y_j$ in $\mathbf{y}$; this is called the **externality** of the two inputs. Although accents are used to distinguish left inputs, right inputs, and outputs (as noted in Section 1), sometimes distinct letters are used instead. It would not be fair to inflict upon the reader an expression like $\acute{\mathbf{x}} = (\acute{x}_1, \dots, \acute{x}_{\acute{n}})$, with a tiny accent on a subscript.

Nothing more needs to be said about the input columns; these are givens in any box product calculation, and will come either from Section 4 or from the result of a previous box product calculation. However, most items in the output column need to be defined carefully in terms of the first two columns. We begin this now by defining $(X, Y, Z)$ in terms of $(\acute{X}, \acute{Y}, \acute{Z})$ and $(\grave{X}, \grave{Y}, \grave{Z})$ through the following formulas:

$$X = \begin{bmatrix} \acute{X} & 0 \\ 0 & \grave{X} \end{bmatrix}, \quad Y = \begin{bmatrix} \acute{Y} & 0 \\ 0 & \grave{Y} \end{bmatrix}, \quad Z = \begin{bmatrix} \acute{Z} & 0 \\ 0 & \grave{Z} \end{bmatrix}. \tag{6.1}$$

This notation has already been anticipated, in a simpler form, in (1.4). The triad $(\mathcal{X}, \mathcal{Y}, \mathcal{Z})$ is defined from $(X, Y, Z)$ in the usual way, as in (2.6), and then we put $\mathcal{I} = \ker \mathcal{X}$. The rest of the items in the output column will be defined gradually over the next few sections, at the point that they arise.

Because there is no overlap between the variables $\mathbf{x}$ and $\mathbf{y}$, the restriction of $\mathcal{X}$ to $\mathbb{R}[[\mathbf{x}]]$ equals $\acute{\mathcal{X}}$, and the restriction to $\mathbb{R}[[\mathbf{y}]]$ equals $\grave{\mathcal{X}}$:

$$\mathcal{X}|_{\mathbb{R}[[\mathbf{x}]]} = \acute{\mathcal{X}}, \quad \mathcal{X}|_{\mathbb{R}[[\mathbf{y}]]} = \grave{\mathcal{X}}, \tag{6.2}$$

and similarly for $\mathcal{Y}$ and $\mathcal{Z}$. Note that (6.1) together with (2.6) imply

$$\mathcal{X} = \acute{\mathcal{X}} + \grave{\mathcal{X}}, \quad \mathcal{Y} = \acute{\mathcal{Y}} + \grave{\mathcal{Y}}, \quad \mathcal{Z} = \acute{\mathcal{Z}} + \grave{\mathcal{Z}}. \tag{6.3}$$

For instance, if $n = 2$ and $m = 3$ with

$$(\acute{X}, \acute{Y}, \acute{Z}) = (M_2^*, N_2^*, H_2) \quad \text{and} \quad (\grave{X}, \grave{Y}, \grave{Z}) = (M_3^*, N_3^*, H_3),$$

then

$$\acute{\mathcal{X}} = \mathcal{D}_{N_2^*} = x_1 \partial/\partial x_2, \quad \grave{\mathcal{X}} = y_1 \partial/\partial y_2 + y_2 \partial/\partial y_3,$$

and

$$\mathcal{X} = x_1 \frac{\partial}{\partial x_2} + y_1 \frac{\partial}{\partial y_2} + y_2 \frac{\partial}{\partial y_3} = \acute{\mathcal{X}} + \grave{\mathcal{X}}.$$

**6.1. Remark.** With differential operators, it can be taken for granted that $\partial/\partial x_i$ is already defined on elements of $\mathbb{R}[[\mathbf{x}, \mathbf{y}]]$ and annihilates functions of $\mathbf{y}$ alone. In a more abstract treatment, $\acute{\mathcal{X}}$ would initially be defined only on $\mathbb{R}[[\mathbf{x}]]$, $\grave{\mathcal{X}}$ only on $\mathbb{R}[[\mathbf{y}]]$, and these would have to be extended in a suitable way to $\mathbb{R}[[\mathbf{x}, \mathbf{y}]]$ before (6.2) would make sense. The usual way of doing this is through tensor products. The map $\psi : \mathbb{R}[[\mathbf{x}]] \times \mathbb{R}[[\mathbf{y}]] \to \mathbb{R}[[\mathbf{x}, \mathbf{y}]]$ defined by $\psi(f, g) = fg$ meets the technical requirements for a tensor product map. (This is a consequence of externality.) This can be seen without using infinite-dimensional tensor products; the restrictions $\psi_{ij} : \mathbb{R}[\mathbf{x}]_i \times \mathbb{R}[\mathbf{y}]_j \to \mathbb{R}[\mathbf{x}, \mathbf{y}]_{ij}$ are tensor product maps of finite-dimensional spaces, where $\mathbb{R}[\mathbf{x}, \mathbf{y}]_{ij}$ is the space of polynomials in $\mathbf{x}$ and $\mathbf{y}$ that are homogeneous in $\mathbf{x}$ with degree $i$ and in $\mathbf{y}$ with degree $j$. That is, $\psi_{ij}$ induces a vector space isomorphism $\widetilde{\psi}_{ij} : \mathbb{R}[\mathbf{x}]_i \otimes \mathbb{R}[\mathbf{y}]_j \to \mathbb{R}[\mathbf{x}, \mathbf{y}]_{ij}$, so $\psi$ extends to an isomorphism $\widetilde{\psi}:\mathbb{R}[[\mathbf{x}]]\otimes\mathbb{R}[[\mathbf{y}]]\to\mathbb{R}[[\mathbf{x},\mathbf{y}]]$. The representation of $\mathfrak{sl}_2$ on $\mathbb{R}[[\mathbf{x},\mathbf{y}]]$ defined by $(\mathcal{X},\mathcal{Y},\mathcal{Z})$ is the tensor product of the representations defined by $(\acute{\mathcal{X}},\acute{\mathcal{Y}},\acute{\mathcal{Z}})$ and $(\grave{\mathcal{X}},\grave{\mathcal{Y}},\grave{\mathcal{Z}})$, and the equation (6.3) takes the form $\mathcal{X}=\acute{\mathcal{X}}\otimes I+I\otimes\grave{\mathcal{X}}$. See [21, §6], [26, §3], and [29, §12.3]. The first of these contains a minor error: on page 228, the condition on $\psi$ should be that its image spans $W$; it is not surjective.

## 7. External (semi)-transvectants

Let $f\in\acute{\mathcal{I}}$ and $g\in\grave{\mathcal{I}}$ be weight invariants (eigenvectors of $\acute{\mathcal{Z}}$ and $\grave{\mathcal{Z}}$ respectively) having weights $\widehat{f}$ and $\widehat{g}$. The **external transvectant** of $f$ and $g$ with **transvectant strength** $s$, denoted $(f,g)^{(s)}$, is defined for integers $s$ satisfying (5.1) by the formula

$$(f,g)^{(s)}=\sum_{j=0}^{s}(-1)^j W_{fg}^{sj}(\acute{\mathcal{Y}}^j f)(\grave{\mathcal{Y}}^{s-j} g), \tag{7.1}$$

where the coefficients are defined, as before, by (5.3). The first few external transvectants are given explicitly by

$$\begin{aligned}
(f,g)^{(0)}&=fg\\
(f,g)^{(1)}&=\widehat{f}f(\grave{\mathcal{Y}}g)-\widehat{g}(\acute{\mathcal{Y}}f)g\\
(f,g)^{(2)}&=\widehat{f}(\widehat{f}-1)f(\grave{\mathcal{Y}}^2g)-2(\widehat{f}-1)(\widehat{g}-1)(\acute{\mathcal{Y}}f)(\grave{\mathcal{Y}}g)\\
&\quad+\widehat{g}(\widehat{g}-1)(\acute{\mathcal{Y}}^2f)(g).
\end{aligned} \tag{7.2}$$

**7.1. Corollary.** *If $f\in\acute{\mathcal{I}}$, $g\in\grave{\mathcal{I}}$ are homogeneous weight invariants satisfying (5.1), then the external transvectant $(f,g)^{(s)}$ is a homogeneous weight invariant in $\mathcal{I}$ with degree and weight*

$$\deg(f,g)^{(s)}=\deg f+\deg g,\qquad \widehat{(f,g)^{(s)}}=\widehat{f}+\widehat{g}-2s.$$

**Proof.** In view of (6.2), (7.1) can be written in the form (5.2), as long as the restrictions $f\in\mathbb{R}[[\mathbf{x}]]$ and $g\in\mathbb{R}[[\mathbf{y}]]$ are observed. This shows that *external* transvectants between $\acute{\mathcal{I}}$ and $\grave{\mathcal{I}}$ form a subset of the *internal* transvectants within $\mathcal{I}$. Therefore Corollary 7.1 follows from Theorem 5.3. □

**7.2. Remark.** External transvectants perfectly express the Clebsch–Gordan theorem for tensor products of $\mathfrak{sl}_2$-representations; the external transvectants give a set of top weight vectors for the irreducible subrepresentations of the tensor product representation mentioned in Remark 6.1. Internal transvectants are subject to many identities, such as $(f,f)^{(s)}=0$ for $s$ odd, that often nullify their connection with Clebsch–Gordan. In particular, Theorem 8.4 below fails if there is any overlap between the variables $\mathbf{x}$ and $\mathbf{y}$. See the references listed in Remark 6.1 for more about tensor products, Clebsch–Gordan, and externality.

The next lemma follows immediately from (7.1) and will be useful later. There is no analog for general internal transvectants.

**7.3. Lemma.** *There is a nonzero constant c such that* $\acute{\mathfrak{X}}^s(f,g)^{(s)} = cfg$*. If* $k > s$, $\acute{\mathfrak{X}}^k(f,g)^{(s)} = 0$.

**Proof.** According to Theorem 2.3, $\acute{\mathfrak{X}}(\acute{\mathfrak{Y}}^i f) = p(\acute{\mathfrak{Y}}^{i-1} f)$ for some positive pressure factor $p$ if $i \geq 1$, and $= 0$ if $i = 0$. Apply $\acute{\mathfrak{X}}$ repeatedly to (7.1), using this fact and remembering that $\acute{\mathfrak{X}}$ treats functions of $\mathbf{y}$ as constants, and that $\acute{\mathfrak{X}} f = 0$ because $f \in \acute{\mathfrak{I}}$. After $s$ repetitions, all terms of (7.1) except the last ($j = s$) have been annihilated, and $s$ positive pressure factors $p_1, \dots, p_s$ have accumulated, so that $\acute{\mathfrak{X}}^s(f,g)^{(s)} = (-1)^s W^{ss}_{fg} p_1 \cdots p_s fg = cfg$. □

Lemma 7.3 remains true if $\acute{\mathfrak{X}}^s$ is replaced by $\grave{\mathfrak{X}}^s$.

## 8. The pre-box product

As outlined in Section 6, Hilbert bases for $\acute{\mathfrak{I}}$ and $\grave{\mathfrak{I}}$ are denoted by $\boldsymbol{\alpha} = (\alpha_1, \dots \alpha_p)$ and $\boldsymbol{\beta} = (\beta_1, \dots, \beta_q)$, and preferred sets of monomials in $\boldsymbol{\alpha}$ and in $\boldsymbol{\beta}$ are denoted by $\acute{A}$ and $\grave{A}$. The following types of transvectants, and products of transvectants, will play a major role in the remainder of this paper.

1. An **$\boldsymbol{\alpha\beta}$-transvectant** is an external transvectant $(f,g)^{(s)}$ in which $f = \boldsymbol{\alpha}^{\mathbf{k}}$ and $g = \boldsymbol{\beta}^{\boldsymbol{\ell}}$ for some $\mathbf{k}$ and $\boldsymbol{\ell}$. In view of (3.1) and (5.1), an $\boldsymbol{\alpha\beta}$-transvectant is well-defined if

$$0 \leq s \leq \widehat{\boldsymbol{\alpha}}\mathbf{k}, \qquad 0 \leq s \leq \widehat{\boldsymbol{\beta}}\boldsymbol{\ell}. \tag{8.1}$$

   The following short notation for $\boldsymbol{\alpha\beta}$-transvectants is sometimes convenient:

$$[\mathbf{k}; \boldsymbol{\ell}; s] = (\boldsymbol{\alpha}^{\mathbf{k}}, \boldsymbol{\beta}^{\boldsymbol{\ell}})^{(s)} \in \mathbb{R}[\mathbf{x}, \mathbf{y}]. \tag{8.2}$$

2. An **admissible transvectant** with respect to preferred sets $\acute{A}$ and $\grave{A}$ is an $\boldsymbol{\alpha\beta}$-transvectant such that $\boldsymbol{\alpha}^{\mathbf{k}} \in \acute{A}$ and $\boldsymbol{\beta}^{\boldsymbol{\ell}} \in \grave{A}$. For the simplest example of a transvectant that is well-defined but not admissible, see Remark 11.6.
3. A **transvectant monomial** is a product of $\boldsymbol{\alpha\beta}$-transvectants (which of course includes powers of transvectants). This product takes place in $\mathbb{R}[[\mathbf{x}, \mathbf{y}]]$ and requires no special definition. An **admissible transvectant monomial** is a product of admissible transvectants.

**8.1. Remark.** Although the only requirement for a transvectant $(\boldsymbol{\alpha}^{\mathbf{k}}, \boldsymbol{\beta}^{\boldsymbol{\ell}})^{(s)}$ to be *admissible* is that $\boldsymbol{\alpha}^{\mathbf{k}}$ and $\boldsymbol{\beta}^{\boldsymbol{\ell}}$ be *preferred*, it would be confusing to call an admissible transvectant a preferred transvectant. Instead, a *subset* of the admissible *transvectant monomials* will be singled out later as *preferred* monomials, in the sense of Section 3.

The **pre-box product** $\acute{A} \,\widehat{\boxtimes}\, \grave{A}$ of $\acute{A}$ and $\grave{A}$ is defined to be the set of all admissible external transvectants between $\acute{A}$ and $\grave{A}$:

$$\acute{A} \,\widehat{\boxtimes}\, \grave{A} = \left\{ [\mathbf{k}; \boldsymbol{\ell}; s] : \boldsymbol{\alpha}^{\mathbf{k}} \in \acute{A},\ \boldsymbol{\beta}^{\boldsymbol{\ell}} \in \grave{A},\ 0 \leq s \leq \widehat{\boldsymbol{\alpha}}\mathbf{k},\ 0 \leq s \leq \widehat{\boldsymbol{\beta}}\boldsymbol{\ell} \right\}. \tag{8.3}$$

This definition is generalized to *subsets* of $\acute{A}$ and $\grave{A}$ in (8.6).

The next task is to show that the pre-box product (8.3) provides a formal vector space basis for $\mathfrak{I}$. This will be proved in Theorem 8.4, after two lemmas.

**8.2. Lemma.** *The set of products $\boldsymbol{\alpha}^{\mathbf{k}}\boldsymbol{\beta}^{\boldsymbol{\ell}}$, for $\boldsymbol{\alpha}^{\mathbf{k}} \in \acute{A}$ and $\boldsymbol{\beta}^{\boldsymbol{\ell}} \in \grave{A}$, is linearly independent.*

**Proof.** The set of $f = \boldsymbol{\alpha}^{\mathbf{k}} \in \acute{A}$ is linearly independent by the definition of preferred set, and similarly for the $g = \boldsymbol{\beta}^{\boldsymbol{\ell}} \in \grave{A}$. To see that the set of products $fg$ is linearly independent, suppose some nontrivial finite linear combination $c_1 f_1 g_1 + \cdots + c_j f_j g_j$ is identically zero, with none of the $c_i$ being zero (since any such terms could be dropped). Since the functions $f_1, \dots, f_j$ are linearly independent, none of them are identically zero. Choose a point $\mathbf{x}$ for which the $f_i(\mathbf{x})$ are not all zero. Then $c_1 f_1(\mathbf{x}) g_1 + \cdots + c_j f_j(\mathbf{x}) g_j$ is a nontrivial linear combination of $g_1, \dots, g_j$ that is identically zero. But this is impossible, since $g_1, \dots, g_j$ are linearly independent. □

**8.3. Lemma.** *The set $\acute{A} \,\widehat{\boxtimes}\, \grave{A}$ of admissible transvectants is linearly independent.*

**Proof.** In view of Remark 3.1, it suffices to prove that for each $d$ and $w$, the transvectants $(f, g)^{(s)} \in \mathfrak{I}_{dw}$ are linearly independent, where we write $f$ and $g$ for $\boldsymbol{\alpha}^{\mathbf{k}}$ and $\boldsymbol{\beta}^{\boldsymbol{\ell}}$. Suppose this is not true. Then for some $k > 0$ there exist $f_i, g_i, s_i \in \mathfrak{I}_{dw}$ for $i = 1, \dots, k$ such that

$$c_1 (f_1, g_1)^{(s_1)} + \cdots + c_k (f_k, g_k)^{(s_k)} = 0, \tag{8.4}$$

where none of the $c_i$ are zero and no transvectants are repeated (since such terms could be combined). Non-repetition means that there do not exist $i \neq j$ in (8.4) such that $f_i = f_j$, $g_i = g_j$, and $s_i = s_j$. But indeed, there do not even exist $i \neq j$ in (8.4) such that $f_i = f_j$ and $g_i = g_j$; if there did, it would follow from what we just said that $s_i \neq s_j$, but then $(f_i, g_j)^{(s_i)}$ and $(f_j, g_j)^{(s_j)}$ would have different weights and would not belong to the same space $\mathfrak{I}_{dw}$.

Let $s$ be the largest value of $s_i$ (which may occur more than once), and apply $\acute{\mathfrak{X}}^s$ to (8.4). By Theorem 7.3, this will produce

$$c_1' f_1 g_1 + \cdots + c_k' f_k g_k = 0, \tag{8.5}$$

where $c_i' = 0$ if $s_i < s$ and $c_i' \neq 0$ if $s_i = s$. By the remarks in the last paragraph, none of the products $f_i g_i$ in (8.5) will be the same. So (8.5) is a nontrivial linear relation among a set of distinct products of the form $\boldsymbol{\alpha}^{\mathbf{k}}\boldsymbol{\beta}^{\boldsymbol{\ell}}$, which is impossible by Lemma 8.2. □

**8.4. Theorem.** *The transvectants in $\acute{A} \,\widehat{\boxtimes}\, \grave{A}$ form a vector space basis for the polynomial invariants $\mathfrak{I} \cap \mathbb{R}[\mathbf{x}, \mathbf{y}]$, and a formal vector space basis for the formal invariants $\mathfrak{I}$.*

**Proof.** By Remark 2.6, applied to $\mathbb{R}[[\mathbf{x}, \mathbf{y}]]$, it suffices to prove the polynomial case. By Corollary 7.1, the transvectants in $\acute{A} \,\widehat{\boxtimes}\, \grave{A}$ lie in $\mathfrak{I} \cap \mathbb{R}[\mathbf{x}, \mathbf{y}]$, and by Lemma 8.3, these are linearly independent. It remains to prove that they span $\mathfrak{I} \cap \mathbb{R}[\mathbf{x}, \mathbf{y}]$. We show that the transvectants of degree $d$ span $\mathfrak{I} \cap R[\mathbf{x}, \mathbf{y}]_d$, by showing that the chains having these transvectants as top weight vectors span $\mathbb{R}[\mathbf{x}, \mathbf{y}]_d$, so that there is no room for any additional top weight vectors besides the transvectants.

For this paragraph, consider a *fixed* pair $f = \boldsymbol{\alpha}^{\mathbf{k}} \in \acute{A}$, $g = \boldsymbol{\beta}^{\boldsymbol{\ell}} \in \grave{A}$, such that $\deg f + \deg g = d$. The lengths of the chains headed by $f$ and $g$ are $\widehat{f} + 1$ and $\widehat{g} + 1$ respectively. Let $V \subset \mathbb{R}[\mathbf{x}, \mathbf{y}]$ be the span of the products $(\acute{\mathfrak{Y}}^i f)(\grave{\mathfrak{Y}}^j g)$ for $i = 0, \dots, \widehat{f}$, $j = 0, \dots, \widehat{g}$. This space of dimension $(\widehat{f} + 1)(\widehat{g} + 1)$ is closed under $\{\mathfrak{X}, \mathfrak{Y}, \mathfrak{Z}\}$ and is therefore a representation space of $\mathfrak{sl}_2$. By (7.1), the well-defined transvectants $(f, g)^{(s)}$ belong to $V$. The strengths $s$ for which the transvectant is well-defined will be $s = 0, 1, \dots, \min\{\widehat{f}, \widehat{g}\}$. Consider cases $\widehat{g} \leq \widehat{f}$ and $\widehat{f} < \widehat{g}$. If $\widehat{g} \leq \widehat{f}$ then $s = 0, 1, \dots, \widehat{g}$, the transvectants have weights $(\widehat{f}+\widehat{g}), (\widehat{f}+\widehat{g}-2), \dots, (\widehat{f}-\widehat{g}-2\widehat{g}) = (\widehat{f}-\widehat{g})$, and the chains headed by these transvectants have lengths $(\widehat{f}+\widehat{g}+1), \dots, (\widehat{f}-\widehat{g}+1)$. This is an arithmetic series having $\widehat{g}$ terms, with sum

$$\frac{\widehat{g}}{2} \cdot \big((\widehat{f}+\widehat{g}+1) + (\widehat{f}-\widehat{g}+1)\big) = (\widehat{f}+1)(\widehat{g}+1).$$

So the element of the chains headed by these transvectants are equal in number to the dimension of $V$. Since the elements of these chains are linearly independent by Theorem 2.3, they span $V$. The other case ($\widehat{f} < \widehat{g}$) gives the same result.

Summing over all such pairs $f, g$, the chains headed by the transvectants span the direct sum of the subrepresentations $V$. This direct sum is $\mathbb{R}[\mathbf{x}, \mathbf{y}]_d$, so the chains span $\mathbb{R}[\mathbf{x}, \mathbf{y}]$ as claimed. □

**8.5. Remark.** Theorem 8.4 is a form of the Clebsch–Gordan theorem, which is often stated only as a count of the number of chains in the tensor product of two irreducible representations of $\mathfrak{sl}_2$. (The number is equal to the number of elements in the irreducible representation with the smaller chain.) That version of the theorem does not give the transvectant formula for the top weight vectors of the chains.

Since $\acute{A}$ and $\grave{A}$ are usually given as disjoint unions of blocks, it is helpful to extend the definition of pre-box product to subsets $\acute{M} \subset \acute{A}$ and $\grave{M} \subset \grave{A}$ as follows:

$$\acute{M} \mathbin{\widehat{\boxtimes}} \grave{M} = \left\{ [\mathbf{k}; \boldsymbol{\ell}; s] : \boldsymbol{\alpha}^{\mathbf{k}} \in \acute{M},\ \boldsymbol{\beta}^{\boldsymbol{\ell}} \in \grave{M},\ 0 \le s \le \widehat{\boldsymbol{\alpha}}\mathbf{k},\ 0 \le s \le \widehat{\boldsymbol{\beta}}\boldsymbol{\ell} \right\}. \tag{8.6}$$

This (extended) pre-box product now distributes over disjoint unions of subsets of $\acute{A}$ and $\grave{A}$; therefore if $\acute{A}$ and $\grave{A}$ are given by block decompositions

$$\acute{A} = \acute{B}^1 \sqcup \cdots \sqcup \acute{B}^{\mu}, \qquad \grave{A} = \grave{B}^1 \sqcup \dots \grave{B}^{\nu}, \tag{8.7}$$

then

$$\acute{A} \mathbin{\widehat{\boxtimes}} \grave{A} = \bigsqcup_{i=1}^{\mu} \bigsqcup_{j=1}^{\nu} \acute{B}^i \mathbin{\widehat{\boxtimes}} \grave{B}^j. \tag{8.8}$$

## 9. The replacement theorem

The pre-box product $\acute{A} \mathbin{\widehat{\boxtimes}} \grave{A}$ provides one formal vector space basis for $\mathfrak{I}$. The replacement theorem (originally proved in [26, Thm. 3]) states that other such bases can be obtained by **replacing** individual transvectants in $\acute{A} \mathbin{\widehat{\boxtimes}} \grave{A}$ by suitable transvectant monomials. If this is done skillfully, the resulting basis can have advantageous properties.

A transvectant monomial

$$[\mathbf{k}_1; \boldsymbol{\ell}_1; s_1]^{m_1} \cdots [\mathbf{k}_u; \boldsymbol{\ell}_u; s_u]^{m_u} \tag{9.1}$$

will be called a **replacement** for the transvectant $[\mathbf{k}; \boldsymbol{\ell}; s]$ if the following linear relation among integer vectors holds:

$$(\mathbf{k}; \boldsymbol{\ell}; s) = m_1(\mathbf{k}_1; \boldsymbol{\ell}_1; s_1) + \cdots + m_u(\mathbf{k}_u; \boldsymbol{\ell}_u; s_u). \tag{9.2}$$

**9.1. Lemma.** *Any replacement for an admissible transvectant is an admissible transvectant monomial.*

**Proof.** We must prove that if $[\mathbf{k}; \boldsymbol{\ell}; s]$ is admissible, then any transvectant $[\mathbf{k}_i; \boldsymbol{\ell}_i; s_i]$ occurring in (9.1) is admissible. Since $\boldsymbol{\alpha}^{\mathbf{k}_i}$ divides $\boldsymbol{\alpha}^{\mathbf{k}} \in \acute{A}$, and $\acute{A}$ is standard, Lemma 3.6 implies that $\boldsymbol{\alpha}^{\mathbf{k}_i} \in \acute{A}$. Similarly, $\boldsymbol{\beta}^{\boldsymbol{\ell}_i} \in \grave{A}$. □

If some or all of the transvectants in $\acute{A} \mathbin{\widehat{\boxtimes}} \grave{A}$ are replaced in this way, the resulting set is called a **replacement set** for $\acute{A} \mathbin{\widehat{\boxtimes}} \grave{A}$. (We may assume that *all* elements of $\acute{A} \mathbin{\widehat{\boxtimes}} \grave{A}$ have been replaced, since each transvectant *is* a transvectant monomial and can be "replaced by itself".)

**9.2. Theorem** *(Replacement theorem). Any replacement set for $\acute{A} \mathbin{\widehat{\boxtimes}} \grave{A}$ is itself a formal vector space basis for $\mathfrak{I}$.*

**Proof.** First we claim that the replacement set for $\acute{A} \mathbin{\widehat{\boxtimes}} \grave{A}$ is linearly independent. The proof parallels that of Lemma 8.3. Since each transvectant $[\mathbf{k}; \boldsymbol{\ell}; s]$ in $\acute{A} \mathbin{\widehat{\boxtimes}} \grave{A}$ is replaced by a transvectant monomial having the same degree and weight, it remains true that no two elements of the replacement set belonging to the same subspace $\mathfrak{I}_{dw}$ can have the same values of $\mathbf{k}$ and $\boldsymbol{\ell}$ (as determined from (9.2)). Next, it is necessary to strengthen Lemma 7.3 so that it applies to transvectant monomials. Let $f_i = \boldsymbol{\alpha}^{\mathbf{k}_i}$ and $g_i = \boldsymbol{\beta}^{\boldsymbol{\ell}_i}$. When the transvectants in (9.1) are written out using the definition (7.1), and these expressions are multiplied together, there will be one term that contains $s$ factors of $\acute{\mathcal{Y}}$ and none of $\grave{\mathcal{Y}}$, where $s = m_1 s_1 + \cdots + m_u s_u$. This term will equal a constant times

$$(\acute{\mathcal{Y}}^{s_1} f_1)^{m_1} \cdots (\acute{\mathcal{Y}}^{s_u} f_u)^{m_u} g_1^{m_1} \cdots g_u^{m_u}.$$

Then, as in Lemma 7.3, $\acute{\mathcal{X}}^s$ applied to this term produces

$$c f_1^{m_1} \cdots f_u^{m_u} g_1^{m_1} \cdots g_u^{m_u},$$

where $c \neq 0$. All other terms have fewer than $s$ factors of $\acute{\mathcal{Y}}$, and are annihilated by $\acute{\mathcal{X}}^s$. The proof of linear independence now goes through as before.

Since the elements of the replacement belonging to any $\mathfrak{I}_{dw}$ are linearly independent and have the same cardinality as the basis for $\mathfrak{I}_{dw}$ given by the pre-box product, they form a basis for $\mathfrak{I}_{dw}$, and the full replacement set forms a formal basis for $\mathfrak{I}$. □

A transvectant is not, in general, equal (when expressed as a polynomial in $\mathbf{x}$ and $\mathbf{y}$) to any of its replacements (other than itself). We say instead that a transvectant is **equivalent** to each of its replacements, and that all the replacements of a given transvectant are equivalent to each other; equivalence is denoted by $\equiv$, and is an equivalence relation on the set of admissible transvectant monomials. Given a transvectant monomial written as in (9.1), it is easy to find the unique transvectant $[\mathbf{k}; \boldsymbol{\ell}; s]$ to which it is equivalent, by using (9.2).

The opposite direction is more difficult. Given a transvectant $[\mathbf{k}; \boldsymbol{\ell}; s]$ in $\acute{A} \mathbin{\widehat{\boxtimes}} \grave{A}$ that we would like to replace, the idea is to split $\mathbf{k}$, $\boldsymbol{\ell}$, and $s$ into two parts, so that $\mathbf{k} = \mathbf{k}' + \mathbf{k}''$, $\boldsymbol{\ell} = \boldsymbol{\ell}' + \boldsymbol{\ell}''$, $s = s' + s''$; this must be done in such a way that the two transvectants $[\mathbf{k}'; \boldsymbol{\ell}'; s']$ and $[\mathbf{k}''; \boldsymbol{\ell}''; s'']$ are both well-defined. If this can be accomplished, we write

$$[\mathbf{k};\boldsymbol{\ell};s]\equiv[\mathbf{k}';\boldsymbol{\ell}';s'][\mathbf{k}'';\boldsymbol{\ell}'';s''], \tag{9.3}$$

and say that we have **factored the transvectant** (up to equivalence). Each factor may be factorable again, until the only remaining transvectants are *prime* (cannot be factored again).

**9.3. Remark.** Later the following related definition will be useful. Let $\delta=[\mathbf{k};\boldsymbol{\ell};s]$ and $\vartheta=[\boldsymbol{\kappa};\boldsymbol{\lambda};\sigma]$ be well-defined transvectants. We say that $\delta$ is **divisible** by $\vartheta$ if $\vartheta$ can be factored from $\delta$ (up to equivalence); this will be the case if and only if $[\mathbf{k}-\boldsymbol{\kappa};\boldsymbol{\ell}-\boldsymbol{\lambda};s-\sigma]$ is well-defined.

In practice it is usually best to use transvectant notation rather than the short notation $[\mathbf{k};\boldsymbol{\ell};s]$. For instance, if $p=1$ and $q=2$ with $\widehat{\alpha}_1=1$, $\widehat{\beta}_1=2$, and $\widehat{\beta}_2=0$, then

$$(\alpha_1^3,\beta_1^2)^{(3)}\equiv(\alpha_1,\beta_1)^{(1)}(\alpha_1^2,\beta_1)^{(2)}$$

is clearly a valid factorization; the corresponding expression

$$[3;2,0;3]\equiv[1;1,0;1][2;1,0;2]$$

is much less clear. It is easier to check well-definedness in the first version, because we associate the weight of a transvectant with its name more easily than with a position in an integer vector. (For a computer, this would be no problem.) Also, the "short" notation forces us to mention $\beta_2$ (by listing its exponent as zero), although it does not appear at all in the "long" notation. Note that the factorization

$$(\alpha_1^3,\beta_1^2)^{(3)}\equiv(\alpha_1,\beta_1)^{(2)}(\alpha_1^2,\beta_1)^{(1)} \qquad \text{(FALSE)}$$

is not valid, even though the values of $\mathbf{k}$, $\boldsymbol{\ell}$, and $s$ add up correctly; $(\alpha_1,\beta_1)^{(2)}$ is not defined, since $\alpha_1$, with weight 1, cannot support strength 2 (see Remark 5.1). The factorization

$$(\alpha_1^2,\beta_2)^{(1)}\equiv\alpha_1(\alpha_1,\beta_2)^{(1)}$$

is valid, since the first factor $\alpha_1$ can be viewed as the zero-transvectant $[1;0,0;0]=(\alpha_1,\beta_1^0\beta_2^0)^{(0)}=(\alpha_1,1)^{(0)}$.

## 10. Box products by inspection, and a note on the expansion method

The pre-box product provides only a formal vector space basis for $\mathcal{I}$, not a Hilbert basis. Furthermore, to compute the pre-box product explicitly (as a set of elements of $\mathbb{R}[\mathbf{x},\mathbf{y}]$) would require the computation of infinitely many transvectants, using the formula (7.1); this cannot be completed in finite time. We now show, in a couple of simple cases, that the replacement theorem, used skillfully, can produce a formal vector space basis for $\mathcal{I}$ that requires only finitely many transvectants, and that this finite set of transvectants constitutes a Hilbert basis for $\mathcal{I}$. The factorizations used for this purpose also reveal a Stanley decomposition for $\mathcal{I}$ in terms of this Hilbert basis; this solves the description problem for $\mathcal{I}$. The *ad hoc* method used for these examples does not extend easily to more complicated cases, and will be replaced by an algorithmic factoring method in Section 15.

**10.1. Example.** $N_{22}$

When $\acute{N} = \grave{N} = N_2$, the input data for both factors of $\acute{A}_2 \boxtimes \grave{A}_2$ is given by (4.4) and (4.5), with $\alpha_1$ changed to $\beta_1$ and $\mathbf{x}$ to $\mathbf{y}$ for the second factor. That is, $\alpha_1 = x_1$ and $\beta_1 = y_1$. Since $\widehat{\alpha}_1 = \widehat{\beta}_2 = 1$,

$$\begin{bmatrix}\infty\\0\end{bmatrix} \widehat{\boxtimes} \begin{bmatrix}\infty\\0\end{bmatrix} = \left\{(\alpha_1^{k_1}, \beta_1^{\ell_1})^{(s)} : 0 \le s \le k_1,\ 0 \le s \le \ell_1\right\}.$$

Each of these transvectants may be factored by inspection as

$$\alpha_1^{k_1-s}\beta_1^{\ell_1-s}(\alpha_1, \beta_1)^{(1)s}. \tag{10.1}$$

According to the replacement theorem, the set of products (10.1) constitutes a formal vector space basis for $\mathcal{I}_{22}$. This basis requires only the finite set of transvectants $\{\alpha_1, \beta_1, \gamma_1\}$, where $\gamma_1 = (\alpha_1, \beta_1)^{(1)}$. Every invariant is a countable linear combination of the invariants in (10.1), and is therefore a formal power series in $(\alpha_1, \beta_1, \gamma_1)$, so this ordered set is a Hilbert basis for $\mathcal{I}_{22}$. Since $k_1$ and $\ell_1$ can be arbitrarily large, there are no restrictions on the powers that appear in (10.1), and all monomials in $(\alpha_1, \beta_1, \gamma_1)$ are preferred monomials for $\mathcal{I}$ in terms of this Hilbert basis. The block decomposition for $\mathcal{I}_{22}$ is the bottom two rows of the "decorated" b-matrix (see Remark 3.4)

$$\left[\begin{array}{ccc} 1 & 1 & 2\\ 1 & 1 & 0\\ \alpha_1 & \beta_1 & \gamma_1\\ \hline \infty & \infty & \infty\\ 0 & 0 & 0\end{array}\right]; \tag{10.2}$$

the first three rows give the degrees, weight, and names of the Hilbert basis elements, collected in one place for convenient reference. The Stanley decomposition is

$$\mathcal{I}_{22} = \mathbb{R}[[\alpha_1, \beta_1, \gamma_1]] = \mathbb{R}[[x_1, y_1, x_1y_2 - x_2y_1]]. \tag{10.3}$$

For completeness we have computed $\gamma_1 = x_1y_2 - x_2y_1$ from (7.2), a step that we usually omit, since it is best done by computer algebra and is mathematically trivial, although it is necessary in order to finalize the result for applications. The Stanley decomposition states that every invariant $f$ can be written uniquely as $f(x_1, x_2, y_1, y_2) = F(\alpha_1, \beta_1, \gamma_1)$, where $F$ is a formal power series in three variables. □

**10.2. Example.** $N_{23}$, first time.

The input data is given by (4.4), (4.5), (4.7), and (4.8), with the appropriate change to $\boldsymbol{\beta}$ and $\mathbf{y}$ for the second factor. The admissible transvectants are

$$\acute{A} \,\widehat{\boxtimes}\, \grave{A} = \{(\alpha_1^{k_1}, \beta_1^{\ell_1}\beta_2^{\ell_2})^{(s)} : 0 \le s \le k_1,\ 0 \le s \le 2\ell_1 + 0\ell_2\}. \tag{10.4}$$

By an *ad hoc* decision, we divide these transvectants into two classes according to whether $s$ is even ($s = 2p$) or odd ($s = 2p + 1$). Introducing arbitrary nonnegative slack variables $i$ and $j$, in the even case the general transvectant is

$$(\alpha_1^{2p+i}, \beta_1^{p+j}\beta_2^{\ell_2})^{(2p)} \equiv \alpha_1^i \beta_1^j \beta_2^{\ell_2} (\alpha_1^2, \beta_1)^{(2)p},$$

while in the odd case

$$(\alpha_1^{2p+1+i}, \beta_1^{p+1+j}\beta_2^{\ell_2})^{(2p+1)} \equiv \alpha_1^i \beta_1^j \beta_2^{\ell_2} (\alpha_1^2, \beta_1)^{(2)p} (\alpha_1, \beta_1)^{(1)}.$$

These two cases give a complete set of replacements. Writing

$$\gamma_1 = (\alpha_1, \beta_1)^{(1)}, \qquad \gamma_2 = (\alpha_1^2, \beta_1)^{(2)}, \tag{10.5}$$

the block decomposition is

$$\left[\begin{array}{c|cc|cc} 1 & 1 & 2 & 2 & 3 \\ 1 & 2 & 0 & 1 & 0 \\ \alpha_1 & \beta_1 & \beta_2 & \gamma_1 & \gamma_2 \\ \hline \infty & \infty & \infty & 1 & \infty \\ 0 & 0 & 0 & 0 & 0 \end{array}\right] \tag{10.6}$$

and the Stanley decomposition is

$$\mathfrak{I}_{23} = \mathbb{R}[[\alpha_1, \beta_1, \beta_2, \gamma_2]] \oplus \mathbb{R}[[\alpha_1, \beta_1, \beta_2, \gamma_2]]\gamma_1. \tag{10.7}$$

The vertical lines in (10.6) separate the $\boldsymbol{\alpha}$, $\boldsymbol{\beta}$, and $\boldsymbol{\gamma}$ portions of the Hilbert basis, for easier reading when the top three "decoration" rows are omitted. Every invariant can be written uniquely as

$$f(x_1, x_2, y_1, y_2, y_3) = F(\alpha_1, \beta_1, \beta_2, \gamma_2) + G(\alpha_1, \beta_1, \beta_2, \gamma_2)\gamma_1,$$

where $F$ and $G$ are formal power series in four variables. □

**10.3. Remark.** The **expansion method** for computing box products, presented in [26] and [29, Ch. 12] begins by assuming that $\acute{A}$ and $\grave{A}$ are given by Stanley decompositions. These Stanley decompositions are repeatedly expanded into longer forms by a procedure that we will not describe here, except to say that each expansion terminates when the least common multiple of two weights is reached. From our present point of view, the purpose of the expansions is to obtain subsets of the pre-box product that can be factored in a natural way, as we have done above by inspection. There are several drawbacks to this method. First, the expansion of the original Stanley decompositions leads to a large number of subsets of the pre-box product, and therefore to a final Stanley decomposition for $\mathfrak{I}$ that can be much longer than necessary. These can be compressed, but there is no algorithm to produce the most effective compression; an initial choice to combine two or more terms can prevent another combination that would have led to a shorter result. The compressions that were done in the references above and in [30] were done by inspection, rather than by any rule. Second, the expansion method can produce preferred sets for $\mathfrak{I}$ that are not standard; see Remark 3.8. These Stanley decompositions are not invalid, but they are more complicated than necessary (in ways that probably cannot be removed by compression). Third, it would be desirable to extend the box product method from the Lie algebra $\mathfrak{sl}_2$ to its $k$-fold Cartesian product $\mathfrak{sl}_2^k$ for integers $k > 1$. The associated Lie group is known in quantum computing, entanglement, and information theory as SLOCC ($k$) ("stochastic local operations with classical communication"). The expansion method cannot be extended to this case, because (with $k > 1$) the weights are integer vectors with $k$ entries, and pairs of these weights need not have least common multiples (which the expansion method requires). The factoring method, however, can be extended to this case. We will not do so here, because it is not relevant to normal forms of dynamical systems.

In our joint work on the expansion method, Jan Sanders and I preferred different notations and emphases; the version in [26] is mine, and that in [29, Ch. 12] is his. The following remarks might be helpful in comparing these. There is a routine procedure to produce a two-variable Hilbert function or "table function" from a Stanley decomposition for $\mathcal{I}$; see [22, p. 256]. This is the generating function for the degrees and weights of the invariants, and it loses information contained in the Stanley decomposition. Sanders [29, p. 291] created a *subscripted* version of the Hilbert function that retains the full information in the Stanley decomposition, and he calculates in a formal way with these functions. I work with the Stanley decomposition directly, and prove the validity of the procedure by a double induction. Sanders compares the two notations in an example on page 305.

## 11. The factoring method: theory

Assume that the content of the two "input" columns of the table in Section 6 is given (and that the preferred sets $\acute{A}$ and $\grave{A}$ are standard). Let $\mathcal{C} = \mathcal{C}(\boldsymbol{\alpha}, \boldsymbol{\beta}) \subset \mathbb{N}^{p+q+1}$ be the set of integer vectors $(\mathbf{k}; \boldsymbol{\ell}; s) \in \mathbb{N}^{p+q+1}$ satisfying (8.1), so that the transvectant $[\mathbf{k}; \boldsymbol{\ell}; s]$ is well-defined. The set $\mathcal{C}$ will be called the **transvectant cone**, because it satisfies the definition of a cone or monoid in integer programming: it contains the origin of $\mathbb{N}^{p+q+1}$ and is closed under vector addition (see [31]) or [33]). It follows that $\mathcal{C}$ is closed under nonnegative integer linear combinations. Let $\mathcal{C}^*$ be the subset of $(\mathbf{k}, \boldsymbol{\ell}, s) \in \mathcal{C}$ such that $\boldsymbol{\alpha}^{\mathbf{k}} \in \acute{A}$ and $\boldsymbol{\beta}^{\boldsymbol{\ell}} \in \grave{A}$, that is, such that $[\mathbf{k}; \boldsymbol{\ell}; s]$ is admissible with respect to $\acute{A}$ and $\grave{A}$. Define the **transvectant map** $\mathbb{T} : \mathcal{C} \to \mathbb{R}[\mathbf{x}, \mathbf{y}]$ by

$$\mathbb{T}(\mathbf{k}; \boldsymbol{\ell}; s) = [\mathbf{k}; \boldsymbol{\ell}; s] = (\boldsymbol{\alpha}^{\mathbf{k}}, \boldsymbol{\beta}^{\boldsymbol{\ell}})^{(s)}. \tag{11.1}$$

**11.1. Lemma.** *The transvectant map, restricted to $\mathcal{C}^*$, is injective, and its (partial) inverse is computable.*

**Proof.** Suppose there exist distinct points $(\mathbf{k}; \boldsymbol{\ell}; s)$ and $(\mathbf{k}'; \boldsymbol{\ell}'; s')$ in $\mathcal{C}$ such that $\mathbb{T}(\mathbf{k}; \boldsymbol{\ell}; s) = \mathbb{T}(\mathbf{k}'; \boldsymbol{\ell}'; s')$ in $\mathbb{R}[[\mathbf{x}, \mathbf{y}]]$. Then $[\mathbf{k}; \boldsymbol{\ell}; s] - [\mathbf{k}'; \boldsymbol{\ell}'; s'] = 0$. This equation has the form (8.4), and leads to a contradiction by the same reasoning as in the proof of Lemma 8.3. (Note that the *proof* of the lemma is needed, not the *statement*. The reason that the statement is insufficient is that set theory does not allow repeated elements in a set. One would like to argue that if $[\boldsymbol{\alpha}, \boldsymbol{\ell}, s] = [\boldsymbol{\alpha}', \boldsymbol{\ell}', s']$ then the set $\{[\boldsymbol{\alpha}; \boldsymbol{\ell}; s], [\boldsymbol{\alpha}'; , \boldsymbol{\ell}'; s']\}$ is linearly dependent, in contradiction to the lemma, so the supposition must be false. But in fact, in this case $\{[\boldsymbol{\alpha}; \boldsymbol{\ell}; s], [\boldsymbol{\alpha}'; , \boldsymbol{\ell}'; s']\}$ contains only one element, not two, and the set is linearly independent, so the lemma is not contradicted.)

For the inverse, the issue is as follows. Each transvectant is a homogeneous polynomial belonging to $\mathbb{R}[\mathbf{x}, \mathbf{y}]_d$ for some $d$. We must show that given any $h \in \mathbb{R}[\mathbf{x}, \mathbf{y}]_d$, it is possible (at least in principle) to determine whether $h$ is an admissible $\boldsymbol{\alpha\beta}$-transvectant, and if so, to recover (uniquely) the values of $\mathbf{k}$, $\boldsymbol{\ell}$, and $s$ such that $h = \mathbb{T}(\mathbf{k}; \boldsymbol{\ell}; s)$. Here is a simple proof: Since there are only finitely many points $(\mathbf{k}; \boldsymbol{\ell})$ such that $\boldsymbol{\alpha}^{\mathbf{k}}\boldsymbol{\beta}^{\boldsymbol{\ell}}$ has degree $d$, we could (in principle) compute, using (7.1), all of the well-defined transvectants $[\mathbf{k}; \boldsymbol{\ell}; s]$ with these values of $\mathbf{k}$ and $\boldsymbol{\ell}$. By the injectivity of $\mathbb{T}$, at most one of these can equal $h$, and if one does, we have detected the values of $\mathbf{k}$, $\boldsymbol{\ell}$, and $s$. Another proof (more effective if the computation is actually needed) is to apply the splitting algorithm of Sanders ([28] and [22, §2.6]). This provides an algorithm to express any polynomial $h \in \mathbb{R}[\mathbf{x}, \mathbf{y}]$ in terms of the chain-weight basis determined by the pre-box product basis for $\mathcal{I}$. This will reveal whether $h$ is a top weight vector, and if so, which one. □

**11.2. Remark.** Lemma 11.1 is used in the proof of Theorem 11.9 below, and elsewhere to see that certain operations on transvectants are well-defined. For instance, in [26] we defined the *stripped form* and *total transvectant strength* of a transvectant monomial without noticing that (in the absence of Lemma 11.1) our definition was only valid for *transvectant symbols*, since it required knowledge of $\mathbf{k}$, $\boldsymbol{\ell}$, and $s$. Lemma 11.1 is not needed for computations, since all of our computations are done using transvectant symbols.

Next, a standard theorem in integer programming clarifies the notion of "prime $\boldsymbol{\alpha\beta}$-transvectant", introduced informally in Section 10. We begin by observing that *factorization* of a transvectant (up to equivalence), as defined in Section 9, corresponds to *decomposition* of its integer vector in $\mathcal{C}$ into a nonnegative linear combination of other integer vectors in $\mathcal{C}$.

**11.3. Theorem.** *With fixed $\boldsymbol{\alpha}$ and $\boldsymbol{\beta}$, there exists a unique finite set of* **prime** *$\boldsymbol{\alpha\beta}$-transvectants, with the following properties: every $\boldsymbol{\alpha\beta}$-transvectant that is not prime is equivalent ($\equiv$) to a product of two or more prime transvectants. No prime transvectant is equivalent to such a product.*

**Proof.** In integer programming, a **spanning set** for an integer cone $\mathcal{C}$ is a subset of $\mathcal{C}$ such that each element of the cone can be written (not necessarily uniquely) as a nonnegative integer linear combination of the vectors in the spanning set. It is known (see [33, §1.4], [31, Thm. 16.4]) that an integer cone has a *unique* finite minimal spanning set. No element of the minimal spanning set can be written as a nonnegative integer linear combination of other points in $\mathcal{C}$. Theorem 11.3 is a restatement of these facts in the appropriate multiplicative language for transvectants. □

**11.4. Remark.** The unique minimal spanning set for $\mathcal{C}$ is known in integer programming as the Hilbert basis for $\mathcal{C}$. We will not use this name because it would cause confusion. The Hilbert basis for $\mathcal{I}$ that we wish to find is generally a *subset* of the primes in $\mathcal{C}$.

Let $\Pi$ be the set of primes in $\mathcal{C}$, and let $\Pi^* = \Pi \cap \mathcal{C}^*$. This finite set of integer vectors will be written as $\Pi^* = \{\boldsymbol{\pi}_1, \dots, \boldsymbol{\pi}_t\}$, with $\boldsymbol{\pi}_i = (\boldsymbol{\kappa}_i; \boldsymbol{\lambda}_i; \sigma_i) \in \mathbb{N}^{p+q+1}$. Also let $\vartheta_i = \mathbb{T}(\boldsymbol{\pi}_i) = [\boldsymbol{\kappa}_i; \boldsymbol{\lambda}_i; \sigma_i] \in \mathbb{R}[\mathbf{x}, \mathbf{y}]$. The $\boldsymbol{\pi}_i$ will be called the **admissible primes** in $\mathcal{C}$ (determined by $\acute{A}$ and $\grave{A}$), while the $\vartheta_i$ will be called the **admissible prime transvectants**. These can be enumerated in any order, but once $\boldsymbol{\vartheta} = (\vartheta_1, \dots, \vartheta_t)$ is fixed, it will be called the **factoring order**, since the box product relative to $\boldsymbol{\vartheta}$ will be found by factoring the primes from all transvectants in this order. Usually we choose a factoring order having the form

$$\boldsymbol{\vartheta} = (\vartheta_1, \dots, \vartheta_t) = (\alpha_1, \dots, \alpha_p, \beta_1, \dots, \beta_q, \gamma_1, \dots, \gamma_r) = (\boldsymbol{\alpha}, \boldsymbol{\beta}, \boldsymbol{\gamma}), \tag{11.2}$$

where $(\gamma_1, \dots, \gamma_r)$ are the admissible prime transvectants of strength greater than zero, enumerated in some order, but factoring orders not of this type are also possible (as in Example 16.3). Note that $t = p + q + r$. Monomials in $\boldsymbol{\vartheta}$ will be written as $\boldsymbol{\vartheta}^{\mathbf{K}} = \vartheta_1^{K_1} \cdots \vartheta_t^{K_t}$, with $\mathbf{K} \in \mathbb{N}^t$, the Newton space for these monomials.

**11.5. Remark.** In writing (11.2) we have assumed that all of the zero-transvectants $\alpha_i = (\alpha_i, 1)^{(0)}$ and $\beta_j = (1, \beta_j)^{(0)}$ are admissible. (It is obvious that they are prime.) This is a natural assumption, which can always be arranged. Suppose, for instance, that $\boldsymbol{\alpha}$ contains an invariant $\alpha_i$ that is not an element of $\acute{A}$. By Lemma 3.6, since $\acute{A}$ is assumed to be standard, no monomial in $\acute{A}$ will contain $\alpha_i$. Therefore $\alpha_i$ is a **redundant element** of the Hilbert basis $\acute{A}$ for $\acute{\mathcal{I}}$, and can be omitted from that Hilbert basis, since it is never used.

**11.6. Remark.** The simplest example having a well-defined transvectant that is not admissible (that is, $\mathcal{C}^* \neq \mathcal{C}$) is $\mathcal{I}_{24} = \mathcal{I}_2 \boxtimes \mathcal{I}_4$. The input data for this problem comes from (4.4) and (4.5) for $\boldsymbol{\alpha}$ and $\acute{A}$, (4.10) and (4.11) for $\boldsymbol{\beta}$ and $\grave{A}$. The general $\boldsymbol{\alpha\beta}$-transvectant is $(\alpha_1^{k_1}, \beta_1^{\ell_1}\beta_2^{\ell_2}\beta_3^{\ell_3}\beta_4^{\ell_4})^{(s)}$, with no restrictions except well-definedness. But the general *admissible* transvectant has the additional restriction $\ell_3 = 0$ or 1. The simplest example having a *prime* transvectant that is not admissible (so that $\mathcal{C}$ contains a prime that is not in $\mathcal{C}^*$) is $\mathcal{I}_{34}$. The input data is (4.7), (4.8), (4.10), and (4.11). The transvectant $(\alpha_1^3, \beta_3^2)^{(6)}$ is well-defined and prime, but is not admissible, since (4.11) does not allow the second power of $\beta_3$. To see that this is prime, note than any viable factorization must have the form $(\alpha_1^2, \beta_3)^{(s')}(\alpha_1, \beta_3)^{(s'')}$ with $s' + s'' = 6$, but there are no such choices of $s'$ and $s''$ that are supported by both of the weights in these transvectants.

Since $\mathcal{C}^*$ may not be a cone, it will not always have a spanning set. But $\Pi^*$ is a **weak spanning set**, in the sense that every element of $\mathcal{C}^*$ is a sum of elements of $\Pi^*$, although not all such sums must belong to $\mathcal{C}^*$.

**11.7. Theorem.** *The set $\Pi^* = \{\boldsymbol{\pi}_1, \dots, \boldsymbol{\pi}_t\}$ of admissible primes is a weak spanning set for $\mathcal{C}^*$, and the ordered set $\boldsymbol{\vartheta} = (\vartheta_1, \dots, \vartheta_t)$ of admissible prime transvectants forms a Hilbert basis for $\mathcal{I}$.*

**Proof.** Consider $(\mathbf{k}; \boldsymbol{\ell}; s) \in \mathcal{C}^*$. Since this also belongs to $\mathcal{C}$, it can be expressed as a nonnegative integer linear combination of elements of $\Pi$, or equivalently, as a sum of elements of $\Pi$ (with repetitions allowed). Suppose that $(\boldsymbol{\kappa}; \boldsymbol{\lambda}; \sigma)$ is a prime that appears in this expression, so that

$$(\mathbf{k}; \boldsymbol{\ell}; s) = (\boldsymbol{\kappa}; \boldsymbol{\lambda}; \sigma) + \cdots \qquad \text{(finite sum)}.$$

Let $\mathbf{t} = (t_1, \dots, t_p)$ and $\mathbf{u} = (u_1, \dots, u_q)$ be indeterminates substituting for $\boldsymbol{\alpha}$ and $\boldsymbol{\beta}$, as in Section 3, and regard $\acute{A}$ and $\grave{A}$ as sets of monomials in $\mathbf{t}$ and $\mathbf{u}$ respectively. Since $(\mathbf{k}; \boldsymbol{\ell}; s) \in \mathcal{C}^*$, we have $\mathbf{t}^{\mathbf{k}} \in \acute{A}$ and $\mathbf{u}^{\boldsymbol{\ell}} \in \grave{A}$. But $\mathbf{t}^{\boldsymbol{\kappa}}$ is a factor of $\mathbf{t}^{\mathbf{k}}$, and $\mathbf{u}^{\boldsymbol{\lambda}}$ is a factor of $\mathbf{u}^{\boldsymbol{\ell}}$. Since $\acute{A}$ and $\grave{A}$ are standard, Lemma 3.6 implies that

$$\mathbf{t}^{\boldsymbol{\kappa}} \in \acute{A} \quad \text{and} \quad \mathbf{u}^{\boldsymbol{\lambda}} \in \grave{A}. \tag{11.3}$$

Therefore $(\boldsymbol{\kappa}; \boldsymbol{\lambda}; \sigma) \in \mathcal{C}^*$. Thus the expression for $(\mathbf{k}; \boldsymbol{\ell}; s)$ as a sum of primes uses only primes in $\mathcal{C}^*$, and $\Pi^*$ is a weak spanning set for $\mathcal{C}^*$. Applying $\mathbb{T}$, it follows that every admissible transvectant can be replaced by a product of the admissible prime transvectants $\vartheta_i$. In other words every admissible transvectant $[\mathbf{k}; \boldsymbol{\ell}; s]$ is equivalent ($\equiv$) to a monomial $\boldsymbol{\vartheta}^{\mathbf{K}}$. By Theorem 9.2, every invariant $f \in \mathcal{I}$ can be written as a (possibly infinite) linear combination of such monomials, so every invariant is a formal power series in $\boldsymbol{\vartheta}$, and $\boldsymbol{\vartheta}$ is a Hilbert basis for $\mathcal{I}$. □

We can now give the precise definition of the box product of preferred sets $\acute{A}$ and $\grave{A}$. Let a factoring order $\boldsymbol{\vartheta}$ be fixed. For each transvectant $\delta = [\mathbf{k}; \boldsymbol{\ell}; s] \in \acute{A} \mathbin{\widehat{\boxtimes}} \grave{A}$, we construct a factorization $\boldsymbol{\vartheta}^{\mathbf{K}}$ of $\delta$ as follows: factor $\vartheta_1$ from $\delta$ as many times as possible, then factor $\vartheta_2$ as many times as possible from what remains, and so forth, until $\vartheta_K$ has been factored. This must result in a complete factorization of $\delta$, since $\boldsymbol{\vartheta}$ is the complete list of admissible prime transvectants; this factorization is called the **preferred factorization** of $\delta$ under the factoring order $\boldsymbol{\vartheta}$. Note that $\delta \equiv \boldsymbol{\vartheta}^{\mathbf{K}}$. The set of all such preferred factorizations of the transvectants in $\acute{A} \mathbin{\widehat{\boxtimes}} \grave{A}$ is denoted by

$$A = \acute{A} \boxtimes_{\boldsymbol{\vartheta}} \grave{A}, \tag{11.4}$$

and is called the **box product of $\acute{A}$ and $\grave{A}$ with factoring order $\boldsymbol{\vartheta}$**. Since $A$ is a replacement for the pre-box product $\acute{A} \mathbin{\widehat{\boxtimes}} \grave{A}$, it follows from Theorem 9.2 that $A$ is a formal vector space basis for $\mathcal{I}$. Clearly $A$ is also a preferred set of monomials in $\boldsymbol{\vartheta}$ (in the sense of Section 3), and we will show in Theorem 11.9 that it is a standard preferred set.

The next theorem characterizes the elements of $A$. Let $\Omega : \mathcal{C}^* \to \mathbb{R}[\mathbf{x}, \mathbf{y}]$ be defined by

$$\Omega(\mathbf{K}) = \mathbb{T}(K_1\boldsymbol{\pi}_1 + \cdots + K_t\boldsymbol{\pi}_t). \tag{11.5}$$

Then $\Omega(\mathbf{K})$ is the unique transvectant that can be factored as $\boldsymbol{\vartheta}^{\mathbf{K}}$, and, given a transvectant $\delta$, $\Omega^{-1}(\delta)$ is the set of all $\mathbf{K}$ such that $\boldsymbol{\vartheta}^{\mathbf{K}}$ is an admissible factorization of $\delta$.

**11.8. Theorem.** *Let $\delta \equiv \boldsymbol{\vartheta}^{\mathbf{K}}$ be the preferred factorization of $\delta$ with order $\boldsymbol{\vartheta}$. Then $\mathbf{K}$ is the largest exponent of $\boldsymbol{\vartheta}$ in lexicographic order among all possible factorizations of $\delta$ in terms of admissible primes. That is,* $\mathbf{K} = \operatorname{maxlex} \Omega^{-1}(\delta)$.

**Proof.** Let $\boldsymbol{\vartheta}^{\mathbf{K}}$ be the preferred factorization of $\delta$, and let $\mathbf{K}' = \operatorname{maxlex} \Omega^{-1}$. Suppose that these are different, and let $j$ be the first index for which $K_j \neq K'_j$; then $K_j < K'_j$. Then consider the partial factorization

$$\delta \equiv (\vartheta_1^{K_1} \cdots \vartheta_{j-1}^{K_{j-1}})\vartheta_j^{K_j} \varepsilon_j$$

obtained at the $j$-th stage of computing the preferred factorization of $\delta$. Here we are supposed to have factored as many copies of $\vartheta_j$ as possible from what remains, so that $\vartheta_j$ cannot be factored again from $\varepsilon_j$. But, by hypothesis, there exists a factorization

$$\delta \equiv (\vartheta_1^{K_1} \cdots \vartheta_{j-1}^{K_{j-1}})\vartheta_j^{K'_j} \varepsilon'_j$$

with $K'_j > K_j$. Therefore $\boldsymbol{\vartheta}^{\mathbf{K}}$ is not the preferred factorization of $\delta$, which is a contradiction. □

The next theorem shows that the box product can be iterated: the output of a box product is suitable as an input to another box product.

**11.9. Theorem.** *The box product $A = \acute{A} \boxtimes_{\boldsymbol{\vartheta}} \grave{A}$ is a standard preferred set with respect to the Hilbert basis $\boldsymbol{\vartheta}$ for $\mathfrak{I}$.*

**Proof.** According to Lemma 3.6, we must prove that any divisor of an element of $A$ also belongs to $A$. The idea is simple: if $\boldsymbol{\vartheta}^K$ is a preferred factorization, then every factor $\boldsymbol{\vartheta}^{K'}$ of $\boldsymbol{\vartheta}^K$ must be preferred, since if not, there would exist a $\boldsymbol{\vartheta}^{K''}$ that is $\equiv \boldsymbol{\vartheta}^{K'}$ but greater in lex order. Then, replacing $\boldsymbol{\vartheta}^{K'}$ with $\boldsymbol{\vartheta}^{K''}$ in $\boldsymbol{\vartheta}^K$ would create a $\boldsymbol{\vartheta}^{K'''}$ that would be $\equiv \boldsymbol{\vartheta}^K$ but greater in lex order, which is impossible. It is worth spelling out the details, since this is a point at which Lemma 11.1 plays a crucial role.

For clarity, let $\mathbf{v} = (v_1, \dots, v_t)$ be a set of indeterminates equal in number to the elements of the Hilbert basis $\boldsymbol{\vartheta}$, and regard $A$ as a set of monomials in $\mathbf{v}$. Let $\mathbf{v}^{\mathbf{K}} \in A$; this means that $\boldsymbol{\vartheta}^{\mathbf{K}}$ is the preferred factorization of the transvectant $\Omega(\mathbf{K})$, which in turn means that $\mathbf{K} = \text{maxlex } \Omega^{-1}\Omega(\mathbf{K})$. Let $\mathbf{v}^{\mathbf{K}'}$ be any divisor of $\mathbf{v}^{\mathbf{K}}$ (under polynomial division in $\mathbb{R}[\mathbf{v}]$); that is, $K_i' \le K_i$ for $i = 1, \dots, t$, or $\mathbf{K} - \mathbf{K}' \in \mathbb{N}^t$. We must show that $\mathbf{v}^{\mathbf{K}'} \in A$; that is, we must show that $\mathbf{K}' = \text{maxlex } \Omega^{-1}\Omega(\mathbf{K}')$.

Suppose (for a proof by contradiction) that this is not true, and that $\mathbf{K}'' \ne \mathbf{K}'$ is instead the largest element of $\Omega^{-1}\Omega(\mathbf{K}')$. Then $\Omega(\mathbf{K}'') = \Omega(\mathbf{K}')$, which (since $\mathbb{T}$ is injective) implies

$$K_1''\boldsymbol{\pi}_1 + \cdots + K_t''\boldsymbol{\pi}_t = K_1'\boldsymbol{\pi}_1 + \cdots + K_t'\boldsymbol{\pi}_t.$$

Also, if $j$ is the first index for which $K_j'' \ne K_j'$, we have $K_j'' > K_j'$. Consider now

$$\mathbf{K}''' = (\mathbf{K} - \mathbf{K}') + \mathbf{K}'' = \mathbf{K} + (\mathbf{K}'' - \mathbf{K}').$$

The first expression for $\mathbf{K}'''$ shows that $\mathbf{K}'''$ is the sum of two elements of $\mathbb{N}^t$, so it belongs to $\mathbb{N}^t$. The second shows that $\mathbf{K}'''$ is greater in lex order than $\mathbf{K}$, since the first nonzero entry in $\mathbf{K}'' - \mathbf{K}'$ is positive. Finally, it follows that

$$K_1'''\boldsymbol{\pi}_1 + \cdots + K_t'''\boldsymbol{\pi}_t = K_1\boldsymbol{\pi}_1 + \cdots + K_t\boldsymbol{\pi}_t,$$

so that $\boldsymbol{\vartheta}^{\mathbf{K}'''} \equiv \boldsymbol{\vartheta}^{\mathbf{K}}$. This contradicts the hypothesis that $\vartheta^{\mathbf{K}}$ is the preferred factorization of $\Omega(\mathbf{K})$. $\square$

## 12. The factoring method: a preliminary example

In subsequent sections we will develop some rather complicated algorithms to calculate box products by the factoring method. But the method itself is, in principle, simple. The example in this section (already done once by an *ad hoc* method in Example 10.1) will be solved here by the most elementary version of the factoring method, to demonstrate this fundamental simplicity. But this "simple" method requires *thought*, and becomes difficult to carry out in hard problems. The algorithmic method, on the other hand, requires only *rule-following*. The "simple" method illustrated here was the starting point for the development of the algorithmic method.

**12.1. Example.** $N_{22}$, second time.

As in Example 10.1, we begin with the general element of $\acute{A}_2 \,\widehat{\boxtimes}\, \grave{A}_2$, written (omitting the subscripts 1) as $(\alpha^k, \beta^\ell)^{(s)}$, with $0 \le s \le k$ and $0 \le s \le \ell$. From this we factor as many copies of $\alpha$ as possible, and write the result as $(\alpha^k, \beta^\ell)^{(s)} \equiv \alpha^i(\alpha^{k-i}, \beta^\ell)^{(s)}$. Since the transvectant in this factorization must be well-defined, we must have $0 \le s \le k - i$ in addition to the previous inequalities. Since $i$ is as large as possible, the transvectant $(\alpha^{k-i-1}, \beta^\ell)^{(s)}$ must be undefined; this implies $k - i - 1 < s$. From these last two inequalities we deduce that $k - i = s$, so the factorization so far achieved is actually $\alpha^i(\alpha^s, \beta^\ell)^{(s)}$. Similar reasoning shows that if we factor out as many $\beta$ as possible, we obtain $\alpha^i\beta^j(\alpha^s, \beta^s)^{(s)}$. If we do not know the next prime $\gamma$, we take the smallest possible nontrivial value of $s$, namely $s = 1$, and observe that $(\alpha, \beta)^{(1)}$ cannot be factored further using the already known primes $\alpha$ and $\beta$; therefore $\gamma = (\alpha, \beta)^{(1)}$ is prime. So we factor as many copies of $\gamma$ as possible from $(\alpha^s, \beta^s)^{(s)}$, obtaining $\gamma^s$. Therefore a complete factorization of the original general transvectant is $\alpha^i\beta^j\gamma^s$; this shows at the same time that $\alpha$, $\beta$, and $\gamma$ are all the primes that exist. But we have not yet established the ranges of the integer variables $i$, $j$, and $s$. Since $k$ can be arbitrarily larger than $s$ initially, $i \in [0, \infty]$; similarly for $j$; and $s$ has no absolute bound since $k$ and $\ell$ have none. Therefore the block decomposition for the factored expressions $\alpha^i\beta^j\gamma^s$ is (10.2). □

We can draw three morals from this example that will motivate the next three sections. **1.** It is necessary to keep track of several integer inequalities at once, and the number of inequalities increases in larger problems. The p-matrix notation of Section 13 handles these inequalities automatically. **2.** It is necessary to carry out operations of propositional logic involving these inequalities. For instance, it is crucial to be able to find the set of transvectants that are well-defined, but are *not* divisible by a particular prime. (In this example, the transvectants $(\alpha^s, \beta^\ell)^{(s)}$ are not divisible by $\alpha$.) In general, such a set will be a disjoint union of several subsets described by p-matrices. The need for a *disjoint* union means that the logic must use the *exclusive* "or" operation $\veebar$, rather than the inclusive operation $\vee$. These logical issues will be handled by Algorithms A and B in Section 14, and the *filtration* process in Section 15. **3.** Finding the ranges of various integers (here $i$, $j$, $s$) is best viewed as a separate problem from finding the factorizations themselves. This is done by Algorithms C and D in Section 14, and the *reconstruction* process in Section 15.

**12.2. Example.** $N_{23}$, second time.

The interested reader is invited to attempt Example 10.2 by the method of this section, using the factoring order

$$\vartheta_1 = \alpha_1,\ \vartheta_2 = \beta_1,\ \vartheta_3 = \beta_2,\ \vartheta_4 = \gamma_1,\ \vartheta_5 = \gamma_2$$

and the following outline. The general transvectant has the form $(\alpha_1^{k_1}, \beta_1^{\ell_1}\beta_2^{\ell_2})^{(s)}$. In the first step one factors as many $\alpha_1$ as possible from this transvectant, to obtain $\alpha_1^i(\alpha_1^{k_1}, \beta_1^{\ell_1}\beta_2^{\ell_2})^{(s)}$, in which $k_1$ has changed. (The new $k_1$ is the old $k_1 - i$.) The fact that no more $\alpha_1$ can be removed implies $s = k_1$, but it is best not to replace $k_1$ by $s$ (as we did in the previous example), but instead to add $s = k_1$ (or $k_1 \le s \le k_1$) to a growing list of inequalities for well-definedness. At the second step (factoring $\beta_1$), we obtain $\alpha_1^i\beta_1^j(\alpha_1^{k_1}, \beta_1^{\ell_1}\beta_2^{\ell_2})^{(s)}$ (where $\ell_1$ has changed), and the condition $(s = 2\ell_1 - 1) \veebar (s = 2\ell_1)$ joins the list. In step 3, it is possible to remove all of the $\beta_2$ from the transvectant, because $\beta_2$ has weight zero. The most interesting step is the fourth, in which $\gamma_1$ is factored. At this point we get $\alpha_1^i \beta_1^j \beta_2^p \gamma_1^q (\alpha_1^{k_1}, \beta_1^{\ell_1})^{(s)}$. The important thing at this stage is to show that $(q = 0) \veebar (q = 1)$; that is, no more than one copy of $\gamma_1$ can be removed, as we already saw in Example 10.2 by an *ad hoc* trick. Readers who do this exercise will be able to compare their solution with the algorithmic one given later in Example 16.2. □

**12.3. Remark.** Hilbert basis elements of weight zero belonging to $\boldsymbol{\alpha}$ or $\boldsymbol{\beta}$ will be called **null input primes**. Null input primes can be factored out of any transvectant at any time, since they do not provide support for the transvectant strength. Therefore their position in the factoring order does not matter. It is simplest to factor them all out at once at the beginning. This can be done by *omitting* the null input primes entirely, and creating a ring $\mathcal{R} = \mathbb{R}[[\cdots]]$, where $\cdots$ is the list of null inputs, to appear in place of $\mathbb{R}$ when the final Stanley decomposition is written down. For instance, in $\mathcal{I}_{23}$ we have $\mathcal{R} = \mathbb{R}[[\beta_2]]$ and the Stanley decomposition (10.7) becomes $\mathcal{R}[[\alpha_1, \beta_1, \gamma_2]] \oplus \mathcal{R}[[\alpha_1, \beta_1, \gamma_1]]\gamma_1$. This is called **suppression of null inputs**.

## 13. The p-matrix notation for sets of transvectants

The algorithms for box products in the next section will call for many sets of transvectants $[\mathbf{k}; \boldsymbol{\ell}; s]$ in which $\mathbf{k}$ and $\boldsymbol{\ell}$ are confined to blocks, but $s$ is confined to an interval that is not fixed but depends on the values of $\mathbf{k}$ and $\boldsymbol{\ell}$. These sets will be represented by p-matrices (parenthesis matrices). The simplest example is the pre-box product of two blocks, which will be written as

$$\begin{bmatrix}\mathbf{b}\\ \mathbf{a}\end{bmatrix} \widehat{\boxtimes} \begin{bmatrix}\mathbf{d}\\ \mathbf{c}\end{bmatrix} = \left(\begin{array}{c|c|cc}\mathbf{b} & \mathbf{d} & \widehat{\boldsymbol{\alpha}}\mathbf{k} & \widehat{\boldsymbol{\beta}}\boldsymbol{\ell}\\ \mathbf{a} & \mathbf{c} & 0 & 0\end{array}\right). \tag{13.1}$$

The three partitions of this matrix give the block for $\mathbf{k}$, the block for $\boldsymbol{\ell}$, and two scalar intervals for $s$, one depending on $\mathbf{k}$ and the other on $\boldsymbol{\ell}$; $s$ must belong to the intersection of these. The general case is defined as follows.

Let $\mathbf{a} \in \mathbb{Z}^p$, $\mathbf{b} \in (\mathbb{Z} \cup \{\infty\})^p$, $\mathbf{c} \in \mathbb{Z}^q$, $\mathbf{d} \in (\mathbb{Z} \cup \{\infty\})^q$. Let $E, F, G, H$ be integer-valued functions

$$E, F : \begin{bmatrix}\mathbf{b}\\ \mathbf{a}\end{bmatrix} \to \mathbb{N}, \qquad G, H : \begin{bmatrix}\mathbf{d}\\ \mathbf{c}\end{bmatrix} \to \mathbb{N}.$$

Then the p-matrix

$$P = \left(\begin{array}{c|c|cc}\mathbf{b} & \mathbf{d} & F & H\\ \mathbf{a} & \mathbf{c} & E & G\end{array}\right) \tag{13.2}$$

stands (initially) for the set of integer vectors

$$P = \left\{(\mathbf{k}; \boldsymbol{\ell}; s) \in \mathbb{N}^{p+q+1} : \quad \mathbf{k} \in \begin{bmatrix}\mathbf{b}\\ \mathbf{a}\end{bmatrix}, \quad \boldsymbol{\ell} \in \begin{bmatrix}\mathbf{d}\\ \mathbf{c}\end{bmatrix}, \right.$$
$$\left. s \in [E(\mathbf{k}), F(\mathbf{k})] \cap [G(\boldsymbol{\ell}), H(\boldsymbol{\ell})] \right\}. \tag{13.3}$$

As with the notation $B$ for blocks, $P$ is used ambiguously for the matrix in (13.2) and the set of integer vectors in (13.3). If $F(\mathbf{k}) \le \widehat{\boldsymbol{\alpha}}\mathbf{k}$ and $G(\boldsymbol{\ell}) \le \widehat{\boldsymbol{\beta}}\boldsymbol{\ell}$, then each $(\mathbf{k}; \boldsymbol{\ell}; s)$ in $P$ belongs to $\mathcal{C}$, and the transvectant $[\mathbf{k}; \boldsymbol{\ell}; s]$ is defined. In this case $P$ will also be used to denote the set of these transvectants.

Note that $\mathbf{a}$, $\mathbf{b}$, $\mathbf{c}$, and $\mathbf{d}$ are allowed to have negative entries. However, the vectors in $P$ contain only nonnegative integers, since according to (13.3) they must belong to $\mathbb{N}^{p+q+1}$. See Remark 3.5.

The functions $E, F, G, H$ in (13.2) appear without variables (not as $E(\mathbf{k}), F(\mathbf{k}), G(\boldsymbol{\ell}), H(\boldsymbol{\ell})$). This is logically correct, because $\mathbf{k}$ and $\boldsymbol{\ell}$ appear in (13.4) as bound or dummy variables, so they are "not really present". However, it is awkward to use "functions without variables" in practice. One could use the "maps-to" symbol $\mapsto$, so that (13.1) would appear as:

$$\left(\begin{array}{c|c|cc} \mathbf{b} & \mathbf{d} & \mathbf{k} \mapsto \widehat{\boldsymbol{\alpha}}\mathbf{k} & \boldsymbol{\ell} \mapsto \widehat{\boldsymbol{\beta}}\boldsymbol{\ell} \\ \mathbf{a} & \mathbf{c} & \mathbf{k} \mapsto 0 & \boldsymbol{\ell} \mapsto 0 \end{array}\right),$$

with $\mathbf{k}$ and $\boldsymbol{\ell}$ again as dummy variables. Another solution is to "decorate" the P-matrix as

$$P = \left(\begin{array}{c|c|cc} \mathbf{k} & \boldsymbol{\ell} & \multicolumn{2}{c}{s} \\ \hline \mathbf{b} & \mathbf{d} & F(\mathbf{k}) & H(\boldsymbol{\ell}) \\ \mathbf{a} & \mathbf{c} & E(\mathbf{k}) & G(\boldsymbol{\ell}) \end{array}\right). \tag{13.4}$$

In this notation the column labels $\mathbf{k}$ and $\boldsymbol{\ell}$ serve as the binding operators; this allows for changes of dummy variables, as in the proof of Algorithm 14.4 below. To save space, our general practice is to omit the column labels. A different type of column labels, as in

$$P = \left(\begin{array}{c|c|cc} \boldsymbol{\alpha} & \boldsymbol{\beta} & \multicolumn{2}{c}{s} \\ \hline \mathbf{b} & \mathbf{d} & F(\mathbf{k}) & H(\boldsymbol{\ell}) \\ \mathbf{a} & \mathbf{c} & E(\mathbf{k}) & G(\boldsymbol{\ell}) \end{array}\right), \tag{13.5}$$

are often more helpful. Here $\boldsymbol{\alpha}$ and $\boldsymbol{\beta}$ indicate the Hilbert bases that receive the exponents $\mathbf{k}$ and $\boldsymbol{\ell}$. (In applications, these may have labels other than $\boldsymbol{\alpha}$ and $\boldsymbol{\beta}$.)

The p-matrix

$$\left(\begin{array}{c|cc|cc} \infty & \infty & 0 & k_1 & 2\ell_1 \\ 0 & 0 & 0 & k_1 & 2\ell_1 - 1 \end{array}\right), \tag{13.6}$$

arises naturally in Example 16.2 below. Notice that $E = F$ in this p-matrix; both are the function $\mathbf{k} \mapsto k_1$. So the fourth column of (13.6) just states that $s = k_1$. To make this easier to see, **whenever the two entries in a column of a p-matrix or a b-matrix are equal, we will replace the two entries by a single entry on an intermediate line.** Thus (13.6) will be written as

$$\left(\begin{array}{c|cc|cc} \infty & \infty & 0 & \multirow{2}{*}{$k_1$} & 2\ell_1 \\ 0 & 0 & 0 & & 2\ell_1 - 1 \end{array}\right). \tag{13.7}$$

This example also illustrates the reason for the convention about negative entries in Remark 3.5. The expression $2\ell_1 - 1$ takes a negative value if $\ell_1 = 0$, so the notation must be prepared to handle this.

The behavior of a p-matrix, compared to a b-matrix, can be somewhat subtle. For instance, in (13.7), it appears that $k_1$ takes every value from zero to infinity, as indicated by the first column. But in fact the fourth column declares that $k_1 = s$ and the fifth column says $s \in [2\ell_1 - 1, 2\ell_1]$, so it follows that $k_1 \in [2\ell_1 - 1, 2\ell_1]$. In fact, $k_1$ does take on every value from zero to infinity, but only when $\ell_1$ is allowed to vary. For any particular value of $\ell_1$, only two values of $k_1$ are allowed. Thus, to understand the true range of $k_1$, it is necessary to take into account *all* of the columns of the p-matrix, and how these separate restrictions interact with each other.

## 14. Algorithms for p-matrices

The factoring method is based on four algorithms, A, B, C, and D, stated and proved in this section. The next section will show how to deploy these algorithms to compute box products. The algorithms are briefly summarized as follows. Recall the definition of *divisible* in Remark 9.3.

**A:** Input: A p-matrix $P$ of well-defined transvectants, and a transvectant $(\vartheta)$. Output: A p-matrix denoted $P(\vartheta)$ giving the elements of $P$ that are divisible by $\vartheta$.
**B:** Input: A p-matrix $P$ of well-defined transvectants, and a p-matrix $P' \subset P$. Output: A disjoint union of p-matrices representing $P \setminus P'$.
**C:** Input: A p-matrix $P$ of points in $\mathbb{N}^{(p+q+1)}$, a point $(\mathbf{k}; \boldsymbol{\lambda}; \sigma) \in \mathbb{N}^{(p+q+1)}$, and a nonnegative integer $j$. Output: A p-matrix denoted $P[j]$ representing the translation of the point set $P$ through the vector $-j(\boldsymbol{\kappa}; \boldsymbol{\lambda}; \sigma)$, intersected with $\mathbb{N}^{(p+q+1)}$.
**D:** Input: Two p-matrices $P$ and $Q$. Output: a p-matrix representing $P \cap Q$.

**14.1. Theorem** *(Algorithm A). Let*

$$P = \left( \begin{array}{c|c|cc} \mathbf{b} & \mathbf{d} & F(\mathbf{k}) & H(\boldsymbol{\ell}) \\ \mathbf{a} & \mathbf{c} & E(\mathbf{k}) & G(\boldsymbol{\ell}) \end{array} \right) \tag{14.1}$$

*be a p-matrix of well-defined transvectants (that is, $F(\mathbf{k}) \le \widehat{\boldsymbol{\alpha}}\mathbf{k}$ and $H(\boldsymbol{\ell}) \le \widehat{\boldsymbol{\beta}}\boldsymbol{\ell}$), and let*

$$\vartheta = [\boldsymbol{\kappa}; \boldsymbol{\lambda}; \sigma] = (\boldsymbol{\alpha}^{\boldsymbol{\kappa}}, \boldsymbol{\beta}^{\boldsymbol{\lambda}})^{(\sigma)}$$

*be any well-defined transvectant. Then the subset $P(\vartheta) \subset P$ consisting of transvectants divisible by $\vartheta$ is given by*

$$P(\vartheta) = \left( \begin{array}{c|c|cc} \mathbf{b} & \mathbf{d} & F'(\mathbf{k}) & H'(\boldsymbol{\ell}) \\ \mathbf{a}' & \mathbf{c}' & E'(\mathbf{k}) & G'(\boldsymbol{\ell}) \end{array} \right), \tag{14.2}$$

*where*

$$\begin{aligned} \mathbf{a}' &= \max\{\mathbf{a}, \boldsymbol{\kappa}\}, \\ \mathbf{c}' &= \max\{\mathbf{c}, \boldsymbol{\lambda}\}, \\ E'(\mathbf{k}) &= \max\{E(\mathbf{k}), \sigma\}, \\ F'(\mathbf{k}) &= \min\{F(\mathbf{k}), \widehat{\boldsymbol{\alpha}}(\mathbf{k} - \boldsymbol{\kappa}) + \sigma\}, \\ G'(\boldsymbol{\ell}) &= \max\{G(\boldsymbol{\ell}), \sigma\}, \\ H'(\boldsymbol{\ell}) &= \min\{H(\boldsymbol{\ell}), \widehat{\boldsymbol{\beta}}(\boldsymbol{\ell} - \boldsymbol{\lambda}) + \sigma\}, \end{aligned} \tag{14.3}$$

*the maxima being computed componentwise in the case of vectors* $\mathbf{a}'$ *and* $\mathbf{c}'$ *(so that* $a'_i = \max\{a_i, \kappa_i\}$ *for* $i = 1, \dots, p$ *and* $c'_j = \max\{c_j, \lambda_j)$ *for* $j = 1, \dots, q$*).*

Note that there are no symbols $\mathbf{b}'$ and $\mathbf{d}'$ (or we could say $\mathbf{b}' = \mathbf{b}$ and $\mathbf{d}' = \mathbf{d}$).

**Proof.** If $\delta = (\boldsymbol{\alpha}^{\mathbf{k}}, \boldsymbol{\beta}^{\boldsymbol{\ell}})^{(s)} \in P$ then

$$\begin{aligned} &a_i \le k_i \le b_i \quad \text{for} \quad i = 1, \dots, p; \\ &c_j \le \ell_j \le d_j \quad \text{for} \quad j = 1, \dots, q; \\ &E(\mathbf{k}) \le s \le F(\mathbf{k}); \\ &G(\boldsymbol{\ell}) \le s \le H(\boldsymbol{\ell}). \end{aligned} \tag{14.4}$$

For $\delta$ to belong to $P(\vartheta)$, it must be possible to write

$$(\boldsymbol{\alpha}^{\mathbf{k}}, \boldsymbol{\beta}^{\boldsymbol{\ell}})^{(s)} \equiv (\boldsymbol{\alpha}^{\boldsymbol{\kappa}}, \boldsymbol{\beta}^{\boldsymbol{\lambda}})^{(\sigma)} (\boldsymbol{\alpha}^{\mathbf{k}-\boldsymbol{\kappa}}, \boldsymbol{\beta}^{\boldsymbol{\ell}-\boldsymbol{\lambda}})^{(s-\sigma)}, \tag{14.5}$$

with the last transvectant being well-defined. That is, in addition to the inequalities already listed, $k_i$, $\ell_j$, and $s$ must satisfy

$$\begin{aligned} &\kappa_i \le k_i \quad \text{for} \quad i = 1, \dots, p; \\ &\lambda_j \le \ell_j \quad \text{for} \quad j = 1, \dots, q; \\ &s \le \widehat{\boldsymbol{\alpha}}(\mathbf{k} - \boldsymbol{\kappa}) + \sigma; \\ &s \le \widehat{\boldsymbol{\beta}}(\boldsymbol{\ell} - \boldsymbol{\lambda}) + \sigma. \end{aligned} \tag{14.6}$$

Putting these conditions together proves (14.3). The transvectants described by (14.3) are well-defined, since $F'(\mathbf{k}) \le F(\mathbf{k}) \le \widehat{\boldsymbol{\alpha}}\mathbf{k}$ and similarly for $H'(\ell)$. □

Algorithm B computes $P \setminus P'$ given p-matrices $P$ and $P' \subset P$. For the statement of this algorithm we introduce some new notations. The columns of $P$ (viewed as nonnegative integer intervals) are denoted by $I_1, \dots, I_{p+q+2}$, and the columns of $P'$ by $J_1, \dots, J_{p+q+2}$. The assumption that $P' \subset P$ implies that $J_i \subset I_i$ for $i = 1, \dots, p+q+2$. In general, $I_i \setminus J_i$ will be a disjoint union of two intervals, one at the "bottom" of $I_i$ and one at the "top:"

$$I_i \setminus J_i = K_i \sqcup L_i.$$

Of course, either or both of $K_i$ and $L_i$ can be empty. In fact, for our usual application $P'$ will be a p-matrix $P(\vartheta)$ produced by Algorithm A; in this case $L_i$ will *always* be empty except in the last two columns ($i = p + q + 1$ or $p + q + 2$), since the top entries in the other columns are equal in $P$ and $P'$.

**14.2. Theorem** *(Algorithm B). Let P and P′ be p-matrices as described above. Create a list of p-matrices as follows. The list is initially empty.*

*For each value of $j = 1, \dots, (p+q+2)$, add*

$$\begin{pmatrix} J_1 & \cdots & J_{j-1} & K_j & I_{j+1} \cdots & I_{p+q+2} \end{pmatrix}$$

*to the list if $K_j \neq \emptyset$, and add*

$$\begin{pmatrix} J_1 & \cdots & J_{j-1} & L_j & I_{j+1} \cdots & I_{p+q+2} \end{pmatrix}$$

*to the list if $L_j \neq \emptyset$.*

*The p-matrices in the final list are disjoint, and their union is $P \setminus P'$.*

**Proof.** To illustrate the notation without adding to its complexity, we prove the case $p = 1$, $q = 2$, which is typical of the general case. We have $(k_1; \ell_1, \ell_2; s) \in P$ if and only if

$$(k_1 \in I_1) \wedge (\ell_1 \in I_2) \wedge (\ell_2 \in I_3) \wedge (s \in I_4) \wedge (s \in I_5), \tag{14.7}$$

where $\wedge$ denotes "and". Similarly, $(k_1; \ell_1, \ell_2; s) \in P'$ if and only if

$$(k_1 \in J_1) \wedge (\ell_1 \in J_2) \wedge (\ell_2 \in J_3) \wedge (s \in J_4) \wedge (s \in J_5). \tag{14.8}$$

We want to find the set of $(k_1; \ell_1, \ell_2; s)$ that satisfy (14.7) but not (14.8). The usual rule for negation of a conjunction gives the negation of (14.8) as

$$(k_1 \notin J_1) \vee (\ell_1 \notin J_2) \vee (\ell_2 \notin J_3) \vee (s \notin J_4) \vee (s \notin J_5), \tag{14.9}$$

but this uses inclusive or and leads to p-matrices that are not disjoint. To obtain a disjoint union, (14.9) must be replaced by a statement using exclusive or (written as $\veebar$). The model for this in propositional logic is (with three propositions)

$$\neg(p \wedge q \wedge r) \equiv (\neg p) \veebar (p \wedge \neg q) \veebar (p \wedge q \wedge \neg r); \tag{14.10}$$

anything that is negated in one disjunct is affirmed in all later disjuncts, to prevent overlap. (The result depends on the order of $p, q, r$, but any order is correct; see Remark 14.3.) Thus the negation of (14.8) is

$$\begin{aligned} &[k_1 \notin J_1] \veebar [(k_1 \in J_1) \wedge (\ell_1 \notin J_2)] \\ &\quad \veebar [(k_1 \in J_1) \wedge (\ell_1 \in J_2) \wedge (\ell_2 \notin J_3)] \\ &\quad \veebar [(k_1 \in J_1) \wedge (\ell_1 \in J_2) \wedge (\ell_2 \in J_3) \wedge (s \notin J_4)] \\ &\quad \veebar [(k_1 \in J_1) \wedge (\ell_1 \in J_2) \wedge (\ell_2 \in J_3) \wedge (s \in J_4) \wedge (s \notin J_5)]. \end{aligned} \tag{14.11}$$

Now we combine (14.7) with each disjunct of (14.11). Each of the five cases split into two, when written in terms of $K_j$ and $L_j$. For instance, the second disjunct,

$$[(k_1 \in J_1) \wedge (\ell_1 \notin J_2)],$$

becomes

$$[(k_1 \in J_1) \wedge (\ell_1 \in K_2)] \veebar [(k_1 \in J_1) \wedge (\ell_1 \in L_2)].$$

These give rise to the p-matrices

$$\left( J_1 \,\middle|\, K_2 \quad I_3 \,\middle|\, I_4 \quad I_5 \right), \quad \left( J_1 \,\middle|\, L_2 \quad I_3 \,\middle|\, I_4 \quad I_5 \right)$$

which will be added to the list generated by the algorithm if they are nonempty. □

The following remark provides variations on Algorithm B useful in calculations.

**14.3. Remark.** Although Algorithm B is stated in the simplest way for the purposes of proof, for it can be simplified further for calculation. It is easy to detect the **noncritical columns** in $P$ and $P'$ for which $I_j = J_j$. For these columns, $K_j = L_j = \emptyset$, and no p-matrices will be produced. So we can skip these columns, and let $j$ vary only over the **critical columns** when applying the algorithm. Furthermore, whatever the current value of $j$ may be, $I_i$ can always be used in the $i$th column of the p-matrix being generated if $i$ is noncritical (since in such cases $I_i = J_i$). Therefore all noncritical columns of the generated p-matrices will agree with $P$, and these can be filled in before beginning the algorithm. Also, we need not move through the critical columns in numerical order, as stated in the algorithm, but can choose a different order. All that is necessary is to keep track of the current **active column** (the current value of $j$), together with the columns that have been **previously active**. When forming a p-matrix, $J_i$ should be entered in each critical column that has been previously active, and $I_i$ in each critical column that has yet to be active. (The active column receives $K_j$ and/or $L_j$, producing one or two matrices.) In examples, we use a downward arrow to indicate the active column in a given p-matrix, and write the p-matrices in the order that they were produced, so that "previously active" columns are those that had downward arrows previously in the list; see (16.2) and (16.4). Changes in the order of the active columns have an effect on the p-matrices obtained, but not on their union. When used to produce a box product by the factoring method, this can affect the *decomposition* obtained for $A$, but not the preferred set $A$ itself, which is determined uniquely by the factoring order according to Theorem 11.8.

**14.4. Theorem** *(Algorithm C). Let $P$ be a p-matrix of points in $\mathbb{N}^{p+q+1}$, given by*

$$P = \left( \begin{array}{c|c|cc} \mathbf{b} & \mathbf{d} & F(\mathbf{k}) & H(\boldsymbol{\ell}) \\ \mathbf{a} & \mathbf{c} & E(\mathbf{k}) & G(\boldsymbol{\ell}) \end{array} \right).$$

*Let $(\boldsymbol{\kappa}; \boldsymbol{\lambda}; \sigma) \in \mathbb{N}^{p+q+1}$ be given, and for any nonnegative integer $j$, let*

$$P[j] = \left( \begin{array}{c|c|cc} \mathbf{b} - j\boldsymbol{\kappa} & \mathbf{d} - j\boldsymbol{\lambda} & F(\mathbf{k} + j\boldsymbol{\kappa}) - j\sigma & H(\boldsymbol{\ell} + j\boldsymbol{\lambda}) - j\sigma \\ \mathbf{a} - j\boldsymbol{\kappa} & \mathbf{c} - j\boldsymbol{\lambda} & E(\mathbf{k} + j\boldsymbol{\kappa}) - j\sigma & G(\boldsymbol{\ell} + j\boldsymbol{\lambda}) - j\sigma \end{array} \right).$$

*Then*

$$P[j] = (P - j(\boldsymbol{\kappa}; \boldsymbol{\lambda}; \sigma)) \cap \mathbb{N}^{(p+q+1)},$$

*that is, $P[j]$ is the result of translating $P$ in $\mathbb{Z}^{p+q+1}$ by the vector $-j(\boldsymbol{\kappa};\boldsymbol{\lambda};\sigma)$ and deleting points with one or more negative entries.*

**14.5. Remark.** Even if $P$ is a p-matrix of well-defined transvectants, this need not be true of $P[j]$. When used in the factoring method, $P[j]$ will always be intersected with a p-matrix $Q$ of well-defined transvectants, and well-definedness comes about automatically at that point.

**Proof.** The elements $(\mathbf{k};\boldsymbol{\ell};s)$ of $P$ satisfy the following inequalities, where $\leq$ is understood to hold componentwise, and to mean $<$ if the right-hand member is $\infty$:

$$\begin{aligned}
&\mathbf{a} \leq \mathbf{k} \leq \mathbf{b},\\
&\mathbf{c} \leq \boldsymbol{\ell} \leq \mathbf{d},\\
&E(\mathbf{k}) \leq s \leq F(\mathbf{k}),\\
&G(\boldsymbol{\ell}) \leq s \leq H(\boldsymbol{\ell}).
\end{aligned}$$

Let the translation through $-j(\boldsymbol{\kappa};\boldsymbol{\lambda};\sigma)$ be written as $\mathbf{k}'=\mathbf{k}-j\boldsymbol{\kappa}$, $\boldsymbol{\ell}'=\boldsymbol{\ell}-j\boldsymbol{\lambda}$, $s'=s-j\sigma$. Substituting these into the inequalities and re-arranging gives

$$\begin{aligned}
&\mathbf{a}-j\boldsymbol{\kappa} \leq \mathbf{k}' \leq \mathbf{b}-j\boldsymbol{\kappa},\\
&\mathbf{c}-j\boldsymbol{\lambda} \leq \boldsymbol{\ell}' \leq \mathbf{d}-j\boldsymbol{\lambda},\\
&E(\mathbf{k}'+j\boldsymbol{\kappa})-\sigma \leq s' \leq F(\mathbf{k}'+j\boldsymbol{\kappa})-\sigma,\\
&G(\boldsymbol{\ell}'+j\boldsymbol{\lambda})-\sigma \leq s' \leq H(\boldsymbol{\ell}'+j\boldsymbol{\lambda})-\sigma.
\end{aligned}$$

Using the notation of (13.4), the set of $(\mathbf{k}';\boldsymbol{\ell}';s')$ in $\mathbb{N}^{p+q+1}$ satisfying these inequalities may be written as

$$\left(\begin{array}{c|c|cc}
\mathbf{k}' & \boldsymbol{\ell}' & \multicolumn{2}{c}{s'}\\
\hline
\mathbf{b}-j\boldsymbol{\kappa} & \mathbf{d}-j\boldsymbol{\lambda} & F(\mathbf{k}'+j\boldsymbol{\kappa})-\sigma & H(\boldsymbol{\ell}'+j\boldsymbol{\lambda})-\sigma\\
\mathbf{a}-j\boldsymbol{\kappa} & \mathbf{c}-j\boldsymbol{\lambda} & E(\mathbf{k}'+j\boldsymbol{\kappa})-\sigma & G(\boldsymbol{\ell}'+j\boldsymbol{\lambda})-\sigma
\end{array}\right).$$

Finally, since $\mathbf{k}'$, $\boldsymbol{\ell}'$, and $s'$ are dummy variables, we are free to replace them everywhere by $\mathbf{k}$, $\boldsymbol{\ell}$, and $s$, and then delete the column headings, as in the discussion of (13.4). This results in the expression given above for $P[j]$. □

**14.6. Theorem** *(Algorithm D). Let $P$ and $Q$ be any two p-matrices, written in the notation of Algorithm B as*

$$P=\left(I_1\cdots I_{p+q+2}\right),\qquad Q=\left(J_1\cdots J_{p+q+2}\right).$$

*Then*

$$P\cap Q=\left((I_1\cap J_1)\cdots(I_{p+q+2}\cap J_{p+q+2})\right).$$

*The intersections of the intervals are given by*

$$[a,b] \cap [c,d] = [\max\{a,c\}, \min\{b,d\}].$$

The proof is obvious and is omitted.

## 15. Factoring method: filtration, reconstruction, and assembly

The box product $\acute{A} \boxtimes_{\boldsymbol{\vartheta}} \grave{A}$, defined in (11.4), will be computed in three steps, called *filtration*, *reconstruction*, and *assembly*. For the following discussion we assume that the admissible primes are known and the factoring order $\boldsymbol{\vartheta} = (\vartheta_1, \dots, \vartheta_t)$ has been decided. The examples in later sections will show how to obtain the primes during the course of the calculation when they are not known in advance.

**The filtration step.** Let $T_0 = \acute{A} \widehat{\boxtimes} \grave{A}$ be the set of all admissible transvectants formed from $\acute{A}$ and $\grave{A}$. The **factoring filtration** of $T_0$ is the finite descending nested sequence of sets

$$T_0 \supset T_1 \supset \cdots \supset T_t = \{1\}, \tag{15.1}$$

where $T_i$ is the set of elements of $T_0$ that are not divisible by any of the first $i$ primes, $\vartheta_1, \dots, \vartheta_i$. The sequence can be defined recursively by

$$T_i = T_{i-1} \setminus T_{i-1}(\vartheta_i).$$

The sequence terminates with $T_t = \{1\}$ because $1 = [\mathbf{0}; \mathbf{0}; 0]$ is the only transvectant that is not divisible by any admissible prime. If $T_i$ is given by a single p-matrix, Algorithm A can be used to compute $T_{i-1}(\vartheta_i)$, followed by Algorithm B to compute the set difference. In general, however, $T_{i-1}$ is a disjoint union of p-matrices. The next lemma shows that the operations $T \mapsto T(\vartheta)$ and $T \mapsto T \setminus T(\vartheta)$ behave nicely with respect to disjoint unions.

**15.1. Lemma.** *If $T = T^1 \sqcup T^2$, then*

$$T(\vartheta) = T^1(\vartheta) \sqcup T^2(\vartheta)$$

*and*

$$T \setminus T(\vartheta) = (T^1 \setminus T^1(\vartheta)) \sqcup (T^2 \setminus T^2(\vartheta)).$$

**Proof.** Any element of $T(\vartheta)$ is an element of $T$, so it belongs to $T^1$ or $T^2$, but not both. Then since the element is divisible by $\vartheta$, it belongs to $T^1(\vartheta)$ or $T^2(\vartheta)$ but not both. The second part is similar. □

The lemma extends to disjoint unions of any number of sets. Suppose that $\acute{A}$ and $\grave{A}$ are given as disjoint unions of b-matrices, as in (8.7). Then $T_0 = \acute{A} \widehat{\boxtimes} \grave{A}$ is given by (8.8), and each term of (8.8) can be written as a p-matrix by (13.1). Therefore $T_0$ is known as a disjoint union of some number $u_0$ of p-matrices:

$$T_0 = P_{01} \sqcup \cdots \sqcup P_{0u_0}. \tag{15.2}$$

Then by Lemma 15.1,

$$T_0(\vartheta_1) = P_{01}(\vartheta_1) \sqcup \cdots \sqcup P_{0u_0}(\vartheta_1);$$

Algorithm A computes each $P_{0j}(\vartheta_1)$ as a single p-matrix, so this is once again a disjoint union of $u_0$ p-matrices (some of which may be empty and can be deleted). Again by Lemma 15.1,

$$T_1 = T_0 \setminus T_0(\vartheta_1) = (P_{01} \setminus P_{01}(\vartheta_1)) \sqcup \cdots \sqcup (P_{0u_0} \setminus P_{0u_0}(\vartheta_1)).$$

Algorithm B computes each term on the right-hand side as a disjoint union of some finite number of p-matrices. Therefore we have

$$T_1 = P_{10} \sqcup \cdots \sqcup P_{1u_1}$$

for some $u_1$. This procedure can be iterated, to obtain

$$T_i = P_{i1} \sqcup \cdots \sqcup P_{iu_i} \tag{15.3}$$

for each $i$. The construction of the factoring filtration ends when $i = t$.

**The reconstruction step.** The $i$th stage of filtration constructed $T_i$ from $T_{i-1}$; the $i$th stage of **reconstruction** reconstructs $T_{i-1}$, in a new form, from $T_i$. Together, these steps amount to factoring each element of $T_{i-1}$ (up to equivalence) into a power of $\vartheta_i$ times an element of $T_{i-1}$. One can either interleave these steps (each stage of filtration being followed by the associated reconstruction), or complete the entire filtration step first and then the reconstructions.

**15.2. Lemma.** *Given $\delta \in T_{i-1}$, there exists a unique $\varepsilon \in T_i$ and $j \in \mathbb{N}$ such that $\delta \equiv \vartheta_i^j \varepsilon$.*

**Proof.** The idea is just to factor $\vartheta_i$ from $\delta$ repeatedly until it cannot be done again. For complete clarity, and to motivate the procedure introduced below, it is best to state this geometrically. Let $\mathbf{p} = (\mathbf{k}; \boldsymbol{\ell}; s) \in \mathcal{C}^*$, and let $\delta = \mathbb{T}(\mathbf{p}) = [\mathbf{k}; \boldsymbol{\ell}; s] \in \mathbb{R}[\mathbf{x}, \mathbf{y}]$ be the associated transvectant. Suppose $\delta \in T_{i-1}$, so that $\delta$ is not divisible by any of $\vartheta_1, \dots, \vartheta_{i-1}$. Then $\mathbf{p}$ is not "subtractible by" any of $\boldsymbol{\pi}_1, \dots, \boldsymbol{\pi}_{i-1}$. That is, none of the integer vectors $\mathbf{p} - \boldsymbol{\pi}_1, \dots, \mathbf{p} - \boldsymbol{\pi}_{i-1}$ belong to $\mathcal{C}$. But some initial segment of the sequence $\mathbf{p}, \mathbf{p} - \boldsymbol{\pi}_i, \mathbf{p} - 2\boldsymbol{\pi}_i, \dots$ (even if only the first term) will belong to $\mathcal{C}$, and Lemma 3.6 implies that this initial segment will in fact belong to $\mathcal{C}^*$, since $\acute{A}$ and $\grave{A}$ are standard. This initial segment must be finite, because the sequence (as a subset of $\mathbb{Z}^{p+q+1}$) will eventually leave $\mathbb{N}^{p+q+1}$. Let $\mathbf{p} - j\boldsymbol{\pi}_i$ be the last element in the initial segment. Since $\pi_i$ cannot be subtracted again from $\mathbf{p} - j\boldsymbol{\pi}_i$, $\varepsilon = \mathbb{T}(\mathbf{p} - j\boldsymbol{\pi}_i)$ cannot be divided by $\vartheta_i = \mathbb{T}\boldsymbol{\pi}_i$. Clearly $\delta = \vartheta_i^j \varepsilon$. □

This lemma gives a *pointwise* reconstruction of $T_{i-1}$ from $T_i$, but does not reconstruct all of $T_{i-1}$ in a finite number of steps. Instead, consider the set $T_{i-1} - j\boldsymbol{\pi}_i \subset \mathbb{Z}^{p+q+1}$ for a fixed $j \in \mathbb{N}$; that is, translate the entire set $T_{i-1}$ through the vector $-j\boldsymbol{\pi}_i$, rather than a single point. Although $T_{i-1} - j\boldsymbol{\pi}_i$ may contain points not in $\mathcal{C}^*$, the intersection

$$(T_{i-1} - j\boldsymbol{\pi}_i) \cap T_i$$

will lie in $\mathcal{C}^*$, and equals the set of all points $\mathbf{q}$ in $T_i$ such that $j\boldsymbol{\pi} + \mathbf{q} \in T_{i-1}$. In terms of transvectants, this is the set of $\varepsilon \in T_i$ such that $\vartheta_i^j \varepsilon$ is equivalent to a transvectant in $T_{i-1}$. The next theorem restates this in terms of p-matrices, so that Algorithms C and D may be used, and also shows that it is only necessary to do this process finitely many times.

**15.3. Theorem.** *Let $P$ and $Q$ be p-matrices of well-defined transvectants, and let $\vartheta = [\boldsymbol{\kappa}; \boldsymbol{\lambda}; \sigma]$ be a well-defined transvectant. Then $P[j] \cap Q$ is the set of transvectants $\varepsilon \in Q$ such that $\vartheta^j \varepsilon \equiv \delta$ for some $\delta \in P$. Furthermore, as $j$ ranges from zero to infinity, $P[j] \cap Q$ takes on only finitely many values, and each of those values is taken on for a finite or infinite interval of integers $j$.*

**Proof.** We have

$$P[j] \cap Q = ((P - j(\boldsymbol{\kappa}; \boldsymbol{\lambda}; \sigma)) \cap \mathbb{N}^{p+q+1}) \cap Q$$
$$= (P - j(\boldsymbol{\kappa}; \boldsymbol{\lambda}; \sigma)) \cap Q$$

since $Q \subset \mathbb{N}^{p+q+1}$. Thus $P[j] \cap Q$ corresponds to the set of transvectants $\varepsilon \in Q$ such that $\vartheta^j \varepsilon$ is a transvectant monomial equivalent to a transvectant in $P$.

For the finiteness, consider first the following one-dimensional situation: given two integer intervals $[a, b]$ and $[c, d]$, and a positive integer $e$, we ask how $[a - je, b - je] \cap [c, d]$ changes as $j$ runs through $j = 0, 1, 2, \dots$. If $b = \infty$, then as soon as $a - je \leq c$ the intersection equals $[c, d]$ for all larger $j$. On the other hand, if $b$ is finite, then as soon as $b - je < c$ the intersection is empty and remains so for larger $j$. So in either case, there are only finitely many values for the intersection. If the intersection $[a - je, b - je] \cap [c, d]$ changes when $j$ is increased by one, it can never return to its previous value, since $[a - je, b - je]$ continues to move to the left. Therefore the set of $j$ for which the intersection takes a particular value must be a finite or infinite interval. These observations transfer immediately to our situation, because the intersection of $P[j]$ and $Q$ is found by intersecting the intervals defined by the corresponding columns. As $j$ increases, $P[j] \cap Q$ changes whenever any of its columns changes, and remains constant whenever all of its columns remain constant. □

By (15.3), $T_{i-1}$ and $T_i$ can be written as disjoint unions of p-matrices:

$$T_{i-1} = \bigsqcup_{r=1}^{u_{(i-1)}} P_{(i-1),r} \quad \text{and} \quad T_i = \bigsqcup_{s=1}^{u_i} P_{is}. \tag{15.4}$$

The following procedure will become clearer through the examples in the next section. Consider $i$ as fixed.

1. For each p-matrix $P_{(i-1),r}$ (for $r = 1, \dots, u_{i-1}$) appearing in $T_{i-1}$, compute $P_{(i-1),r}[j]$ by Algorithm C, taking $\vartheta = \vartheta_i$.
2. Compute the intersections $P_{(i-1),r}[j] \cap P_{is}$ by Algorithm D, for each combination of numerical values for $r$ and $s$, leaving $j$ as a variable.
3. According to Theorem 15.3, for each pair $(r, s)$ there will be finitely many values of the set $R = P_{(i-1),r}[j] \cap P_{is}$, each occurring for an interval $J_R$ of values of $j$. Record the various sets $R$ and their associated intervals $J_R$, for all $r$ and $s$.

4. The reconstruction of $T_{i-1}$ is now given by

$$T_{i-1} = \bigsqcup_R J_R R, \tag{15.5}$$

where each term $J_R R$ in (15.5) consists of an interval $J_R$ of integers representing powers of $\vartheta_i$, multiplied by a p-matrix $R$ of transvectants contained in $P_{is}$. Thus each element of $J_R R$ represents a set of products $\vartheta_i^j \varepsilon$, with $j \in J_R$ and $\varepsilon \in R$. The union in (15.5) is over all of the sets $R$, for all $r$ and $s$ and (intervals of) $j$; this is of course a finite union. Especially in hard problems, it can be more convenient to reconstruct the individual components $P_{(i-1),r}$ of $T_{i-1}$ rather than the entire set $T_{i-1}$. In this case the union in (15.5) will be over the sets $R$ having a fixed $r$.

For instance, in Example 16.2 below, with $i = 4$, we will see that the only nonempty $R$ is $T_4$, with $J_R = [0, 1]$. Therefore (15.5) reduces to

$$T_3 = \begin{bmatrix} 1 \\ 0 \end{bmatrix} T_4 = \begin{bmatrix} 1 \\ 0 \end{bmatrix} \left( \begin{array}{c|cc|cc} \infty & \infty & 0 & k_1 & 2\ell_1 \\ 0 & 0 & & & \end{array} \right).$$

This equation states that every element of $T_4$ can be factored as $\vartheta_4^j (\alpha_1^{k_1}, \beta_1^{\ell_1})^{(s)}$, where $j = 0$ or 1 and where $s = k_1 = 2\ell_1$.

In general, each $P_{is}$ may have several subsets $R$ appearing in (15.5). Those that are proper subsets of $P_{is}$ will be called **special subsets**. Thus, although we have described reconstruction as *reconstruction of* $T_{i-1}$ *from* $T_i$, it would be more accurate to call it *reconstruction of* $T_{i-1}$ *from the p-matrices* $P_{is}$ *of* $T_i$ *and from special subsets of these.*

**The assembly step.** The assembly step is simple and requires no special algorithms. In the simplest examples, Examples 16.1 and 12.2, where no special subsets appear, it is enough to "plug" one reconstruction into another. For instance, in Example 12.2 we obtain

$$\begin{aligned} T_0 &= \begin{bmatrix} \infty \\ 0 \end{bmatrix} T_1 = \begin{bmatrix} \infty \\ 0 \end{bmatrix} \begin{bmatrix} \infty \\ 0 \end{bmatrix} T_2 = \cdots \\ &= \begin{bmatrix} \infty \\ 0 \end{bmatrix} \begin{bmatrix} \infty \\ 0 \end{bmatrix} \begin{bmatrix} \infty \\ 0 \end{bmatrix} \begin{bmatrix} 1 \\ 0 \end{bmatrix} \begin{bmatrix} \infty \\ 0 \end{bmatrix} T_5, \end{aligned}$$

where $T_5 = \{1\}$. The intervals here are (from left to right) intervals of $\vartheta_1, \vartheta_2, \ldots, \vartheta_5$. Considering the meaning of the "product" notation here, it is clear that this just says

$$T_0 = \left[ \begin{array}{c|cc|cc} \infty & \infty & \infty & 1 & \infty \\ 0 & 0 & 0 & 0 & 0 \end{array} \right], \tag{15.6}$$

as written in (16.5). When there are special sets, these must be given their own reconstructions along with the other p-matrices, but this is simple because they are subsets of p-matrices that have already been reconstructed.

## 16. Examples $N_{22}$, $N_{23}$ with two factoring orders, and a first look at $N_{223}$

The first two examples in this section, $N_{22}$ and $N_{23}$, have been worked out earlier, in Examples 10.1 and 10.2, using *ad hoc* methods. The first of these has also been worked out in Example 12.1 by a preliminary version of the factoring method, without the algorithms. Here we work them using the algorithms. (The computation will be long and tedious, but requires only rule-following, since the thought required for the previous methods has gone into the proofs of the algorithms. The hope is that the method should be programmable.) Then we consider $N_{23}$ again with a different factoring order, and then consider $N_{223}$, which received a somewhat unsatisfactory solution in [26].

**16.1. Example.** $N_{22}$, third time.

As in Example 10.1, $p = q = 1$, $\boldsymbol{\alpha} = \alpha$, $\boldsymbol{\beta} = \beta$ (omitting the subscript 1), and $\widehat{\alpha} = \widehat{\beta} = 1$. We do not assume that the primes are known in advance, so we adopt the factoring order (11.2) in the form $\boldsymbol{\vartheta} = (\alpha, \beta, \cdots)$ and proceed with the first two steps. The set of all transvectants is

$$T_0 = \begin{bmatrix} \infty \\ 0 \end{bmatrix} \mathbin{\widehat{\boxtimes}} \begin{bmatrix} \infty \\ 0 \end{bmatrix} = \left( \begin{array}{c|c|cc} \infty & \infty & k & \ell \\ 0 & 0 & 0 & 0 \end{array} \right). \tag{16.1}$$

From (16.1), $E$ and $G$ in $T_0$ are the zero functions, and $F$ and $H$ are identity functions ($F(k) = k$, $H(\ell) = \ell$). We have $\vartheta_1 = \alpha_1 = [1; 0; 0]$, so $\kappa = 1$, $\lambda = 0$, and $\sigma = 0$. Therefore by Algorithm A (14.1),

$$T_0(\vartheta_1) = \left( \begin{array}{c|c|cc} \infty & \infty & k-1 & \ell \\ 1 & 0 & 0 & 0 \end{array} \right).$$

For instance, $F'(k) = \min\{F(\mathbf{k}), \widehat{\boldsymbol{\alpha}}(\mathbf{k} - \boldsymbol{\kappa}) + \sigma\} = \min\{k, 1(k-1) + 0\} = \min\{k, k-1\} = k - 1$.

To compute $T_1$ we use Algorithm B (14.2) together with Remark 14.3 and the filtration step described in Section 15. Comparing $T_0$ and $T_0(\vartheta_1)$, we see that the first and third columns are critical. Taking these columns, in this natural order, as active, and using the arrow notation described in the remark, we have

$$T_1 = T_0 \setminus T_0(\vartheta_1) = \left( \begin{array}{c|c|cc} \downarrow & & & \\ 0 & \infty & k & \ell \\ & 0 & 0 & 0 \end{array} \right) \sqcup \left( \begin{array}{c|c|cc} & & \downarrow & \\ \infty & \infty & \mathrm{k} & \ell \\ 1 & 0 & & 0 \end{array} \right). \tag{16.2}$$

In the first p-matrix, the active (first) column is $K_1 = I_1 \setminus J_1 = [0, \infty] \setminus [1, \infty] = \{0\}$, and the rest is from $T_0$. In the second, the active (third) column is $K_3 = [0, k] \setminus [0, k-1] = \{k\}$, the previously active (first) column is $J_1 = [1, \infty]$, and the rest is from $T_0$. For a strictly algorithmic method, equation (16.2) can be used as it stands, but it is best to simplify the result as follows. The first column of the first p-matrix states that $k = 0$, so the third column is actually $[0, k] = [0, 0] = \{0\}$; this can then be combined with the second p-matrix to produce a single p-matrix for $T_1$:

$$T_1 = \left( \begin{array}{c|c|cc} 0 & \infty & 0 & \ell \\ & 0 & & 0 \end{array} \right) \sqcup \left( \begin{array}{c|c|cc} \infty & \infty & \mathrm{k} & \ell \\ 1 & 0 & & 0 \end{array} \right) = \left( \begin{array}{c|c|cc} \infty & \infty & \mathrm{k} & \ell \\ 0 & 0 & & 0 \end{array} \right). \tag{16.3}$$

Such procedures seem easy, but we have not formalized them algorithmically.

In Remark 14.3, an alternative order is described for the active columns. Using this order, the third column is active first, and we obtain

$$T_1 = T_0 \setminus T_0(\vartheta_1) = \left(\begin{array}{c|c|cc} & & \downarrow & \\ \infty & \infty & \mathrm{k} & \ell \\ 0 & 0 & & 0 \end{array}\right) \sqcup \left(\begin{array}{c|c|cc} \downarrow & & & \\ 0 & \infty & k-1 & \ell \\ & 0 & 0 & 0 \end{array}\right). \tag{16.4}$$

This time the second p-matrix is clearly empty, since $k = 0$ (from the first column) implies that the third column is $[0, k-1] = [0, -1] = \emptyset$. So we delete the last column and obtain the same result as in (16.3) without having to combine columns.

The calculations for $T_1(\beta_1)$ and $T_2$ follow the same pattern, and we obtain

$$T_2 = \left(\begin{array}{c|c|cc} \infty & \infty & k & \ell \\ 0 & 0 & & \end{array}\right).$$

This p-matrix implies that $s = k = \ell$; all of the elements of $T_2$ have the form $(\alpha^s, \beta^s)^{(s)}$. Since we do not know the next prime, we take the first nontrivial case in $T_2$ ($s = k = \ell = 1$) to get $\vartheta_3 = \gamma = (\alpha, \beta)^{(1)}$; this is clearly prime, but it is important to notice that it *must* be prime since (having survived the factoring of $\alpha$ and $\beta$) it is not divisible by these primes, and (because it is the "smallest" element of $T_2$) cannot be divided by any other element of $T_2$. Now use Algorithm A with $\kappa = 1, \lambda = 1, \sigma = 1$ to find that

$$T_2(\vartheta_3) = \left(\begin{array}{c|c|cc} \infty & \infty & k & \ell \\ 1 & 1 & & \end{array}\right).$$

The calculations for this go as follows. First, $a' = \max\{a, \kappa\} = \max\{0, 1\} = 1$, and similarly $c' = 1$. From $T_2$ we see that $E(k) = F(k) = k$ and $G(\ell) = H(\ell) = \ell$. Now we must calculate $E', F', G', H'$ separately. (We cannot assume that $E' = F'$ just because $E = F$, since the formulas in Algorithm A for $E'$ and $F'$ are different.) We find that $E'(k) = \max\{E(k), \sigma\} = \max\{k, 1\} = k$, since $a' = 1$ implies that $k \geq 1$ in $T_2(\vartheta_3)$. Also $F'(k) = \min\{k, 1(k-1)+1\} = k$, so in fact $E' = F'$ after all. The calculations for $G$ and $H$ are similar. Comparing $T_2$ with $T_2(\vartheta_3)$, only the first two columns are critical. We find

$$T_3 = \left(\begin{array}{c|c|cc} \downarrow & & & \\ 0 & \infty & k & \ell \\ & 0 & & \end{array}\right) \sqcup \left(\begin{array}{c|c|cc} & \downarrow & & \\ \infty & 0 & k & \ell \\ 1 & & & \end{array}\right) = \{1\} \sqcup \emptyset = \{1\}.$$

This completes the filtration step.

For the reconstruction step we use Algorithms C and D as outlined in Section 15. Since $\vartheta_1 = \alpha_1 = [1; 0; 0]$, we have $(\kappa; \lambda; \sigma) = (1; 0; 0)$ for the reconstruction of $T_0$ from $T_1$. Algorithm C gives

$$T_0[j] = \left(\begin{array}{c|c|cc} \infty & \infty & k+j & \ell \\ -j & 0 & 0 & 0 \end{array}\right).$$

Algorithm D, and the p-matrix above for $T_1$, now show that

$$T_0[j] \cap T_1 = T_1 \quad \text{for all} \quad j \in [0, \infty].$$

That is, in the notation of (15.5), there is only one set $R = T_1$, and $J_R = [0, \infty]$. Therefore the reconstruction of $T_0$ from $T_1$ is

$$T_0 = \begin{bmatrix} \alpha_1 \\ \hline \infty \\ 0 \end{bmatrix} T_1.$$

Similarly,

$$T_1 = \begin{bmatrix} \beta_1 \\ \hline \infty \\ 0 \end{bmatrix} T_2.$$

Finally, with $(\kappa; \lambda; \sigma) = (1; 1; 1)$ we find

$$T_2[j] = \left( \begin{array}{c|c|cc} \infty & \infty & k & \ell \\ -j & -j & & \end{array} \right)$$

and

$$T_2[j] \cap T_3 = T_3 = \{1\} \quad \text{for all} \quad j \in [0, \infty].$$

Again there is only one $R = \{1\}$ with $J_R = [0, \infty]$, so

$$T_2 = \begin{bmatrix} \gamma_1 \\ \hline \infty \\ 0 \end{bmatrix} \{1\}.$$

This completes the reconstruction.

For the assembly step, we easily see that

$$T_0 = \begin{bmatrix} \alpha_1 \\ \hline \infty \\ 0 \end{bmatrix} \begin{bmatrix} \beta_1 \\ \hline \infty \\ 0 \end{bmatrix} \begin{bmatrix} \gamma_1 \\ \hline \infty \\ 0 \end{bmatrix} \{1\} = \begin{bmatrix} \alpha & \beta & \gamma \\ \hline \infty & \infty & \infty \\ 0 & 0 & 0 \end{bmatrix},$$

in agreement with (10.2). □

**16.2. Example.** $N_{23}$, third time.

To recompute Example 10.2 by the factoring method, we begin with the pre-box product

$$T_0 = \begin{bmatrix} \infty \\ 0 \end{bmatrix} \widehat{\boxtimes} \begin{bmatrix} \infty & \infty \\ 0 & 0 \end{bmatrix} = \left( \begin{array}{c|cc|cc} 1 & 2 & 0 & & \\ \alpha_1 & \beta_1 & \beta_2 & s & s \\ \hline \infty & \infty & \infty & k_1 & 2\ell_1 \\ 0 & 0 & 0 & 0 & 0 \end{array} \right),$$

using $\widehat{\alpha}_1 = 1$, $\widehat{\beta}_1 = 2$, and $\widehat{\beta}_2 = 0$. Again, we do not assume that the primes are known, and we begin with $\boldsymbol{\vartheta} = (\alpha_1, \beta_1, \beta_2, \cdots)$. During the calculation, we will find two more primes, $\vartheta_4 = \gamma_1 = (\alpha_1, \beta_1)^{(1)}$ and $\vartheta_5 = \gamma_2 = (\alpha_1^2, \beta_1)^{(2)}$, so that $\boldsymbol{\vartheta} = (\alpha_1, \beta_1, \beta_2, \gamma_1, \gamma_2)$ as in (11.2). This order will give the same result as in Example 10.2. Later, in Example 16.3, we show that a different factoring order leads to a different block decomposition. Since $\widehat{\beta}_2 = 0$, we could suppress the null input prime $\beta_2$ according to Remark 12.3, by setting $\mathcal{R} = \mathbb{R}[[\beta_2]]$ and omitting the $\beta_2$ column from $T_0$ and all of the following calculations. The result would then appear as (16.7) below.

The calculation of the factoring filtration goes as follows.

$$T_0(\vartheta_1) = \left(\begin{array}{c|cc|cc} \infty & \infty & \infty & k_1 - 1 & 2\ell_1 \\ 1 & 0 & 0 & 0 & 0 \end{array}\right)$$

$$T_1 = \left(\begin{array}{c|cc|cc} & & & \downarrow & \\ \infty & \infty & \infty & & 2\ell_1 \\ 0 & 0 & 0 & k_1 & 0 \end{array}\right) \sqcup \left(\begin{array}{c|cc|cc} \downarrow & & & & \\ & \infty & \infty & k_1 - 1 & 2\ell_1 \\ 0 & 0 & 0 & k_1 & 0 \end{array}\right)$$

$$= \left(\begin{array}{c|cc|cc} & & & \downarrow & \\ \infty & \infty & \infty & & 2\ell_1 \\ 0 & 0 & 0 & k_1 & 0 \end{array}\right) \sqcup \emptyset$$

$$T_1(\vartheta_2) = \left(\begin{array}{c|cc|cc} \infty & \infty & \infty & & 2\ell_1 - 2 \\ 0 & 1 & 0 & k_1 & 0 \end{array}\right)$$

$$T_2 = \left(\begin{array}{c|cc|cc} \infty & \infty & \infty & & 2\ell_1 \\ 0 & 0 & 0 & k_1 & 2\ell_1 - 1 \end{array}\right)$$

$$T_2(\vartheta_3) = \left(\begin{array}{c|cc|cc} \infty & \infty & \infty & & 2\ell_1 \\ 0 & 0 & 1 & k_1 & 2\ell_1 - 1 \end{array}\right)$$

$$T_3 = \left(\begin{array}{c|cc|cc} \infty & \infty & & & 2\ell_1 \\ 0 & 0 & 0 & k_1 & 2\ell_1 - 1 \end{array}\right).$$

At this point we must determine the next prime, so we look for the simplest element of $T_3$, which then cannot be divided by any other element. Taking $k_1 = \ell_1 = 1$ and noting (from the fourth column) that this implies $s = 1$, we have $\vartheta_4 = \gamma_1 = (\alpha, \beta)^{(1)}$. Then we compute

$$T_3(\vartheta_4) = \left(\begin{array}{c|cc|cc} \infty & \infty & & & \\ 1 & 1 & 0 & k_1 & 2\ell_1 - 1 \end{array}\right).$$

The details of the calculation of $T_3(\vartheta_4)$ are as follows. The first step is to identify the functions $E, F, G, H$ in $T_3$:

$$E(\mathbf{k}) = F(\mathbf{k}) = k_1, \quad G(\boldsymbol{\ell}) = 2\ell_1 - 1, \quad H(\boldsymbol{\ell}) = 2\ell_1.$$

Then we use Algorithm 14.1 with $\kappa_1 = \lambda_1 = \sigma = 1$ and $\lambda_2 = 0$ (taken from $\vartheta_4$). This gives

$$E'(\mathbf{k}) = \max\{k_1, 1\} = k_1 \quad \text{since} \quad k_1 \geq 1;$$

$$F'(\mathbf{k}) = \min\{k_1, k_1 - 1 + 1\} = k_1;$$

$$G'(\boldsymbol{\ell}) = \max\{2\ell_1 - 1, 1\} = 2\ell_1 - 1 \quad \text{since} \quad \ell_1 \geq 1;$$
$$H'(\boldsymbol{\ell}) = \min\{2\ell_1, (2,0)\cdot(\ell_1 - 1, 0) + 1\} = 2\ell_1 - 1,$$

which gives the expression for $T_3(\vartheta_4)$.

Since $T_4$ implies $s = k_1 = 2\ell_1$, the simplest element has $\ell_1 = 1$, $s = k_1 = 2$. Therefore $\vartheta_5 = (\alpha^2, \beta)^{(2)}$. Now

$$T_4 = \left(\begin{array}{c|cc|cc} \infty & \infty & 0 & k_1 & 2\ell_1 \\ 0 & 0 & & & \end{array}\right)$$
$$T_4(\vartheta_5) = \left(\begin{array}{c|cc|cc} \infty & \infty & 0 & k_1 & 2\ell_1 \\ 2 & 1 & & & \end{array}\right)$$
$$T_5 = \left(\begin{array}{c|cc|cc} 1 & \infty & 0 & k_1 & 2\ell_1 \\ 0 & 0 & & & \end{array}\right) = \{1\},$$

and the filtration is finished.

The reconstruction phase goes as follows. The interesting step here is the reconstruction of $T_3$ from $T_4$, where, in the notation of (15.5) there are two sets $R$ (one empty) and two intervals $J_R$, namely $[0, 1]$ and $[2, \infty]$. Also, in the first expression for $T_3[j] \cap T_4$, it is convenient to add a third column, temporarily, to the two final columns that are ordinarily used to control $s$.

$$T_0[j] = \left(\begin{array}{c|cc|cc} \infty & \infty & \infty & k_1 + j & 2\ell \\ -j & 0 & 0 & 0 & 0 \end{array}\right)$$
$$T_1 \cap T_0[j] = T_1 \qquad j \in [0, \infty]$$
$$T_0 = \begin{bmatrix} \infty \\ 0 \end{bmatrix} T_1$$
$$T_1[j] = \left(\begin{array}{c|cc|cc} \infty & \infty & \infty & k_1 & 2\ell_1 \\ 0 & -j & 0 & & 0 \end{array}\right)$$
$$T_1[j] \cap T_2 = T_2 \qquad j \in [0, \infty]$$
$$T_1 = \begin{bmatrix} \infty \\ 0 \end{bmatrix} T_2$$
$$T_2[j] = \left(\begin{array}{c|cc|cc} \infty & \infty & \infty & k_1 & 2\ell_1 \\ 0 & 0 & -j & & 2\ell_1 - 1 \end{array}\right)$$
$$T_2[j] \cap T_3 = T_3 \qquad j \in [0, \infty]$$
$$T_2 = \begin{bmatrix} \infty \\ 0 \end{bmatrix} T_3$$
$$T_3[j] = \left(\begin{array}{c|cc|cc} \infty & \infty & 0 & k_1 & 2\ell_1 + j \\ -j & -j & & & 2\ell_1 + j - 1 \end{array}\right)$$
$$T_3[j] \cap T_4 = \left(\begin{array}{c|cc|ccc} \infty & \infty & 0 & k_1 & 2\ell_1 & 2\ell_1 + j \\ 0 & 0 & & & & 2\ell_1 + j - 1 \end{array}\right)$$
$$= \begin{cases} \left(\begin{array}{c|cc|cc} \infty & \infty & 0 & k_1 & 2\ell_1 \\ 0 & 0 & & & \end{array}\right) = T_4 & j \in [0, 1] \\ \emptyset & j \in [2, \infty] \end{cases}$$

$$T_3 = \begin{bmatrix} 1 \\ 0 \end{bmatrix} T_4$$

$$T_4[j] = \left( \begin{array}{c|cc|cc} \infty & \infty & 0 & k_1 & 2\ell_1 \\ -2j & -j & & & \end{array} \right)$$

$$T_4[j] \cap T_5 = \left( \begin{array}{c|cc|cc} 0 & 0 & 0 & 0 & 0 \end{array} \right) = T_5 = \{1\} \qquad j \in [0, \infty]$$

$$T_4 = \begin{bmatrix} \infty \\ 0 \end{bmatrix} T_5$$

The assembly stage is simple, since the only "special subset" that arises is the empty set associated with $[2, \infty]$ in the reconstruction of $T_3$, which (of course) does not appear in the expression for $T_3$. The result of the assembly was found in (15.6):

$$T_0 = \left[ \begin{array}{c|cc|cc} \infty & \infty & \infty & 1 & \infty \\ 0 & 0 & 0 & 0 & 0 \end{array} \right], \tag{16.5}$$

where the columns represent powers of $\alpha, \beta_1, \beta_2, \gamma_1, \gamma_2$. This is the same as the "decorated" block (10.6). The Stanley decomposition is

$$\mathcal{I}_{23} = \mathbb{R}[[\alpha_1, \beta_1, \beta_2, \gamma_2]] \oplus \mathbb{R}[[\alpha_1, \beta_1, \beta_2, \gamma_2]]\gamma_1, \tag{16.6}$$

as previously found in (10.7). If we had suppressed the null input prime $\beta_2$ at the beginning of the calculation, putting $\mathcal{R} = \mathbb{R}[[\beta_2]]$, the third column of (16.5) would be missing, and the Stanley decomposition would be written

$$\mathcal{I}_{23} = \mathcal{R}[[\alpha_1, \beta_1, \gamma_2] \oplus \mathcal{R}[[\alpha_1, \beta_1, \gamma_2]]\gamma_1. \tag{16.7}$$

Note that $\gamma_2$ is a **null output prime**, that is, a prime of weight zero occurring in $\boldsymbol{\gamma}$ (rather than $\boldsymbol{\alpha}$ or $\boldsymbol{\beta}$), and happens also to appear in every coefficient ring in the Stanley decomposition. (This is not always the case for null outputs, as it is for null inputs. See (16.15).) So we could introduce $\mathcal{R}^+ = \mathbb{R}[[\beta_2, \gamma_2]]$ and write the Stanley decomposition as

$$\mathcal{I}_{23} = \mathcal{R}^+[[\alpha_1, \beta_1] \oplus \mathcal{R}^+[[\alpha_1, \beta_1]]\gamma_1. \qquad \square \tag{16.8}$$

**16.3. Example.** $N_{23}$, third time, with a different factoring order.

To see that the factoring order matters in the case of $N_{23}$, consider these two factorizations:

$$(\alpha_1^2, \beta_1^2)^{(2)} \equiv \beta_1(\alpha_1^2, \beta_1)^{(2)} = \beta_1\gamma_2,$$
$$(\alpha_1^2, \beta_1^2)^{(2)} \equiv (\alpha_1, \beta_1)^{(1)2} = \gamma_1^2. \tag{16.9}$$

The first is obtained if either $\beta_1$ or $\gamma_2$ is factored first, the second if $\gamma_1$ is factored first. Another way to say this is that the first factorization is lex maximal unless $\gamma_1$ occurs before both $\beta_1$ and $\gamma_2$ in $\boldsymbol{\vartheta}$, in which case the second factorization is lex maximal.

The factoring order $\vartheta = (\beta_2, \gamma_1, \gamma_2, \alpha_1, \beta_1)$ gives the following factoring filtration (omitting three empty p-matrices in $T_2$, and two in $T_3$, that are produced by Algorithm B):

$$
\begin{aligned}
T_0 &= \left(\begin{array}{c|cc|cc} \alpha_1 & \beta_1 & \beta_2 & s & s \\ \hline \infty & \infty & \infty & k_1 & 2\ell_1 \\ 0 & 0 & 0 & 0 & 0 \end{array}\right) \\
T_1 &= \left(\begin{array}{c|cc|cc} \infty & \infty & 0 & k_1 & 2\ell_1 \\ 0 & 0 & & 0 & 0 \end{array}\right) \\
T_2 &= \left(\begin{array}{c|cc|cc} \infty & \infty & 0 & 0 & 0 \\ 0 & 0 & & & \end{array}\right) \sqcup \left(\begin{array}{c|cc|cc} \infty & \infty & 0 & k_1 & 2\ell_1 \\ 1 & 1 & & 1 & \end{array}\right) = P_{21} \sqcup P_{22} \\
T_3 &= \left(\begin{array}{c|cc|cc} \infty & \infty & 0 & 0 & 0 \\ 0 & 0 & & & \end{array}\right) \\
T_4 &= \left(\begin{array}{c|cc|cc} 0 & \infty & 0 & 0 & 0 \\ & 0 & & & \end{array}\right) \\
T_5 &= \{1\}.
\end{aligned}
$$

The reconstruction stage is more interesting, since it is the first example with a (nontrivial) special set. The first two steps are easy:

$$
T_0 = \begin{bmatrix} \beta_2 \\ \hline \infty \\ 0 \end{bmatrix} T_1, \qquad T_1 = \begin{bmatrix} \gamma_1 \\ \hline \infty \\ 0 \end{bmatrix} T_2.
$$

This already shows that there is no limit on $\gamma_1$, in contrast to Example 12.2 and in agreement with the second factoring in (16.9). For the next step,

$$
\begin{aligned}
P_{21}[j] &= \left(\begin{array}{c|cc|cc} \infty & \infty & 0 & -2j & -2j \\ -2j & -j & & & \end{array}\right) \\
P_{22}[j] &= \left(\begin{array}{c|cc|cc} \infty & \infty & 0 & k_1 & -2\ell_1 \\ 1-2j & 1-j & & 1-2j & \end{array}\right).
\end{aligned}
$$

Then $T_2[j] = P_{21}[j] \sqcup P_{22}[j]$, and it follows that

$$
T_2[j] \cap T_3 = \begin{cases} \left(\begin{array}{c|cc|cc} \infty & \infty & 0 & 0 & 0 \\ 0 & 0 & & & \end{array}\right) \sqcup \emptyset & j = 0 \\ \emptyset \sqcup \left(\begin{array}{c|cc|cc} \infty & 0 & 0 & 0 & 0 \\ 0 & & & & \end{array}\right) & j \geq 1. \end{cases}
$$

That is,

$$
T_2[j] \cap T_3 = \begin{cases} T_3 & j = 0 \\ \left(\begin{array}{c|cc|cc} \infty & 0 & 0 & 0 & 0 \\ 0 & & & & \end{array}\right) & j \geq 1. \end{cases}
$$

Here the first case ($T_3$) is associated with the interval $J = [0, 0] = \{0\}$ of $j$, and the second case (a special subset of $T_3$) is associated with $J = [1, \infty]$. The reconstruction of $T_2$ from $T_3$ is therefore given by

$$T_2 = \left[\begin{array}{c} \gamma_2 \\ \hline 0 \end{array}\right] T_3 \sqcup \left[\begin{array}{c} \gamma_2 \\ \hline \infty \\ 1 \end{array}\right] \left( \begin{array}{c|cc|cc} \infty & \multirow{2}{*}{0} & \multirow{2}{*}{0} & \multirow{2}{*}{0} & \multirow{2}{*}{0} \\ 0 & & & & \end{array} \right).$$

The remaining steps of reconstruction are easy, giving

$$T_3 = \left[\begin{array}{c} \alpha_1 \\ \hline \infty \\ 0 \end{array}\right] T_1, \qquad T_1 = \left[\begin{array}{c} \beta_1 \\ \hline \infty \\ 0 \end{array}\right] T_5.$$

Finally, the assembly gives

$$T_0 \equiv \left[\begin{array}{ccccc} \beta_2 & \gamma_1 & \gamma_2 & \alpha_1 & \beta_1 \\ \hline \infty & \infty & \multirow{2}{*}{0} & \infty & \infty \\ 0 & 0 & & 0 & 0 \end{array}\right] \sqcup \left[\begin{array}{ccccc} \beta_2 & \gamma_1 & \gamma_2 & \alpha_1 & \beta_1 \\ \hline \infty & \infty & \infty & \infty & \multirow{2}{*}{0} \\ 0 & 0 & 1 & 0 & \end{array}\right].$$

Written as a Stanley decomposition, this becomes

$$\mathcal{J}_{23} = \mathbb{R}[[\alpha_1, \beta_1, \beta_2, \gamma_1]] \oplus \mathbb{R}[[\alpha_1, \beta_2, \gamma_1, \gamma_2]]\gamma_2. \qquad \square$$

**16.4. Example.** $N_{223}$ by Algorithms A–D.

This example was first done in [26], but the result obtained there by the expansion method was long and awkward. The solution obtained here is much shorter; see also Remark 3.8. A block decomposition for $\mathcal{J}_{223}$ can be found in two ways, as $\mathcal{J}_{22} \boxtimes \mathcal{J}_3$ or as $\mathcal{J}_2 \boxtimes \mathcal{J}_{23}$. Here we use the first of these. The starting point is the pre-box product (decorated with degrees and weights)

$$T_0 = A_{22} \mathbin{\widehat{\boxtimes}} A_3 = \left[\begin{array}{ccc} 1 & 1 & 2 \\ 1 & 1 & 0 \\ \alpha_1 & \alpha_2 & \alpha_3 \\ \hline \infty & \infty & \infty \\ 0 & 0 & 0 \end{array}\right] \mathbin{\widehat{\boxtimes}} \left[\begin{array}{cc} 1 & 2 \\ 2 & 0 \\ \beta_1 & \beta_2 \\ \hline \infty & \infty \\ 0 & 0 \end{array}\right]. \tag{16.10}$$

We immediately suppress the null input primes $\alpha_3$ and $\beta_2$, by Remark 12.3, writing $\mathcal{R} = \mathbb{R}[[\alpha_3, \beta_2]]$ and replacing (16.10) by

$$T_0 = \left[\begin{array}{cc} \alpha_1 & \alpha_2 \\ \hline \infty & \infty \\ 0 & 0 \end{array}\right] \mathbin{\widehat{\boxtimes}} \left[\begin{array}{c} \beta_1 \\ \hline \infty \\ 0 \end{array}\right]. \tag{16.11}$$

The set of admissible primes, in the order in which they are discovered during the calculation, is

$$\boldsymbol{\vartheta} = (\alpha_1, \alpha_2, \beta_1, \gamma_1, \gamma_2, \gamma_3, \gamma_4, \gamma_5),$$

where $\gamma_1 = (\alpha_1, \beta_1)^{(1)}$, $\gamma_2 = (\alpha_2, \beta_1)^{(1)}$, $\gamma_3 = (\alpha_1^2, \beta_1)^{(2)}$, $\gamma_4 = (\alpha_1\alpha_2, \beta_1)^{(2)}$, and $\gamma_5 = (\alpha_2^2, \beta_1^2)^{(2)}$. Notice that since $\alpha_3$ and $\beta_2$ are suppressed, they do not appear in the list of primes. In particular, $\vartheta_3 = \beta_1$, not $\alpha_3$.

The result found from the algorithms is as follows. Degrees, weights, and labels given in the first b-matrix apply throughout.

$$\begin{aligned} &\left[\begin{array}{cc|c|ccccc} 1 & 1 & 1 & 2 & 2 & 3 & 3 & 3 \\ 1 & 1 & 2 & 1 & 1 & 0 & 0 & 0 \\ \alpha_1 & \alpha_2 & \beta_1 & \gamma_1 & \gamma_2 & \gamma_3 & \gamma_4 & \gamma_5 \\ \hline \begin{matrix}\infty\\0\end{matrix} & 0 & \begin{matrix}\infty\\0\end{matrix} & \begin{matrix}1\\0\end{matrix} & 0 & \begin{matrix}\infty\\0\end{matrix} & \begin{matrix}1\\0\end{matrix} & \begin{matrix}\infty\\0\end{matrix} \end{array}\right] \\ &\sqcup \left[\begin{array}{cc|c|ccccc} \begin{matrix}\infty\\0\end{matrix} & 0 & \begin{matrix}\infty\\0\end{matrix} & 0 & 1 & 0 & 0 & \begin{matrix}\infty\\0\end{matrix} \end{array}\right] \\ &\sqcup \left[\begin{array}{cc|c|ccccc} \begin{matrix}\infty\\0\end{matrix} & \begin{matrix}\infty\\1\end{matrix} & \begin{matrix}\infty\\0\end{matrix} & 0 & 0 & 0 & 0 & \begin{matrix}\infty\\0\end{matrix} \end{array}\right] \\ &\sqcup \left[\begin{array}{cc|c|ccccc} \begin{matrix}\infty\\0\end{matrix} & \begin{matrix}\infty\\1\end{matrix} & \begin{matrix}\infty\\0\end{matrix} & 0 & 1 & 0 & 0 & \begin{matrix}\infty\\0\end{matrix} \end{array}\right]. \end{aligned} \tag{16.12}$$

But the fourth b-matrix can be combined with the second to give

$$\begin{aligned} &\left[\begin{array}{cc|c|ccccc} \begin{matrix}\infty\\0\end{matrix} & 0 & \begin{matrix}\infty\\0\end{matrix} & \begin{matrix}1\\0\end{matrix} & 0 & \begin{matrix}\infty\\0\end{matrix} & \begin{matrix}1\\0\end{matrix} & \begin{matrix}\infty\\0\end{matrix} \end{array}\right] \\ &\sqcup \left[\begin{array}{cc|c|ccccc} \begin{matrix}\infty\\0\end{matrix} & \begin{matrix}\infty\\0\end{matrix} & \begin{matrix}\infty\\0\end{matrix} & 0 & 1 & 0 & 0 & \begin{matrix}\infty\\0\end{matrix} \end{array}\right] \\ &\sqcup \left[\begin{array}{cc|c|ccccc} \begin{matrix}\infty\\0\end{matrix} & \begin{matrix}\infty\\1\end{matrix} & \begin{matrix}\infty\\0\end{matrix} & 0 & 0 & 0 & 0 & \begin{matrix}\infty\\0\end{matrix} \end{array}\right]. \end{aligned} \tag{16.13}$$

Written as a Stanley decomposition with the coefficient rings distributing over the direct sums, this becomes

$$\mathcal{R}[\alpha_1, \beta_1, \gamma_3, \gamma_5](1 \oplus \gamma_1 \oplus \gamma_4 \oplus \gamma_1\gamma_4) \oplus \mathcal{R}[[\alpha_1, \alpha_2, \beta_1, \gamma_5]](\alpha_2 \oplus \gamma_2). \tag{16.14}$$

Since the primes $\alpha_1, \beta_1, \gamma_5$ appear in every coefficient ring, we can put $\mathcal{R}^+ = \mathbb{R}[[\alpha_1, \alpha_3, \beta_1, \beta_2, \gamma_5]]$ and abbreviate this further as

$$\mathcal{R}^+[[\gamma_3]](1 \oplus \gamma_1 \oplus \gamma_4 \oplus \gamma_1\gamma_4) \oplus \mathcal{R}^+[[\alpha_2]](\alpha_2 \oplus \gamma_2). \tag{16.15}$$

(Notice that $\gamma_3$ and $\gamma_4$, although they are null output primes, do not occur in every coefficient ring.)

Now we present the first few steps of the calculation, leaving the rest to the reader. This time we interleave the filtration and reconstruction steps. With $\alpha_3$ and $\beta_2$ suppressed, and $\boldsymbol{\vartheta}$ as above, we have

$$
\begin{aligned}
T_0 &= \left(\begin{array}{cc|c|cc}
\alpha_1 & \alpha_2 & \beta_1 & s & s \\ \hline
\infty & \infty & \infty & k_1+k_2 & 2\ell_1 \\
0 & 0 & 0 & 0 & 0
\end{array}\right) \\
T_0(\vartheta_1) &= \left(\begin{array}{cc|c|cc}
\infty & \infty & \infty & k_1+k_2-1 & 2\ell_1 \\
1 & 0 & 0 & 0 & 0
\end{array}\right) \\
T_1 &= \left(\begin{array}{cc|c|cc}
 & & & \downarrow & \\
\infty & \infty & \infty & k_1+k_2 & 2\ell_1 \\
0 & 0 & 0 & & 0
\end{array}\right) \sqcup \left(\begin{array}{cc|c|cc}
\downarrow & & & & \\
0 & \infty & \infty & k_2-1 & 2\ell_1 \\
 & 0 & 0 & 0 & 0
\end{array}\right) \\
&= P_{11} \sqcup P_{12}.
\end{aligned}
$$

In the last equation we have used the notation of (15.3). Notice that in the fourth column of the last p-matrix, $k_1$ drops out because the first column implies $k_1 = 0$. Also notice that this second p-matrix is not empty (as it was in the similar case of $T_1$ in Example 16.2). For the reconstruction of $T_0$ from $T_1$ we have $\boldsymbol{\kappa} = (1, 0)$, $\boldsymbol{\lambda} = 0$, $\sigma = 0$,

$$
\begin{aligned}
T_0[j] &= \left(\begin{array}{cc|c|cc}
\infty & \infty & \infty & k_1+k_2+j & 2\ell_1 \\
-j & 0 & 0 & 0 & 0
\end{array}\right) \\
T_0[j] \cap P_{11} &= P_{11} \quad j \in [0, \infty] \\
T_0[j] \cap P_{12} &= P_{12} \quad j \in [0, \infty] \\
T_0 &\equiv \left[\begin{array}{c}
\alpha_1 \\ \hline \infty \\ 0
\end{array}\right] T_1.
\end{aligned}
$$

This means that each element of $T_0$ can be factored (up to equivalence) in one and only one way as a power of $\alpha_1$ times an element of $T_1$.

Next we have

$$
\begin{aligned}
T_1(\vartheta_2) = &\left(\begin{array}{cc|c|cc}
\infty & \infty & \infty & k_1+k_2-1 & 2\ell_1 \\
0 & 1 & 0 & k_1+k_2 & 0
\end{array}\right) \\
&\sqcup \left(\begin{array}{cc|c|cc}
0 & \infty & \infty & k_2-1 & 2\ell_1 \\
 & 1 & 0 & 0 & 0
\end{array}\right).
\end{aligned}
$$

The two p-matrices here are subsets of $P_{11}$ and $P_{12}$ respectively, and the first of these is empty. Therefore

$$
\begin{aligned}
T_2 &= (P_{11} \setminus \emptyset) \sqcup \left(P_{12} \setminus \left(\begin{array}{cc|c|cc}
0 & \infty & \infty & k_2-1 & 2\ell_1 \\
 & 1 & 0 & 0 & 0
\end{array}\right)\right) \\
&= P_{11} \sqcup \left(\begin{array}{cc|c|cc}
 & \downarrow & & & \\
0 & 0 & \infty & -1 & 2\ell_1 \\
 & & 0 & 0 & 0
\end{array}\right) \\
&= P_{11} \setminus \emptyset = P_{11}.
\end{aligned}
$$

Notice that, comparing $P_{12}$ with the p-matrix to be subtracted from it, there is only one critical column, so there is only one p-matrix (with $\downarrow$ in the active position) when the difference is computed. This p-matrix has $-1$ in the $F$ position because $k_2 = 0$. For the reconstruction of $T_1$ from $T_2$ we have $\boldsymbol{\kappa} = (0, 1)$, $\boldsymbol{\lambda} = 0$, $\sigma = 0$,

$$P_{11}[j] = \left(\begin{array}{cc|c|cc} \infty & \infty & \infty & k_1 + k_2 + j & 2\ell_1 \\ 0 & -j & 0 & & 0 \end{array}\right)$$

$$P_{11}[j] \cap T_2 = \begin{cases} T_2 & j = 0 \\ \emptyset & j \in [1, \infty] \end{cases}$$

$$P_{12}[j] = \left(\begin{array}{cc|c|cc} 0 & \infty & \infty & k_2 - 1 + j & 2\ell_1 \\ & -j & 0 & 0 & 0 \end{array}\right)$$

$$P_{12}[j] \cap T_2 = \begin{cases} \emptyset & j = 0 \\ S_2 = \left(\begin{array}{cc|c|cc} 0 & \infty & \infty & k_2 & 2\ell_1 \\ & 0 & 0 & & 0 \end{array}\right) & j \in [1, \infty]. \end{cases}$$

It is now clear that there are four subsets involved in the reconstruction of $T_1$ from $T_2$; these sets $R$ and their associated intervals $J_R$, for use in (15.5), are as follows:

1. $R = T_2$ with $J_R = \{0\}$,
2. $R = \emptyset$ with $J_R = [1, \infty]$,
3. $R = \emptyset$ with $J_R = \{0\}$,
4. $R = S_2$ with $J_R = [1, \infty]$.

Since the empty sets can be omitted, equation (15.5) reduces to

$$\begin{aligned} T_1 &= \begin{bmatrix} \alpha_2 \\ \hline 0 \end{bmatrix} T_2 \sqcup \begin{bmatrix} \alpha_2 \\ \hline \infty \\ 1 \end{bmatrix} S_2 \\ &= T_2 \sqcup \begin{bmatrix} \alpha_2 \\ \hline \infty \\ 1 \end{bmatrix} S_2. \end{aligned}$$

**16.5. Remark.** We showed every detail of this reconstruction step in order to illustrate equation (15.5). In this instance we could have skipped some steps by noticing that $P_{11} = T_2$, so it is not actually necessary to compute $P_{11}[j]$; the final conclusion can be reached by doing only the calculations involving $P_{12}$. The notation $S_2$ for the required special subset of $T_2$ is an *ad hoc* notation. We have not attempted to find a systematic notation for the special subsets that may arise in a given calculation. In harder problems there will often be nested families of special subsets within a particular $P_{is}$.

The remaining steps are similar and are left to the reader. After the *method of tokens* is developed in Part 2, it will be possible to do this example more easily, and full details will be given. (The final answer is the same as above.) □

# Part 3. Boosting

## 17. Boosting to equivariants

The goal of boosting is, given $N$, a normal form style, and a Stanley decomposition for $\mathcal{I}$, to construct a Stanley decomposition for $\mathcal{E}$. In equation (1.6), it was suggested that this process is a modified box product of the form $\mathcal{E} = \mathcal{I} \boxtimes \mathbb{R}^n$. Here $\mathcal{I}$ is an algebra of invariants, expressed as usual by a block decomposition $A$. The right-hand factor $\mathbb{R}^n$ is just Euclidean $n$-space, which is a vector space but not an algebra; this will result in some minor differences between the box product for boosting and the box product for invariants. Before we can proceed, an $\mathfrak{sl}_2$ triad for $\mathbb{R}^n$ must be specified so that transvectants can be defined. There are already two natural $\mathfrak{sl}_2$ triads on $\mathbb{R}^n$, namely $(N, M, H)$ and $(M^*, N^*, H)$ (see Section 2), but neither of these is suitable, and we must instead use $(-N^*, -M^*, -H)$. Understanding the reason for this is the first issue that must be addressed. Next, one should expect that the subspace $\ker(-N^*) = \ker N^* \subset \mathbb{R}^n$ will play the role of the "invariants" of the triad $(-N^*, -M^*, -H)$; but this, being only a vector space, has only a vector space basis $\mathbf{B}$, not a Hilbert basis $\boldsymbol{\beta}$. Since nothing can be built from $\boldsymbol{\beta}$ by multiplication of its elements (which is not defined), there is no need for a block decomposition for $\ker N^*$, and the second factor of the box product will simply be $\mathbf{B}$. Thus the box product that produces a description of $\mathcal{E}$ will have the form $A \boxtimes \mathbf{B}$ (rather than $\acute{A} \boxtimes \grave{A}$). The full table of notations for this section, replacing the one in Section 6, will be as follows.

| left input | right input | output |
|---|---|---|
| $\mathbf{x} \in \mathbb{R}^n$ | | |
| $(X, Y, Z) = (M^*, N^*, H)$ | | |
| $\mathbb{R}[[\mathbf{x}]]$ | $\mathbb{R}^n$ | $\mathcal{V}^n$ |
| $(\mathcal{X}, \mathcal{Y}, \mathcal{Z})$ | $(-N^*, -M^*, -H)$ | $(\mathsf{X}, \mathsf{Y}, \mathsf{Z})$ |
| $\mathcal{I} = \ker \mathcal{X}$ | $\ker N^*$ | $\mathcal{E} = \ker \mathsf{X}$ |
| $\boldsymbol{\alpha} = (\alpha_1, \dots, \alpha_p)$ | $\mathbf{B} = (\mathbf{b}_1, \dots, \mathbf{b}_q)$ | |
| $\mathbf{k} \in \mathbb{N}^p$ | | |
| $\boldsymbol{\alpha}^{\mathbf{k}} = \alpha_1^{k_1} \cdots \alpha_p^{k_p}$ | | |
| $A \subset \mathbb{R}[\boldsymbol{\alpha}]$ | | $A \boxtimes \mathbf{B}$ |
| $A = B^1 \sqcup \cdots \sqcup B^\mu$ | | |

The left-hand input column uses the notation of Sections 2–5. The information in this column is assumed to be given, either by the results of Section 4 or as the result of a previous box product calculation. (In that case the left input column of this table would be the same as the output column of the table in Section 6, except that $n + m$ would be renamed as $n$, $(\mathbf{x}, \mathbf{y})$ as $\mathbf{x}$, and $\boldsymbol{\vartheta}$ as $\boldsymbol{\alpha}$.)

The output column is easy to understand, because the goal of boosting is exactly to obtain this column. We know from Section 2 that the space $\mathcal{V}^n$ of formal vector fields on $\mathbb{R}^n$ carries a natural $\mathfrak{sl}_2$ triad $(\mathsf{X}, \mathsf{Y}, \mathsf{Z})$ given by (2.9), and that the equivariants that constitute the normal form are given by $\mathcal{E} = \ker \mathsf{X}$. All that remains to be explained in this column is the definition of the box product $A \boxtimes \mathbf{B}$. But the middle column must first be discussed in greater detail.

Any formal vector field in $\mathcal{V}^n$ can be written uniquely as

$$\mathbf{v}(\mathbf{x}) = f_1(\mathbf{x})\mathbf{e}_1 + \cdots + f_n(\mathbf{x})\mathbf{e}_n,$$

where $f_i \in \mathbb{R}[[\mathbf{x}]]$. Thus, a formal basis for $\mathcal{V}^n$ consists of the products of the monomials $\mathbf{x}^{\mathbf{m}}$ (a formal basis for $\mathbb{R}[[\mathbf{x}]]$) and the standard basis vectors for $\mathbb{R}^n$. This explains the third line of the table (reading across).

**17.1. Remark.** More technically, let $\psi : \mathbb{R}[[\mathbf{x}]] \times \mathbb{R}^n \to \mathcal{V}^n$ be defined by $\psi(f, \mathbf{w}) = f\mathbf{w}$. Then $\psi$ induces an isomorphism $\widetilde{\psi} : \mathbb{R}[[\mathbf{x}]] \otimes \mathbb{R}^n \to \mathcal{V}^n$. That is, the ordinary product of a scalar function and a fixed vector satisfies the requirements of a tensor product. (Compare Remark 6.1).

To establish that $(-N^*, -M^*, -H)$ is the triad on $\mathbb{R}^n$ that will correctly produce $(\mathsf{X}, \mathsf{Y}, \mathsf{Z})$ in column 3, we begin from $(\mathsf{X}, \mathsf{Y}, \mathsf{Z})$ and compute its effect on $f\mathbf{e}_i$ using (2.3) and (2.13) to find that

$$\begin{aligned} \mathsf{X}(f\mathbf{e}_i) &= (\mathcal{X}f)\mathbf{e}_i + (f)(-N^*)\mathbf{e}_i, \\ \mathsf{Y}(f\mathbf{e}_i) &= (\mathcal{Y}f)\mathbf{e}_i + (f)(-M^*)\mathbf{e}_i, \\ \mathsf{Z}(f\mathbf{e}_i) &= (\mathcal{Z}f)\mathbf{e}_i + (f)(-H)\mathbf{e}_i. \end{aligned} \tag{17.1}$$

With $f = 1$ this gives $\mathsf{X}\mathbf{e}_i = -N^*\mathbf{e}_i$, $\mathsf{Y}\mathbf{e}_i = -M^*\mathbf{e}_i$, and $\mathsf{Z}\mathbf{e}_i = -Z^*\mathbf{e}_i$.

**17.2. Remark.** The tensor product of the Lie algebra representations of $\mathfrak{sl}_2$ defined by $(\mathcal{X}, \mathcal{Y}, \mathcal{Z})$ and $(-N^*, -M^*, -H)$ is given in tensor notation by $(\mathcal{X} \otimes I + I \otimes (-N^*), \mathcal{Y} \otimes I + I \otimes (-M^*), \mathcal{Z} \otimes I_I \otimes (-H))$. This reduces to the right-hand side of (17.1) under the isomorphism of Remark 17.1.

Continuing down the middle column, ker $(-N^*) =$ ker $N^*$ is the top space of $\mathbb{R}^n$ under $(-N^*, -M^*, -H)$, and it has a natural basis $\mathbf{B} = (\mathbf{b}_1, \dots, \mathbf{b}_q)$, where $q$ is the number of Jordan blocks in $N$ and each $\mathbf{b}_j$ equals the Euclidean basis vector $\mathbf{e}_r$, where $r$ is the row number of the bottom row of the $j$-th block. This is illustrated in Remark 17.3. The weight of $\mathbf{e}_r$ (its eigenvalue with respect to $-H$) is one less than the block size of the block with bottom row $r$. (See Theorem 2.3, $\ell_j = \widehat{\mathbf{v}}_j + 1$.)

**17.3. Remark.** To compare the triads $(M^*, N^*, H)$ and $(-N^*, -M^*, -H)$, consider the case of $N_{23}$. We first get $(N_{23}, M_{23}, H_{23})$ as follows:

$$N = \begin{bmatrix} 0 & 1 & & & \\ 0 & 0 & & & \\ & & 0 & 1 & 0 \\ & & 0 & 0 & 1 \\ & & 0 & 0 & 0 \end{bmatrix}, \quad M = \begin{bmatrix} 0 & 0 & & & \\ 1 & 0 & & & \\ & & 0 & 0 & 0 \\ & & 2 & 0 & 0 \\ & & 0 & 2 & 0 \end{bmatrix}, \quad H = \begin{bmatrix} 1 & & & & \\ & -1 & & & \\ & & 2 & & \\ & & & 0 & \\ & & & & -2 \end{bmatrix}.$$

Then for $(M^*, N^*, H)$ the top weight vectors (in ker $M^*$) are $\mathbf{e}_1$ and $\mathbf{e}_3$, whose subscripts are the row numbers of the top rows of the Jordan blocks; the weights of these chain tops (eigenvalues of $H$) are 1 and 2; and the chains are $\{\mathbf{e}_1, N^*\mathbf{e}_1 = \mathbf{e}_2\}$ and $\{\mathbf{e}_3, \mathbf{e}_4, \mathbf{e}_5\}$. For $(-N^*, -M^*, -H)$ the top weight vectors (in ker $(-N^*)$) are $\mathbf{e}_2$ and $\mathbf{e}_5$, whose subscripts correspond to the bottom rows of the Jordan blocks; the weights (eigenvalues of $-H$) are again 1 and 2, as expected, since top weight vectors must be positive; and the chains are $\{\mathbf{e}_2, -\mathbf{e}_1\}$ and $\{\mathbf{e}_5, -\mathbf{e}_4, \mathbf{e}_3\}$. Notice that the top weight vectors can be read off from the positions of the highest eigenvalues of the strings on the diagonal of $-H$.

We are now ready to introduce the external transvectants between the invariants $\mathcal{I} = \ker \mathcal{X} \subset \mathbb{R}[[\mathbf{x}]]$ and the elements of $\ker(-N^*) = \ker N^* \subset \mathbb{R}^n$. This follows the pattern of (7.1), replacing $\acute{\mathcal{Y}}$ by $\mathcal{Y}$, $\grave{\mathcal{Y}}$ by $(-M^*)$, and $g$ by $\mathbf{w} \in \ker N^*$. Thus for integers $s$ satisfying

$$0 \le s \le \widehat{f} \quad \text{and} \quad 0 \le s \le \widehat{\mathbf{w}} \tag{17.2}$$

we have

$$\begin{aligned}(f, \mathbf{w})^{(s)} &= \sum_{j=0}^{s} (-1)^j W^{sj}_{f\mathbf{w}} (\mathcal{Y}^j f)(-M^{*(s-j)}\mathbf{w}) \\ &= (-1)^s \sum_{j=0}^{s} W^{sj}_{f\mathbf{w}} (\mathcal{Y}^j f)(M^{*(s-j)}\mathbf{w}),\end{aligned} \tag{17.3}$$

where

$$W^{sj}_{f\mathbf{w}} = \binom{s}{j} \cdot \frac{(\widehat{f}-j)!}{(\widehat{f}-s)!} \cdot \frac{(\widehat{\mathbf{w}}-s+j)!}{(\widehat{\mathbf{w}}-s)!}.$$

In particular,

$$(f, \mathbf{w})^{(1)} = -(\widehat{f})(f)(M^*\mathbf{w}) - (\widehat{\mathbf{w}})(\mathcal{Y} f)\mathbf{w}. \tag{17.4}$$

Note that because of the minus sign in front of $M^*$ in the first line of (17.3), the signs in the second line do not alternate; instead, each odd-strength transvectant has only minus signs, and each even-strength transvectant has only plus signs.

The remaining definitions are straightforward analogs of earlier definitions. An $\boldsymbol{\alpha}, \mathbf{B}$**-transvectant** is one of the form $(\boldsymbol{\alpha}^{\mathbf{k}}, \mathbf{b}_j)^{(s)}$, and is **admissible** if $\boldsymbol{\alpha}^{\mathbf{k}} \in A$. The weight of each $\mathbf{b}_j$ imposes an absolute limit to $s$, as opposed to the case of $(\boldsymbol{\alpha}^{\mathbf{k}}, \boldsymbol{\beta}^{\boldsymbol{\ell}})^{(s)}$, where the weight of $\boldsymbol{\beta}^{\boldsymbol{\ell}}$ can become arbitrarily large. The set of all admissible $\boldsymbol{\alpha}, \mathbf{B}$-transvectants is the **pre-box product** $A \mathbin{\widehat{\boxtimes}} \mathbf{B}$.

**17.4. Remark.** Since the pre-box product is used primarily as a basis, it is safe to change the signs of the elements, and we usually delete the negative signs of the odd-strength transvectants. We do not do this in the definition of the transvectants, because the signs may be significant in other uses, especially in identities. Similarly, in [29, Ch. 12], Sanders often omits numerical factors that occur in every term of a transvectant, when it is used as a basis element.

Since $\ker N^*$ is not a polynomial space, it has no Newton space, and the transvectants do not form a cone. Factorization makes sense only in the first argument $\boldsymbol{\alpha}^{\mathbf{k}}$. A **prime transvectant** is an $\boldsymbol{\alpha}\mathbf{B}$-transvectant which becomes ill-defined if any factor is removed from $\boldsymbol{\alpha}^{\mathbf{k}}$. Note that all transvectants are vectors; thus $\alpha_1, \dots, \alpha_p$ (which in the invariant case are transvectants of the form $(\alpha_i, 1)^{(0)}$) cannot be considered as primes in the current setting. The only primes of strength zero are $(1, \mathbf{e}_r)^{(0)} = \mathbf{e}_r$ for each bottom row $r$ of a Jordan block.

The following theorems correspond to the indicated theorems in previous sections, and the proofs, which are omitted, are similar. Replacements (in Theorem 17.6) are indicated by equivalence ($\equiv$) as before.

**17.5. Theorem** *(Compare Theorem 8.4). The set $A \,\widehat{\boxtimes}\, B$ of admissible $\boldsymbol{\alpha}, \mathbf{B}$-transvectants is a formal vector space basis for $\mathcal{E} = \ker \mathsf{X}$.*

**17.6. Theorem** *(Compare Theorems 9.1 and 9.2). Any element $(\boldsymbol{\alpha}^{\mathbf{k}}, \mathbf{e}_r)^{(s)}$ of $A \,\widehat{\boxtimes}\, B$ may be replaced by a product $\boldsymbol{\alpha}^{\mathbf{k}'}(\boldsymbol{\alpha}^{\mathbf{k}''}, \mathbf{e}_r)^{(s)}$ provided that $\mathbf{k}' + \mathbf{k}'' = \mathbf{k}$ and that $(\boldsymbol{\alpha}^{\mathbf{k}''}, \mathbf{e}_r)^{(s)}$ is well-defined; the resulting set will remain a formal vector space basis for $\mathcal{E}$. The transvectant remaining in any replacement is admissible.*

To define the **box product** $A \boxtimes \mathbf{B}$, a **factoring order** must be assigned; this is simply an ordering of the elements of $\boldsymbol{\alpha}$ (which may or may not coincide with the initially given ordering). Notice that the factoring order here is entirely separate from the list of primes; for invariants, these are the same. To compute the box product, we consider the general admissible transvectant and factor as many as possible of each $\alpha_i$ from the transvectant, proceeding in the chosen factoring order and stopping when the remaining transvectant is prime. We do not use p-matrices (but see Remark 17.9), but instead factor directly, as in Example 12.1. There are two reasons for this. First, only the individual $\alpha_i$ need to be removed, not complicated transvectants of the form $(\boldsymbol{\alpha}^{\boldsymbol{\kappa}}, \boldsymbol{\beta}^{\boldsymbol{\lambda}})^{(\sigma)}$; the difficulties of the algorithms A–D arise only because such primes exist in the case of invariants. Second, the absolute upper bound to $s$ coming from $\mathbf{e}_r$ keeps the number of cases finite. (In the expansion method of [26], boosting never required "recycling", which was the technique used there to handle arbitrarily large values of $s$.)

**17.7. Example.** $\mathcal{E}_2$, the normal form for $N_2$

Recall that for $N_2$, $\boldsymbol{\alpha}$ is the single invariant $\alpha_1 = x_1$; we will omit the subscript on $\alpha_1$. The preferred set $A = A_2$ given by (4.5) consists of all powers of $\alpha$. From the discussion above, $\mathbf{B} = (\mathbf{b}_1) = (\mathbf{e}_2)$. The pre-box product $A \,\widehat{\boxtimes}\, \mathbf{B}$ then consists of all $(\alpha^k, \mathbf{e}_2)^{(s)}$ with $s$ limited to 0 and 1 because $\widehat{\mathbf{e}}_2 = 1$. Since $\widehat{\alpha} = 1$, the factorizations of these transvectants are

$$(\alpha^k, \mathbf{e}_2)^{(0)} = \alpha^k \mathbf{e}_2$$
$$(\alpha^k, \mathbf{e}_2)^{(1)} \equiv \alpha^{k-1}(\alpha, \mathbf{e}_2)^{(1)}.$$

(The first of these is an equality, not just an equivalence.) The only prime is $(\alpha, \mathbf{e}_2)^{(1)}$. To illustrate the explicit computation of such transvectants, we have, from (17.4),

$$(\alpha, \mathbf{e}_2)^{(1)} = -(\widehat{\alpha})(\alpha)(M^* \mathbf{e}_2) - (\widehat{\mathbf{e}}_2)(\mathsf{Y}\alpha)\mathbf{e}_2 = -(x_1 \mathbf{e}_1 + x_2 \mathbf{e}_2),$$

since $\alpha = x_1$, $\widehat{\alpha} = \widehat{e_2} = 1$, $\mathcal{Y}\alpha = (x_2\partial/\partial x_1)x_1 = x_2$, and $M^*\mathbf{e}_2 = \mathbf{e}_1$. Therefore, changing the sign of the transvectant by Remark 17.4, the equivariants $\mathcal{E}_2 = \ker \mathsf{X}_2 \subset \mathcal{V}^2$ have the Stanley decomposition

$$\mathcal{E}_2 = \mathbb{R}[[\mathbf{x}_1]] \begin{bmatrix} 0 \\ 1 \end{bmatrix} \oplus \mathbb{R}[[\mathbf{x}_1]] \begin{bmatrix} x_1 \\ x_2 \end{bmatrix}.$$

Since for normal form purposes we should restrict to quadratic and higher terms, this becomes

$$\mathcal{E}_2 = \mathbb{R}[[\mathbf{x}_1]]x_1^2 \begin{bmatrix} 0 \\ 1 \end{bmatrix} \oplus \mathbb{R}[[\mathbf{x}_1]]x_1 \begin{bmatrix} x_1 \\ x_2 \end{bmatrix},$$

in agreement with (1.3). □

**17.8. Example.** $\mathcal{E}_{23}$, the normal form for $N_{23}$

The box product for $\mathcal{I}_{23}$ was worked out in Examples 10.2, 12.2, and 16.2; we rename the Hilbert basis $\boldsymbol{\vartheta} = (\alpha_1, \beta_1, \beta_2, \gamma_1, \gamma_2)$ as $\boldsymbol{\alpha} = (\alpha_1, \alpha_2, \alpha_3, \alpha_4, \alpha_5)$, and write the basis for $\ker N^*$ as $\mathbf{B} = (\mathbf{b}_1, \mathbf{b}_2) = (\mathbf{e}_2, \mathbf{e}_5)$. To facilitate comparison with the treatment by Sanders in [29, §12.6.7], we will use his notations for this problem, which are

$$\mathfrak{a} = \alpha_1, \quad \mathfrak{b} = \alpha_2, \quad \mathfrak{c} = \alpha_3, \quad \mathfrak{d} = \alpha_4, \quad \mathfrak{e} = \alpha_5,$$
$$\mathfrak{u} = \mathbf{e}_2, \quad \mathfrak{v} = \mathbf{e}_5.$$

The block decomposition for $\mathcal{I}_{23}$ (with degrees, weights, and labels) is then

$$A = \left[\begin{array}{c|cc|cc} 1 & 1 & 2 & 2 & 3 \\ 1 & 2 & 0 & 1 & 0 \\ \mathfrak{a} & \mathfrak{b} & \mathfrak{c} & \mathfrak{d} & \mathfrak{e} \\ \hline \infty & \infty & \infty & 1 & \infty \\ 0 & 0 & 0 & 0 & 0 \end{array}\right]. \tag{17.5}$$

The pre-box product $A \mathbin{\widehat{\boxtimes}} \mathbf{B}$ consists (in our notation) of all transvectants of the forms $(\boldsymbol{\alpha}^{\mathbf{k}}, \mathbf{e}_2)^{(s=0,1)}$ or $(\boldsymbol{\alpha}^{\mathbf{k}}, \mathbf{e}_5)^{(s=0,1,2)}$, since $\widehat{\mathbf{e}}_2 = 1$ and $\widehat{\mathbf{e}}_5 = 2$.

**17.9. Remark.** By slightly modifying the p-matrix notation, this can be written as follows (with both Sanders' notations and ours indicated):

$$A \mathbin{\widehat{\boxtimes}} B = \left(\begin{array}{ccccc|cc|cc} \mathfrak{a} & \mathfrak{b} & \mathfrak{c} & \mathfrak{d} & \mathfrak{e} & \mathfrak{u} & \mathfrak{v} & & \\ \alpha_1 & \alpha_2 & \alpha_3 & \alpha_4 & \alpha_5 & \mathbf{b}_1 & \mathbf{b}_2 & s & s \\ \hline \infty & \infty & \infty & 1 & \infty & \multirow{2}{*}{1} & \multirow{2}{*}{0} & k_1 + 2k_2 + k_4 & 1 \\ 0 & 0 & 0 & 0 & 0 & & & 0 & 0 \end{array}\right)$$
$$\sqcup \left(\begin{array}{ccccc|cc|cc} \alpha_1 & \alpha_2 & \alpha_3 & \alpha_4 & \alpha_5 & \mathbf{b}_1 & \mathbf{b}_2 & s & s \\ \hline \infty & \infty & \infty & 1 & \infty & \multirow{2}{*}{0} & \multirow{2}{*}{1} & k_1 + 2k_2 + k_4 & 2 \\ 0 & 0 & 0 & 0 & 0 & & & 0 & 0 \end{array}\right)$$

Since $\mathbf{b}_1$ and $\mathbf{b}_2$ cannot be multiplied, the $\mathbf{b}_1$ and $\mathbf{b}_2$ columns do not call for $\mathbf{b}_1^1\mathbf{b}_2^0$ and $\mathbf{b}_1^0\mathbf{b}_2^1$. Instead, 1 indicates the presence of a basis element and 0 its absence. It would be possible in principle to apply Algorithms A–D to these p-matrices, but this would be harder than doing the factoring directly.

Since $\widehat{\mathfrak{c}} = \widehat{\mathfrak{e}} = 0$, these can be factored at any time, so (by Remark 12.3) we suspend them and write $\mathcal{R} = \mathbb{R}[[\mathfrak{c}, \mathfrak{e}]]$. There are six possible factoring orders for $\mathfrak{a}$, $\mathfrak{b}$, and $\mathfrak{d}$, and these give quite different results for the final Stanley decomposition for $\mathcal{E}_{23}$. The following table shows the number of terms in the Stanley decomposition that results from each factoring order.

| | | | |
|---|---|---|---|
| $(\mathfrak{a}, \mathfrak{b}, \mathfrak{d})$ | 14 | $(\mathfrak{b}, \mathfrak{a}, \mathfrak{d})$ | 15 |
| $(\mathfrak{a}, \mathfrak{d}, \mathfrak{b})$ | 16 | $(\mathfrak{b}, \mathfrak{d}, \mathfrak{a})$ | 18 |
| $(\mathfrak{d}, \mathfrak{a}, \mathfrak{b})$ | 19 | $(\mathfrak{d}, \mathfrak{b}, \mathfrak{a})$ | 21 |

The columns show that for any ordering of $\mathfrak{a}$ and $\mathfrak{b}$, the best result is obtained if $\mathfrak{d}$ is placed last, the next best if $\mathfrak{d}$ is second, and the worst if $\mathfrak{d}$ is first. The rows show that for any position for $\mathfrak{d}$, the best result is obtained if $\mathfrak{a}$ is factored before $\mathfrak{b}$. We will carry out the factorization leading to the best result, 14 terms with order $(\mathfrak{a}, \mathfrak{b}, \mathfrak{d})$. First, consider transvectants of the form $(\mathfrak{a}^i\mathfrak{b}^j\mathfrak{d}^k, \mathfrak{u})^{(s)}$, with $i, j \in [0, \infty]$ but $k \in [0, 1]$. Since $\widehat{\mathfrak{u}} = 1$, $s$ is restricted to 0, 1. The subscript on $\equiv$ indicates which of $\mathfrak{a}$, $\mathfrak{b}$, $\mathfrak{d}$ is being factored out at each stage.

$$
\begin{aligned}
(\mathfrak{a}^i\mathfrak{b}^j\mathfrak{d}^k, \mathfrak{u})^{(0)} &= \mathfrak{a}^i\mathfrak{b}^j\mathfrak{d}^k\mathfrak{u} = \begin{cases} \mathfrak{a}^i\mathfrak{b}^j\mathfrak{u} & k=0 \\ \mathfrak{a}^i\mathfrak{b}^j\mathfrak{d}\mathfrak{u} & k=1, \end{cases} \\
(\mathfrak{a}^i\mathfrak{b}^j\mathfrak{d}^k, \mathfrak{u})^{(1)} &\equiv_{\mathfrak{a}} \begin{cases} \mathfrak{a}^i(\mathfrak{b}^j\mathfrak{d}^k, \mathfrak{u})^{(1)} & (j,k) \neq (0,0) \\ \mathfrak{a}^{i-1}(\mathfrak{a}, \mathfrak{u})^{(1)} & j=k=0 \end{cases} \\
&\equiv_{\mathfrak{b}} \begin{cases} \begin{cases} \mathfrak{a}^i\mathfrak{b}^j(\mathfrak{d}, \mathfrak{u})^{(1)} & k=1 \\ \mathfrak{a}^i\mathfrak{b}^{j-1}(\mathfrak{b}, \mathfrak{u})^{(1)} & k=0, j \geq 1 \end{cases} \\ \mathfrak{a}_{i-1}(\mathfrak{a}, \mathfrak{u})^{(1)} & j=k=0 \end{cases} \\
&\equiv_{\mathfrak{d}} \text{ no change since } \mathfrak{d} \text{ cannot be factored.}
\end{aligned} \tag{17.6}
$$

Taking these terms in the order that they are produced, the $\mathfrak{u}$ portion of the Stanley decomposition is

$$
\begin{aligned}
&\mathcal{R}[[\mathfrak{a}, \mathfrak{b}]]\mathfrak{u} \oplus \mathcal{R}[[\mathfrak{a}, \mathfrak{b}]]\mathfrak{d}\mathfrak{u} \oplus \mathcal{R}[[\mathfrak{a}, \mathfrak{b}]](\mathfrak{d}, \mathfrak{u})^{(1)} \\
&\quad \oplus \mathcal{R}[[\mathfrak{a}, \mathfrak{b}]](\mathfrak{b}, \mathfrak{u})^{(1)} \oplus \mathcal{R}[[\mathfrak{a}]](\mathfrak{a}, \mathfrak{u})^{(1)}.
\end{aligned} \tag{17.7}
$$

Next, consider transvectants of the form $(\mathfrak{a}^i\mathfrak{b}^j\mathbf{d}^k, \mathfrak{v})^{(s)}$ with $k \in [0, 1]$. Since $\widehat{\mathfrak{v}} = 2$, $s$ is restricted to 0, 1, 2. The cases $s = 0, 1$ go exactly like the $\mathfrak{u}$ portion and produce

$$
\begin{aligned}
&\mathcal{R}[[\mathfrak{a}, \mathfrak{b}]]\mathfrak{v} \oplus \mathcal{R}[[\mathfrak{a}, \mathfrak{b}]]\mathfrak{d}\mathfrak{v} \oplus \mathcal{R}[[\mathfrak{a}, \mathfrak{b}]](\mathfrak{d}, \mathfrak{v})^{(1)} \\
&\quad \oplus \mathcal{R}[[\mathfrak{a}, \mathfrak{b}]](\mathfrak{b}, \mathfrak{v})^{(1)} \oplus \mathcal{R}[[\mathfrak{a}]](\mathfrak{a}, \mathfrak{v})^{(1)}.
\end{aligned} \tag{17.8}
$$

For $s = 2$ we have

$$(\mathfrak{a}^i\mathfrak{b}^j\mathfrak{d}^k, \mathfrak{v})^{(2)} \equiv_{\mathfrak{a}} \begin{cases} \mathfrak{a}^i(\mathfrak{b}^j\mathfrak{d}^k, \mathfrak{v})^{(2)} & 2j+k \geq 2 \text{ i.e. } j \geq 1, k = 0, 1 \\ \mathfrak{a}^{i-1}(\mathfrak{a}\mathfrak{d}, \mathfrak{v})^{(2)} & 2j+k = 1 \text{ i.e. } j = 0, k = 1 \\ \mathfrak{a}^{i-2}(\mathfrak{a}^2, \mathfrak{v})^{(2)} & 2j+k = 0 \text{ i.e. } j = k = 0, i \geq 2 \end{cases}$$

$$\equiv_{\mathfrak{b}} \begin{cases} \mathfrak{a}^i\mathfrak{b}^{j-1}(\mathfrak{b}\mathfrak{d}^k, \mathfrak{v})^{(2)} & j \geq 1, k = 0, 1, \\ \text{no change} \\ \text{no change} \end{cases}$$

$$\equiv_{\mathfrak{d}} \begin{cases} \begin{cases} \mathfrak{a}^i\mathfrak{b}^{j-1}(\mathfrak{b}, \mathfrak{v})^{(2)} & j \geq 1, k = 0 \\ \mathfrak{a}^i\mathfrak{b}^{j-1}\mathfrak{d}(\mathfrak{b}, \mathfrak{v})^{(2)} & j \geq 1, k = 1 \end{cases} \\ \text{no change} \\ \text{no change} \end{cases} \tag{17.9}$$

This gives the following terms for the Stanley decomposition of the $\mathfrak{v}$ part with $s = 2$:

$$\mathcal{R}[[\mathfrak{a}, \mathfrak{b}]](\mathfrak{b}, \mathfrak{v})^{(2)} \oplus \mathcal{R}[[\mathfrak{a}, \mathfrak{b}]]\mathfrak{d}(\mathfrak{b}, \mathfrak{v})^{(2)}$$
$$\oplus \mathcal{R}[[\mathfrak{a}]](\mathfrak{a}\mathfrak{d}, \mathfrak{v})^{(2)} \oplus \mathcal{R}[[\mathfrak{a}]](\mathfrak{a}^2, \mathfrak{v})^{(2)}. \tag{17.10}$$

The full Stanley decomposition for $\mathcal{E}_{23}$ with factoring order $(\mathfrak{a}, \mathfrak{b}, \mathfrak{d})$ consists of the 14 terms in (17.7), (17.8), and (17.10). Observe that each of the terms has a different Stanley basis element. Here is the full list of 14 such elements:

$$\mathfrak{u}, \quad \mathfrak{d}\mathfrak{u}, \quad (\mathfrak{d}, \mathfrak{u})^{(1)}, \quad (\mathfrak{b}, \mathfrak{u})^{(1)}, \quad (\mathfrak{a}, \mathfrak{u})^{(1)},$$
$$\mathfrak{v}, \quad \mathfrak{d}\mathfrak{v}, \quad (\mathfrak{d}, \mathfrak{v})^{(1)}, \quad (\mathfrak{b}, \mathfrak{v})^{(1)}, \quad (\mathfrak{a}, \mathfrak{v})^{(1)},$$
$$(\mathfrak{b}, \mathfrak{v})^{(2)}, \quad \mathfrak{d}(\mathfrak{b}, \mathfrak{v}), \quad (\mathfrak{a}\mathfrak{d}, \mathfrak{v})^{(2)}, \quad (\mathfrak{a}^2, \mathfrak{v}). \tag{17.11}$$

It was claimed, above, that the Stanley decomposition with factoring order $(\mathfrak{a}, \mathfrak{b}, \mathfrak{d})$ is the best among the six possible factoring orders. We will now prove that claim, *without* doing the factorizations for the other orders. Namely, we will prove that all 14 Stanley basis elements listed in (17.11) must also appear in each of the other five Stanley decompositions, along with whatever "extra" Stanley basis elements are present. (The presence of extra basis elements implies that the coefficient rings must be smaller for at least some of the 14 "necessary" basis elements. These coefficient rings can only be found by carrying out the full factorization for each factoring order.)

First, notice that all of the elements in (17.11) are prime except $\mathfrak{d}\mathfrak{u}$, $\mathfrak{d}\mathfrak{v}$, and $\mathfrak{d}(\mathfrak{b}, \mathfrak{v})^{(2)}$. Each of these primes is a well-defined transvectant and therefore must belong to some term of each possible Stanley decomposition. But, being prime, it can only belong to a term of the Stanley decomposition if it is the Stanley basis element of that term. The three exceptional elements are not prime, and their necessity must be proved differently. The transvectant $(\mathfrak{d}\mathfrak{b}, \mathfrak{v})^{(2)}$ is well-defined, so it must belong to some term in each Stanley decomposition. It cannot be the Stanley basis element of its term, since the factoring method does not terminate until every transvectant is completely factored. The only way to factor this transvectant further is $(\mathfrak{d}\mathfrak{b}, \mathfrak{v})^{(2)} \equiv \mathfrak{d}(\mathfrak{b}, \mathfrak{v})^{(2)}$. Since $\mathfrak{d}$ can appear only to the first power, it cannot belong to the coefficient ring of the term, but must be retained as a factor in the Stanley basis element, which must be $\mathfrak{d}(\mathfrak{b},\mathfrak{v})^{(2)}$. The arguments for $\mathfrak{d}\mathfrak{u}$ and $\mathfrak{d}\mathfrak{v}$ are similar (but simpler).

Although the Stanley decomposition given by (17.7), (17.8), and (17.10) is the best that can be obtained *by the factoring method*, it is still possible to achieve some improvement by going outside of the framework of the factoring method. We now show that it is possible to reduce the number of Stanley terms from 14 to 12, using a trick introduced by Sanders in [29, §12.6.7].

**17.10. Remark.** Sanders obtained a Stanley decomposition for $\mathcal{E}_{23}$ containing 17 terms, using the expansion method, and reduced this to 12 terms by five applications of this "Sanders trick". We obtained one with 14 terms by the factoring method, and will reduce it to 12 terms by two applications of the Sanders trick. It is easy to see that Sanders' Stanley decomposition with 17 terms is not "finished" (according to the factoring method) because it contains transvectants that are not prime, that is, the factorization is incomplete. So this example confirms that the factoring method is better than the expansion method, while still allowing improvements in the case of boosting. The Stanley decomposition in [29, p. 310] contains an error; $\mathcal{R}[[\mathfrak{b}]](\mathfrak{b},\mathfrak{u})^{(1)}$ should be replaced by $\mathcal{R}[[\mathfrak{a},\mathfrak{b}]](\mathfrak{b},\mathfrak{u})^{(1)}$. (This was copied incorrectly from page 309, line 6.) With this correction, it agrees with (17.12) below.

We claim that the terms $\mathcal{R}[[\mathfrak{a},\mathfrak{b}]]\mathfrak{d}\mathfrak{u}$ and $\mathcal{R}[[\mathfrak{a},\mathfrak{b}]]\mathfrak{d}\mathfrak{v}$ can be removed from our Stanley decomposition for $\mathcal{E}_{23}$ provided that the coefficient rings of the term $\mathcal{R}[[\mathfrak{a}]](\mathfrak{a},\mathfrak{u})^{(1)}$ and $\mathcal{R}[[\mathfrak{a}]](\mathfrak{a},\mathfrak{v})^{(1)}$ are increased to $\mathcal{R}[[\mathfrak{a},\mathfrak{b}]]$. This is proved in Lemma 17.12 below. It leaves us with the 12-term Stanley decomposition

$$\begin{aligned}
&\mathcal{R}[[\mathfrak{a},\mathfrak{b}]]\mathfrak{u}\oplus\mathcal{R}[[\mathfrak{a},\mathfrak{b}]](\mathfrak{d},\mathfrak{u})^{(1)}\\
&\quad\oplus\mathcal{R}[[\mathfrak{a},\mathfrak{b}]](\mathfrak{b},\mathfrak{u})^{(1)}\oplus\mathcal{R}[[\mathfrak{a},\mathfrak{b}]](\mathfrak{a},\mathfrak{u})^{(1)}\\
&\quad\oplus\mathcal{R}[[\mathfrak{a},\mathfrak{b}]]\mathfrak{v}\oplus\mathcal{R}[[\mathfrak{a},\mathfrak{b}]](\mathfrak{d},\mathfrak{v})^{(1)}\\
&\quad\oplus\mathcal{R}[[\mathfrak{a},\mathfrak{b}]](\mathfrak{b},\mathfrak{v})^{(1)}\oplus\mathcal{R}[[\mathfrak{a},\mathfrak{b}]](\mathfrak{a},\mathfrak{v})^{(1)}\\
&\quad\oplus\mathcal{R}[[\mathfrak{a},\mathfrak{b}]](\mathfrak{b},\mathfrak{v})^{(2)}\oplus\mathcal{R}[[\mathfrak{a},\mathfrak{b}]]\mathfrak{d}(\mathfrak{b},\mathfrak{v})^{(2)}\\
&\quad\oplus\mathcal{R}[[\mathfrak{a}]](\mathfrak{a}\mathfrak{d},\mathfrak{v})^{(2)}\oplus\mathcal{R}[[\mathfrak{a}]](\mathfrak{a}^2,\mathfrak{v})^{(2)}. \qquad \square
\end{aligned} \tag{17.12}$$

**17.11. Remark.** This decomposition not only does not follow from the factoring method, but indeed is outside the framework of the replacement theorem. For example, the second line of (17.12) contains $\mathfrak{a}(\mathfrak{b},\mathfrak{u})^{(1)}$ and $\mathfrak{b}(\mathfrak{a},\mathfrak{u})^{(1)}$, both of which are replacements for $(\mathfrak{a}\mathfrak{b},\mathfrak{u})^{(1)}$. This would not be permitted if the replacement theorem were being used.

**17.12. Lemma.** *The following equality of sets holds, and remains true if* $\mathfrak{u}$ *is replaced by* $\mathfrak{v}$.

$$\begin{aligned}
&\mathcal{R}[[\mathfrak{a},\mathfrak{b}]]\mathfrak{d}\mathfrak{u}\oplus\mathcal{R}[[\mathfrak{a},\mathfrak{b}]](\mathfrak{b},\mathfrak{u})^{(1)}\oplus\mathcal{R}[[\mathfrak{a}]](\mathfrak{a},\mathfrak{u})^{(1)}\\
&\quad=\mathcal{R}[[\mathfrak{a},\mathfrak{b}]](\mathfrak{b},\mathfrak{u})^{(1)}\oplus\mathcal{R}[[\mathfrak{a},\mathfrak{b}]](\mathfrak{a},\mathfrak{u})^{(1)}.
\end{aligned} \tag{17.13}$$

**Proof.** Sanders observes in [29, Cor. 12.3.11] that

$$\mathfrak{d}\mathfrak{u}=2\mathfrak{b}(\mathfrak{a},\mathfrak{u})^{(1)}-\mathfrak{a}(\mathfrak{b},\mathfrak{u})^{(1)}. \tag{17.14}$$

As a result, any element $f(\mathfrak{a},\mathfrak{b})\mathfrak{d}\mathfrak{u} \in \mathcal{R}[[\mathfrak{a},\mathfrak{b}]]\mathfrak{d}\mathfrak{u}$ can be written as

$$f(\mathfrak{a},\mathfrak{b})\mathfrak{d}\mathfrak{u} = 2f(\mathfrak{a},\mathfrak{b})\mathfrak{b}(\mathfrak{a},\mathfrak{u})^{(1)} - f(\mathfrak{a},\mathfrak{b})\mathfrak{a}(\mathfrak{b},\mathfrak{u})^{(1)},$$

showing that the right-hand side of (17.13) is contained in the right-hand side. For the reverse inclusion, expand the right-hand side as $\mathcal{R}[[\mathfrak{a},\mathfrak{b}]](\mathfrak{b},\mathfrak{u})^{(1)} \oplus \mathcal{R}[[\mathfrak{a}]](\mathfrak{a},\mathfrak{u})^{(1)} \oplus \mathcal{R}[[\mathfrak{a},\mathfrak{b}]]\mathfrak{b}(\mathfrak{a},\mathfrak{u})^{(1)}$ and note that the first two terms are contained in the left-hand side; it remains to check the third term. But by another use of (17.14), $f(\mathfrak{a},\mathfrak{b})\mathfrak{b}(\mathfrak{a},\mathfrak{u})^{(1)} = (1/2)f(\mathfrak{a},\mathfrak{b})[\mathfrak{a}(\mathfrak{b},\mathfrak{u})^{(1)} + \mathfrak{d}\mathfrak{u}] \in \mathcal{R}[[\mathfrak{a},\mathfrak{b}]]\mathfrak{d}\mathfrak{u} \oplus \mathcal{R}[[\mathfrak{a},\mathfrak{b}]](\mathfrak{b},\mathfrak{u})^{(1)}$. The same argument works with $\mathfrak{v}$ in place of $\mathfrak{u}$, if equation (17.14) is replaced by $2\mathfrak{d}\mathfrak{v} = 2\mathfrak{b}(\mathfrak{a},\mathfrak{v})^{(1)} - \mathfrak{a}(\mathfrak{b},\mathfrak{v})^{(1)}$. □

**17.13. Example.** $\mathcal{E}_4$, the equivariants of $N_4$.

Using the Hilbert basis $\boldsymbol{\alpha}$ for $\mathcal{I}_4$ given in (4.10), and the preferred set $A$ from (4.11), with the weights, our notation, and that of Sanders [29, §12.6.3], we have

$$A_4 = \left[\begin{array}{cccc} 3 & 2 & 3 & 0 \\ \alpha_1 & \alpha_2 & \alpha_3 & \alpha_4 \\ \mathfrak{a} & \mathfrak{b} & \mathfrak{c} & \mathfrak{d} \\ \hline \infty & \infty & 1 & \infty \\ 0 & 0 & 0 & 0 \end{array}\right].$$

As noted above, $B = (\mathbf{e}_4) = (\mathfrak{u})$, with $\widehat{\mathfrak{u}} = 3$. Suppressing $\mathfrak{d}$ and setting $\mathcal{R} = \mathbb{R}[[\mathfrak{d}]]$, the pre-box product is

$$A \,\widehat{\boxtimes}\, B = \left\{(\mathfrak{a}^i\mathfrak{b}^j\mathfrak{c}^k,\mathfrak{u})^{(s)} : i,j \in [0,\infty], k \in [0,1], s \in [0, 3i+2j+3k\right\}.$$

From the experience about factoring order gained in Example 17.8, we try the factoring order $(\mathfrak{b},\mathfrak{a},\mathfrak{c})$, with $\mathfrak{b}$ before $\mathfrak{a}$ so that the weights increase (which was good in $N_{23}$), but keeping $\mathfrak{c}$ till last because $k$ is restricted to $[0,1]$ (as was the case for $\mathfrak{d}$ in $N_{23}$). The factoring proceeds as follows. For $s = 0$,

$$(\mathfrak{a}^i\mathfrak{b}^j\mathfrak{c}^k,\mathfrak{u})^{(0)} = \mathfrak{a}^i\mathfrak{b}^j\mathfrak{c}^k\mathfrak{u} = \begin{cases} \mathfrak{a}^i\mathfrak{b}^j\mathfrak{u} & k=0 \\ \mathfrak{a}^i\mathfrak{b}^j\mathfrak{c}\mathfrak{u} & k=1. \end{cases}$$

The cases $s = 1$ and $s = 2$ can be treated together; if $s = 1$ then the requirement for removing all $\mathfrak{b}$ is $3i + 3k \geq 1$ is, and if $s = 2$, the requirement is $3i + 3k \geq 2$, but both of these are equivalent to $(i,k) \neq (0,0)$, which may be written as $(k = 1) \veebar (i \geq 1, k = 0)$. Then

$$\begin{aligned}(\mathfrak{a}^i\mathfrak{b}^j\mathfrak{c}^k,\mathfrak{u})^{(s)} &\equiv_{\mathfrak{b}} \begin{cases} \mathfrak{b}^j(\mathfrak{a}^i\mathfrak{c}^k,\mathfrak{u})^{(s)} & (k=1)\veebar(i\geq 1, k=0) \\ \mathfrak{b}^{j-1}(\mathfrak{b},\mathfrak{u})^{(s)} & i=k=0, j\geq 1 \end{cases} \\ &\equiv_{\mathfrak{a}} \begin{cases} \begin{cases} \mathfrak{a}^i\mathfrak{b}^j(\mathfrak{c},\mathfrak{u})^{(2)} & k=1 \\ \mathfrak{a}^{i-1}\mathfrak{b}^j(\mathfrak{a},\mathfrak{u})^{(s)} & i\geq 1, k=0 \end{cases} \\ \mathfrak{b}^{j-1}(\mathfrak{b},\mathfrak{u})^{(s)} & i=k=0, j\geq 1 \end{cases} \\ &\equiv_{\mathfrak{c}} \text{no change.}\end{aligned}$$

If $s = 3$ the requirement for removing all $\mathfrak{b}$ is $3i + 3k \geq 3$, which is again equivalent to $(k = 1) \veebar (i \geq 1, k = 0)$. However, this time if $i = k = 0$ the original transvectant $(\mathfrak{a}^i \mathfrak{b}^j \mathfrak{c}^k, \mathfrak{u})^{(3)}$ is defined only if $j \geq 2$, and the factoring goes a little differently.

$$\begin{aligned}(\mathfrak{a}^i \mathfrak{b}^j \mathfrak{c}^k, \mathfrak{u})^{(3)} &\equiv_{\mathfrak{b}} \begin{cases} \mathfrak{b}^j (\mathfrak{a}^i \mathfrak{c}^k, \mathfrak{u})^{(3)} & (k=1) \veebar (i \geq 1, k=0) \\ \mathfrak{b}^{j-2}(\mathfrak{b}^2, \mathfrak{u})^{(3)} & i = k = 0, j \geq 2 \end{cases} \\ &\equiv_{\mathfrak{b}} \begin{cases} \begin{cases} \mathfrak{a}^i \mathfrak{b}^j (\mathfrak{c}, \mathfrak{u})^{(3)} & k = 1 \\ \mathfrak{a}^{i-1} \mathfrak{b}^j (\mathfrak{a}, \mathfrak{u})^{(3)} & k = 0, i \geq 1 \end{cases} \\ \mathfrak{b}^{j-2}(\mathfrak{b}^2, \mathfrak{u})^{(3)} & i = k = 0, j \geq 2 \end{cases} \\ &\equiv_{\mathfrak{c}} \text{no change.}\end{aligned}$$

Thus there exist eleven Stanley basis elements, with the following Stanley decomposition:

$$\begin{aligned}&\mathcal{R}[[\mathfrak{a}, \mathfrak{b}]]\mathfrak{u} \oplus \mathcal{R}[[\mathfrak{a}, \mathfrak{b}]]\mathfrak{c}\mathfrak{u} \\ &\quad \oplus \mathcal{R}[[\mathfrak{a}, \mathfrak{b}]](\mathfrak{c}, \mathfrak{u})^{(1)} \oplus \mathcal{R}[[\mathfrak{a}, \mathfrak{b}]](\mathfrak{a}, \mathfrak{u})^{(1)} \oplus \mathcal{R}[[\mathfrak{b}]](\mathfrak{b}, \mathfrak{u})^{(1)} \\ &\quad \oplus \mathcal{R}[[\mathfrak{a}, \mathfrak{b}]](\mathfrak{c}, \mathfrak{u})^{(2)} \oplus \mathcal{R}[[\mathfrak{a}, \mathfrak{b}]](\mathfrak{a}, \mathfrak{u})^{(2)} \oplus \mathcal{R}[[\mathfrak{b}]](\mathfrak{b}, \mathfrak{u})^{(2)} \\ &\quad \oplus \mathcal{R}[[\mathfrak{a}, \mathfrak{b}]](\mathfrak{c}, \mathfrak{u})^{(3)} \oplus \mathcal{R}[[\mathfrak{a}, \mathfrak{b}]](\mathfrak{a}, \mathfrak{u})^{(3)} \oplus \mathcal{R}[[\mathfrak{b}]](\mathfrak{b}^2, \mathfrak{u})^{(3)}.\end{aligned}$$

All of these Stanley basis elements are prime except for $\mathfrak{c}\mathfrak{u}$, which is clearly unavoidable, as in Example 17.8; therefore no shorter Stanley decomposition than this one can be obtained by the factoring method. However, using the Sanders trick it is possible to reduce the number of Stanley terms to nine. Specifically, we can:

1. Drop $\mathcal{R}[[\mathfrak{a}, \mathfrak{b}]]\mathfrak{c}\mathfrak{u}$ by enlarging $\mathcal{R}[[\mathfrak{b}]](\mathfrak{b}, \mathfrak{u})^{(1)}$ to $\mathcal{R}[[\mathfrak{a}, \mathfrak{b}]](\mathfrak{b}, \mathfrak{u})^{(1)}$ using the identity
$$3\mathfrak{c}\mathfrak{u} = 2\mathfrak{b}(\mathfrak{a}, \mathfrak{u})^{(1)} - 3\mathfrak{a}(\mathfrak{b}, \mathfrak{u})^{(1)}.$$
2. Drop $\mathcal{R}[[\mathfrak{a}, \mathfrak{b}]](\mathfrak{c}, \mathfrak{u})^{(1)}$ by enlarging $\mathcal{R}[[\mathfrak{b}]](\mathfrak{b}, \mathfrak{u})^{(2)}$ to $\mathbb{R}[[\mathfrak{a}, \mathfrak{b}]](\mathfrak{b}, \mathfrak{u})^{(2)}$ using the identity
$$4(\mathfrak{c}, \mathfrak{u})^{(1)} = 3\mathfrak{b}(\mathfrak{a}, \mathfrak{u})^{(2)} - 3\mathfrak{a}(\mathfrak{b}, \mathfrak{u})^{(2)}.$$

This produces

$$\begin{aligned}&\mathcal{R}[[\mathfrak{a}, \mathfrak{b}]]\mathfrak{u} \\ &\quad \oplus \mathcal{R}[[\mathfrak{a}, \mathfrak{b}]](\mathfrak{a}, \mathfrak{u})^{(1)} \oplus \mathcal{R}[[\mathfrak{a}, \mathfrak{b}]](\mathfrak{b}, \mathfrak{u})^{(1)} \\ &\quad \oplus \mathcal{R}[[\mathfrak{a}, \mathfrak{b}]](\mathfrak{c}, \mathfrak{u})^{(2)} \oplus \mathcal{R}[[\mathfrak{a}, \mathfrak{b}]](\mathfrak{a}, \mathfrak{u})^{(2)} \oplus \mathcal{R}[[\mathfrak{a}, \mathfrak{b}]](\mathfrak{b}, \mathfrak{u})^{(2)} \\ &\quad \oplus \mathcal{R}[[\mathfrak{a}, \mathfrak{b}]](\mathfrak{c}, \mathfrak{u})^{(3)} \oplus \mathcal{R}[[\mathfrak{a}, \mathfrak{b}]](\mathfrak{a}, \mathfrak{u})^{(3)} \oplus \mathcal{R}[[\mathfrak{b}]](\mathfrak{b}^2, \mathfrak{u})^{(3)}. \qquad \square \qquad (17.15)\end{aligned}$$

**17.14. Remark.** Sanders obtained a Stanley decomposition with 13 terms by the expansion method on page 300 of [29], reduced it to 11 terms on page 301 using two other identities, and reduced this further to 9 terms on page 302 using the identities we used above. His result agrees with ours.

## Acknowledgments

This paper has been under development for several years. I would like to thank Jan Sanders, Theodore Murdock, and David Malonza for discussions in person and by email during this time.